\documentclass[11pt]{amsart}
\usepackage{latexsym, graphics, amscd, amssymb}

\newcommand{\cA}{\mathcal{A}}

\newcommand{\cC}{\mathcal{C}}

\newcommand{\cG}{\mathcal{G}}
\newcommand{\cM}{\mathcal{M}}
\newcommand{\cN}{\mathcal{N}}

\newcommand{\cP}{\mathcal{P}}

\newcommand{\bC}{\mathbb{C}}

\newcommand{\bI}{\mathbb{I}}
\newcommand{\bQ}{\mathbb{Q}}
\newcommand{\bR}{\mathbb{R}}
\newcommand{\bZ}{\mathbb{Z}}

\newcommand{\RP}{\mathbb{RP}}

\newcommand{\ad}{\mathrm{ad}}
\newcommand{\Ad}{\mathrm{Ad}}
\newcommand{\Hom}{\mathrm{Hom}}

\newcommand{\Map}{\mathrm{Map}}
\newcommand{\Aut}{\mathrm{Aut}}

\newcommand{\rank}{\mathrm{rank}}
\newcommand{\Prin}{\mathrm{Prin}}
\newcommand{\diag}{\mathrm{diag}}
\newcommand{\Ker}{\mathrm{Ker}}

\newcommand{\fa}{\mathfrak{a}}
\newcommand{\fb}{\mathfrak{b}}
\newcommand{\fg}{\mathfrak{g}}

\newcommand{\ft}{\mathfrak{t}}
\newcommand{\fh}{\mathfrak{h}}
\newcommand{\fu}{\mathfrak{u}}
\newcommand{\fp}{\mathfrak{p}}

\newcommand{\fsl}{\mathfrak{sl}}
\newcommand{\fso}{\mathfrak{so}}
\newcommand{\fsp}{\mathfrak{sp}}
\newcommand{\fz}{\mathfrak{z}}
\newcommand{\fP}{\mathfrak{P}}
\newcommand{\fT}{\mathfrak{T}}

\newcommand{\tSi}{\tilde{\Si}}

\newcommand{\tP}{\tilde{P}}

\newcommand{\tV}{\tilde{V}}

\newcommand{\tg}{\tilde{g}}

\newcommand{\ta}{\tilde{a}}

\newcommand{\tb}{\tilde{b}}

\newcommand{\tal}{ {\tilde{\alpha}} }

\newcommand{\hLambda}{\widehat{\Lambda}}

\newcommand{\Si}{\Sigma}

\newcommand{\ep}{\epsilon}

\newcommand{\ab}{a_1,b_1,\ldots,a_\ell,b_\ell}

\newcommand{\pab}{\prod_{i=1}^\ell[a_i,b_i]}

\newcommand{\BC}{\overline{C}}

\newcommand{\ym}[2]{ X_{\mathrm{YM} }^{{#1},{#2}}(G_\bR) }

\newcommand{\Vym}[2]{ V_{\mathrm{YM} }^{{#1},{#2}}(G_\bR) }

\newcommand{\flU}[2]{ X_{\mathrm{flat}}^{ {#1}, {#2}}(U(n)) }

\newcommand{\ymU}[2]{ X_{\mathrm{YM} }^{{#1},{#2}}(U(n)) }

\newcommand{\ymSOO}[2]{ X_{\mathrm{YM} }^{{#1},{#2}}(SO(2n+1)) }

\newcommand{\ymSOE}[2]{ X_{\mathrm{YM} }^{{#1},{#2}}(SO(2n)) }

\newcommand{\VymSOO}[2]{ V_{\mathrm{YM} }^{{#1},{#2}}(SO(2n+1)) }

\newcommand{\VymSOE}[2]{ V_{\mathrm{YM} }^{{#1},{#2}}(SO(2n)) }

\newcommand{\VymSp}[2]{ V_{\mathrm{YM} }^{{#1},{#2}}(Sp(n)) }

\newcommand{\flS}[2]{ X_{\mathrm{flat}}^{ {#1}, {#2}}(U(1)) }

\newcommand{\dbar}{\bar{\partial}}

\newtheorem{thm}{Theorem}[section]

\newtheorem{fact}[thm]{Fact}
\newtheorem{lm}[thm]{Lemma}
\newtheorem{rem}[thm]{Remark}

\newtheorem{pro}[thm]{Proposition}
\newtheorem{ex}[thm]{Example}
\newtheorem{df}[thm]{Definition}
\newtheorem{no}[thm]{Notation}

\begin{document}

\title{Yang-Mills Connections on Orientable and Nonorientable Surfaces}
\author{Nan-Kuo Ho}
\address{ Department of Mathematics,
National Cheng-Kung University, Taiwan}
\email{nankuo@mail.ncku.edu.tw}
\thanks{The first author was supported by Grant NSC 95-2115-M-006-012-MY2 and
NSERC Postdoctoral Fellowship}

\author{Chiu-Chu Melissa Liu}
\address{Department of Mathematics, Northwestern University \\ and
Department of Mathematics, Columbia University}
\email{ccliu@math.northwestern.edu}

\subjclass[2000]{Primary 53D20; Secondary 58E15}

\keywords{}

\begin{abstract}
In \cite{ym},
Atiyah and Bott studied Yang-Mills functional over a Riemann surface from the point of view of Morse theory.
In \cite{HL4},
we study Yang-Mills functional on the space of
connections on a principal $G_\bR$-bundle over a closed, connected, nonorientable
surface, where $G_\bR$ is any compact connected Lie group.
In this paper, we generalize the discussion
in \cite{ym} and \cite{HL4}. We obtain explicit descriptions
of equivariant Morse stratification of Yang-Mills functional on orientable
and nonorientable surfaces for non-unitary classical groups $SO(n)$ and $Sp(n)$.
When the surface is orientable,
we use Laumon and Rapoport's method \cite{LR} to invert the
Atiyah-Bott recursion relation, and write down explicit formulas of rational
equivariant Poincar\'{e} series of the semistable stratum of the space of
holomorphic structures on a principal $SO(n,\bC)$-bundle or
a principal $Sp(n,\bC)$-bundle.
\end{abstract}

\maketitle
\tableofcontents

\section{Introduction}

Let $G_\bR$ be a compact, connected Lie group.
The complexification $G$ of $G_\bR$ is a connected reductive algebraic
group over $\bC$. For example, when $G_\bR=U(n)$, then $G=GL(n,\bC)$.
Let $P$ be a $C^\infty$ principal $G_\bR$-bundle over a Riemann surface
$\Si$, and let $\xi_0=P\times_{G_\bR} G$ be the associated $C^\infty$ principal
$G$-bundle. The space $\cA(P)$ of $G_\bR$-connections on $P$ is
isomorphic to the space $\cC(\xi_0)$ of $(0,1)$-connections
($\dbar$ operators) on $\xi_0$ as infinite dimensional complex
affine spaces. In the seminal paper \cite{ym}, Atiyah and Bott
obtained results on the topology of the moduli space $\cM(\xi_0)$ of
($S$-equivalence classes)
of semi-stable holomorphic structures on $\xi_0$ by studying
the Morse theory of the Yang-Mills functional on $\cA(P)$.
The absolute minimum of Yang-Mills functional is achieved by
{\em central} Yang-Mills connections, and $\cM(\xi_0)$ can
be identified with the moduli space of gauge equivalence
classes of central Yang-Mills connections on $P$.
When the absolute minimum of the Yang-Mills functional
is zero, which happens exactly when the obstruction class
$o(P)\in H^2(\Si,\pi_1(G))\cong \pi_1(G)$ is torsion,
the central Yang-Mills connections are flat connections, and $\cM(\xi_0)$ can
be identified with the moduli space of gauge equivalence
classes of flat connections on $P$.

Atiyah and Bott provided an algorithm
of computing the equivariant Poincar\'{e} series
$P_t^{\cG^\bC}(\cC_{ss};\bQ)$, where $\cC_{ss}$ is
the semi-stable stratum in $\cC(\xi_0)$ and
$\cG^\bC=\Aut(\xi_0)$ is the gauge group.
They proved that the stratification of $\cC(\xi_0)$ is
$\cG^\bC$-equivariantly perfect, so
$$
P_t^{\cG^\bC}(\cC(\xi_0);\bQ)
=P_t^{\cG^\bC}(\cC_{ss};\bQ) + \sum_{\lambda\in \Xi_{\xi_0}'}t^{2d_\mu}P_t^{\cG^\bC}(\cC_\mu;\bQ)
$$
where $d_\mu$ is the complex codimension of the stratum $\cC_\mu$, which
is a complex submanifold of $\cC(\xi_0)$, and
the sum is over all strata except for the top one $\cC_{ss}$.
The left hand side can be identified with $P_t(B\cG;\bQ)$, the rational Poincar\'{e}
series of the classifying space $B\cG$ of the gauge group $\cG=\Aut(P)$.
On the right hand side, $P_t^{\cG^\bC}(\cC_\mu;\bQ)$ can be related to the equivariant
Poincar\'{e} series of the top stratum of the space of connections on a principal
$G_\mu$-bundle, where $G_\mu$ is a subgroup of $G$.
So once $P_t(B\cG;\bQ)$ is computed, $P_t^{\cG^\bC}(\cC_{ss};\bQ)$ can
be computed recursively. Zagier solved
the recursion relation for $G=GL(n,\bC)$ in \cite{Za}, and Laumon and Rapoport
solved the recursion relation for a general connected reductive algebraic
group $G$ over $\bC$ in \cite{LR}.
The series $P_t^{\cG^\bC}(\cC_{ss};\bQ)$ can be identified
with $P_t^{G_\bR}(V_{ss}(P);\bQ)$, where $V_{ss}(P)$ is the representation
variety of central Yang-Mills connections on $P$. When the obstruction
class $o_2(P)\in H^2(\Si;\pi_1(G))\cong \pi_1(G)$ is torsion,
$V_{ss}(P)$ is the representation variety of flat Yang-Mills connections
on $P$, which is a connected component of $\Hom(\pi_1(\Si),G_\bR)$.

In \cite{HL4}, we study Yang-Mills functional on the space
of connections on a principal $G_\bR$-bundle $P$ over a closed,
connected, nonorientable surface $\Si$. By pulling back
connections to the orientable double cover $\pi:\tSi\to \Si$,
one gets an inclusion $\cA(P)\hookrightarrow\cA(\tP)$ from the space of connections
on $P$ to the space of connections on $\tilde{P}$,
where $\tP=\pi^*P\to \tSi$.
The Yang-Mills functional on $\cA(P)$ is the restriction of the
Yang-Mills functional
on $\cA(\tP)$. For nonorientable surfaces, the absolute
minimum of the Yang-Mills functional is zero for any $P$, achieved
by flat connections. The moduli space of gauge equivalence classes
of flat $G_\bR$-connections on $P$ can be identified with
a connected component of $\Hom(\pi_1(\Si), G_\bR)/G_\bR$, where
$G_\bR$ acts by conjugation.

In this paper, we generalize the discussion
in \cite{ym} and \cite{HL4} in the following directions:

\begin{enumerate}
\item In Section \ref{sec:gauge_group}, we compute the rational Poincar\'{e} series $P_t(B\cG;\bQ)$
of the classifying space of the gauge group $\cG$ of a principal
$G_\bR$-bundle over any closed connected (orientable or
nonorientable) surface.
The case where $\Si$ is orientable
is known (see \cite[Theorem 2.15]{ym}, \cite[Theorem 3.3]{LR}).

\item  When $\Si$ is orientable and $G_{ss}=[G,G]$ is not simply connected
(for example, when $G=G_{ss}=SO(n,\bC)$, $n>2$), the recursion relation
\cite[Theorem 3.2]{LR} that Laumon and Rapoport solved
in \cite{LR} is not exactly the Atiyah-Bott recursion relation
\cite[Theorem 10.10]{ym}. As a result,
their formula for $P_t^{ss}(G,\nu_G')$ \cite[Theorem 3.4]{LR}
is not exactly $P_t^{\cG^\bC}(\cC_{ss}(\xi_0);\bQ)$
when $G_{ss}$ is not simply connected.
In Appendix \ref{sec:LR},
we show that the method in \cite{LR} inverts
the Atiyah-Bott recursion relation and yields a closed
formula for $P_t^{\cG^\bC}(\cC_{ss}(\xi_0);\bQ)=
P_t^{G_\bR}(V_{ss}(P);\bQ)$, where $G_\bR$ is any
compact connected real Lie group (Theorem \ref{thm:Pt},
Theorem \ref{thm:LRformula}).

\item In \cite{HL4}, we established an exact correspondence
between the gauge equivalence classes of Yang-Mills
$G_\bR$-connections on $\Si$ and conjugacy classes of
representations $\Gamma_\bR(\Si)\to G_\bR$, where $\Gamma_\bR(\Si)$
is the super central extension of $\pi_1(\Si)$. This correspondence
allows us to obtain explicit description of $\cG$-equivariant Morse
stratification by studying the corresponding representation variety
of Yang-Mills connections. In Section \ref{sec:YMrepresentation}, we
recover the description in terms of Atiyah-Bott points for
orientable $\Si$, and determine candidates of Atiyah-Bott points for
nonorientable $\Si$.

\item In Section \ref{sec:SOodd}, Section \ref{sec:SOeven}, and Section \ref{sec:Sp}, we
give explicit descriptions of $\cG$-equivariant Morse strata of
Yang-Mills functional on orientable and non-orientable surfaces for
non-unitary classical groups $SO(2n+1)$, $SO(2n)$, and $Sp(n)$. When
$\Si$ is non-orientable, some twisted representation  varieties
(introduced and studied in Section \ref{sec:twisted} and Section
\ref{sec:twistedO}) arise in the reduction of these non-unitary
classical groups. This is new: in the $U(n)$ case (see \cite[Section
6, 7]{HL4}), the reduction involves only representation varieties of
$U(m)$, where $m<n$, of the nonorientable surface and of its double
cover.

\item When $\Si$ is orientable, we use the closed formula in (2) to write down explicit formulas
for $P_t^{G_\bR}(V_{ss}(P);\bQ)$ for non-unitary
classical groups (Theorem \ref{thm:PtSOodd}, Theorem \ref{thm:PtSOeven}, and Theorem \ref{thm:PtSp}).
These formulas are analogues of Zagier's formula for $U(n)$.
\end{enumerate}

The topology of $\Hom(\pi_1(\Si), G_\bR)/G_\bR$
is largely unknown when $\Si$ is nonorientable.
Using algebraic topology methods, T. Baird computed the $SU(2)$-equivariant
cohomology of $\Hom(\pi_1(\Si), SU(2))$ and the ordinary cohomology
of the quotient space $\Hom(\pi_1(\Si), SU(2))/SU(2)$ for any closed nonorientable surface $\Si$ \cite{B}.
He also proposed conjectures for general $G$.

For the purpose of Morse theory we should consider the Sobolev
space of $L_{k-1}^2$ connections  $\cA(P)^{k-1}$ and the group
of $L_k^2$ gauge transformations $\cG(P)^k$ and $\cG^\bC(P)^k$,
where $k\geq 2$. We will not emphasize the regularity issues
through out the paper, but refer the reader to \cite[Section 14]{ym}
and \cite{Da} for details.

\subsection*{Acknowledgments}
We thank Tom Baird, Ralph Cohen, Robert Friedman, Paul Goerss, Lisa Jeffrey, Eckhard Meinrenken,
John Morgan, David Nadler, Daniel Ramras,
Paul Selick, Michael Thaddeus, Jonathan Weitsman, and Christopher Woodward
for helpful conversations. We thank G\'{e}rard Laumon and Michael Rapoport for
confirming our understanding of their paper \cite{LR}.

\section{Topology of Gauge Group} \label{sec:gauge_group}

Let $\Si$ be a closed connected surface.
By classification of surfaces, $\Si$ is homeomorphic
to a Riemann surface of genus $\ell\geq 0$ if it
is orientable, and $\Si$ is homeomorphic to
the connected sum of $m >0$ copies of $\RP^2$ if
it is nonorientable.

Let $G_\bR$ be a compact connected Lie group.
Let $P$ be a principal $G_\bR$-bundle over $\Si$, and let
$\Aut(P)=\cG(P)$ be the gauge group. When $\Si$ is orientable,
the rational Poincar\'{e} series  $P_t(B\cG(P);\bQ)$ was computed in
\cite[Section 2]{ym} for $G_\bR=U(n)$. The computation can be
generalized to any general compact connected Lie group (see \cite[Theorem 3.3]{LR}).
In this section, we will compute $P_t(B\cG(P);\bQ)$ when
$G_\bR$ is any compact connected Lie group and $\Si$ is
any closed connected (orientable or non-orientable)
surface.

Following the strategy in \cite[Section 2]{ym}, we first find
the rational homotopy type of the classifying space $BG_\bR$ of $G_\bR$
(see \cite{Se}). Note that $BG_\bR$ is homotopic to $BG$, where
$G$ is the complexification of $G_\bR$.
Let $H_\bR$ be a maximal torus of $G_\bR$.
Then $H_\bR \cong U(1)^n$, and
$$
H^*(BH_\bR;\bZ)\cong \bZ[u_1,\ldots, u_n],
$$
where $u_i\in H^2(BH_\bR;\bZ)$. The Weyl group
$W$ acts on $H^*(BH_\bR;\bQ)\cong \bQ[u_1,\ldots,u_n]$, and
$$
H^*(BG_\bR;\bQ)\cong H^*(BH_\bR;\bQ)^W
\cong \bQ[I_1,\ldots, I_n]
$$
where $I_k$ is a homogeneous polynomial of degree $d_k$ in $u_1,\ldots,u_n$.
We may take $I_k\in \bZ[u_1,\ldots,u_n]$, so that
$I_k\in H^{2d_k}(BG_\bR;\bZ)$.
We may assume that $d_1=\cdots = d_r=1$, and
$d_k>1$ for $k>r$. Then $r=\dim_\bR(Z(G_\bR))$,
where $Z(G_\bR)$ is the center of $G_\bR$. In particular,
$r=0$ if and only if $G_\bR$ is semisimple.
The classes $I_1,\ldots, I_n$ are the universal characteristic classes
of principal $G_\bR$-bundles.  Each  $I_k\in H^{2d_k}(BG_\bR;\bZ)$
induces a continuous map
$I_k^*:BG_\bR \to K(\bZ;2d_k)$ to an Eilenberg-MacLane space,
so we have a continuous map
$$
BG_\bR \to \prod_{k=1}^n K(\bZ,2d_k).
$$
This is a rational homotopy equivalence.
\begin{fact}\label{thm:BG} Let $\stackrel{\bQ}{\simeq}$ denote rational homotopy equivalence. Then
$$
BG_\bR \stackrel{\bQ}{\simeq}\prod_{k=1}^n K(\bZ,2d_k)
$$
\end{fact}

In addition to Fact \ref{thm:BG}, we need the following two results:
\begin{pro}[{\cite[Proposition 2.4]{ym}}] \label{thm:BcG}
$$
B\cG(P) \simeq \Map_P(\Si,BG_\bR),
$$
where  the subscript $P$ denotes the component of a map of
$\Si$ into $BG_\bR$ which induces $P$.
\end{pro}

\begin{thm}[Thom]\label{thm:thom}
$$
\Map(X,K(A,n))=\prod_q K(H^q(X,A),n-q)
$$
where $K(A,n)$ is the Eilenberg-MacLane space characterized
by
$$
\pi_q(K(A,n))=\left\{\begin{array}{ll} 0 & q\neq n\\ A & q=n
\end{array}\right.
$$
\end{thm}

Since $\pi_q(X\times Y)=\pi_q(X)\times \pi_q(Y)$, we have
$$
K(A_1\times A_2,n)=K(A_1,n)\times K(A_2,n).
$$

Let $\Si$ be a Riemann surface of genus $\ell$. Then
\begin{eqnarray*}
&&\Map\Bigl(\Si,\prod_{k=1}^{n} K(\bZ,2d_k)\Bigr)
=  \prod_{k=1}^{n} \Map(\Si, K(\bZ,2d_k))\\
&=& \prod_{k=1}^{n} \Bigl( K(H^2(\Si;\bZ),2d_k-2)\times
K(H^1(\Si;\bZ),2d_k-1)\times K(H^0(\Si,\bZ),2d_k)\Bigr)\\
&=& \Bigl(\bZ\times K(\bZ,1)^{2\ell}\times K(\bZ,2)\Bigr)^r \\
&& \times \prod_{k=r+1}^{n}
\Bigl( K(\bZ,2d_k-2)\times K(\bZ,2d_k-1)^{2\ell}\times K(\bZ,2d_k) \Bigr)
\end{eqnarray*}
where the factor $\bZ^r$ corresponds to different connected components. So
\begin{eqnarray*}
\Map_P(\Si,BG_\bR)&\stackrel{\bQ}{\simeq}&
 \Bigl( K(\bZ,1)^{2\ell}\times K(\bZ,2)\Bigr)^r\\
&&\times \prod_{k=r+1}^{n} \Bigl( K(\bZ,2d_k-2)\times K(\bZ,2d_k-1)^{2\ell}\times K(\bZ,2d_k) \Bigr).
\end{eqnarray*}
It follows that
\begin{thm}[{\cite[Theorem 3.3]{LR}}]  \label{thm:BG_orientable}
Let $B\cG$ be the classifying space of  the gauge group $\cG$ of a principal $G_\bR$-bundle
over a Riemann surface of genus $\ell$. Then
$$
P_t(B\cG;\bQ)=\left(\frac{(1+t)^{2\ell}}{1-t^2}\right)^r\prod_{k=r+1}^{n}
\frac{(1+t^{2d_k-1})^{2\ell}}{(1-t^{2d_k-2})(1-t^{2d_k})}.
$$
\end{thm}
Note that $P_t(B\cG;\bQ)$ does not depend on the topological type
of the underlying principal $G_\bR$-bundle.

Let $\Si$ be the connected sum of $m>0$ copies of $\RP^2$. Then
\begin{eqnarray*}
&& \Map\Bigl(\Si,\prod_{k=1}^{n} K(\bZ,2d_k)\Bigr) = \prod_{k=1}^{n} \Map(\Si, K(\bZ,2d_k))\\
&=&  \prod_{k=1}^{n} \Bigl( K(H^2(\Si;\bZ),2d_k-2)\times K(H^1(\Si;\bZ),2d_k-1)\times
K(H^0(\Si,\bZ),2d_k)\Bigr)\\
&=& \prod_{k=1}^r\Bigl(\bZ/2\bZ\times K(\bZ,1)^{m-1}\times K(\bZ,2)\Bigr)\\
& & \times \prod_{k=r+1}^{n}
\Bigl( K(\bZ/2\bZ,2d_k-2)\times K(\bZ,2d_k-1)^{m-1}\times K(\bZ,2d_k) \Bigr)
\end{eqnarray*}
where the factor $(\bZ/2\bZ)^r$ corresponds to different connected components. So
$$
\Map_P(\Si,BG_\bR)\stackrel{\bQ}{\simeq}  \prod_{k=1}^{n}
\Bigl( K(\bZ,2d_k-1)^{m-1}\times K(\bZ,2d_k) \Bigr)
$$
It follows that
\begin{thm}\label{thm:BG_nonorientable}
Let $B\cG$ be the classifying space of the gauge group $\cG$ of a principal $G_\bR$-bundle
over a non-orientable surface which is diffeomorphic to
the connected sum of $m>0$ copies of $\RP^2$. Then
$$
P_t(B\cG;\bQ)=\prod_{k=1}^{n}
\frac{(1+t^{2d_k-1})^{m-1}}{(1-t^{2d_k})}.
$$
\end{thm}

For classical groups we have:
\begin{enumerate}
\item[(A)] $G_\bR=U(n)$: $W\cong S(n)$, the symmetric group, so
$$
H^*(BU(n);\bQ)=\bQ[u_1,\ldots,u_n]^{S(n)}=\bQ[c_1,\ldots,c_n],
$$
where $c_k$ is the $k$-th elementary symmetric function in $u_1,\ldots,u_n$.
In fact, the generator $c_k\in H^{2k}(BU(n);\bQ)$ is the universal rational $k$-th Chern class.
So $d_k=k$, $k=1,\ldots,n$.

\item[(B)] $G_\bR=SO(2n+1)$: $W=G(n)$, the wreath product of $\bZ/2\bZ$ by $S(n)$, so
$$
H^*(BSO(2n+1);\bQ)=\bQ[u_1,\ldots,u_n]^{G(n)}=\bQ[p_1,\ldots,p_n],
$$
where $p_k$ is the $k$-th elementary symmetric function in $u_1^2,\ldots,u_n^2$.
In fact, $p_k\in H^{4k}(BU(n);\bQ)$ is the universal rational $k$-th Pontrjagin class.
So $d_k=2k$, $k=1,\ldots,n$.

\item[(C)] $G_\bR=Sp(n)$: $W=G(n)$, the wreath product of $\bZ/2\bZ$ by $S(n)$, so
$$
H^*(BSp(n);\bQ)=\bQ[u_1,\ldots,u_n]^{G(n)}=\bQ[\sigma_1,\ldots,\sigma_n],
$$
where $\sigma_k$ is the $k$-th elementary symmetric function in $u_1^2,\ldots,u_n^2$.
So $d_k=2k$, $k=1,\ldots,n$.

\item[(D)] $G_\bR=SO(2n)$: $W=SG(n)$, the subgroup of $G(n)$ consisting of
even permutations, so
$$
H^*(BSO(2n);\bQ)=\bQ[u_1,\ldots,u_n]^{SG(n)}=\bQ[p_1,\ldots,p_{n-1},e],
$$
where $p_k$ is the $k$-th elementary symmetric function in $u_1^2,\ldots,u_n^2$, and
$e=u_1\cdots u_n$. In fact, $p_k\in H^{4k}(BU(n);\bQ)$ is the universal rational $k$-th Pontrjagin class,
and $e\in H^{2n}(BSO(2n);\bQ)$ is the universal rational Euler class.
So $d_k=2k$, $k=1,\ldots,n-1$, and $d_n=n$.
\end{enumerate}

\section{Holomorphic Principal Bundles over Riemann Surfaces}
\label{sec:algebraic_geometry}

Let $G$ be the complexification of a compact, connected real Lie group $G_\bR$.
Then  $G$ is a reductive algebraic group over $\bC$.
For example, if $G_\bR=U(n)$ then $G=GL(n,\bC)$.
We fix a topological principal $G_\bR$-bundle $P$ over a Riemann surface $\Si$,
and let $\xi_0= P\times_{G_\bR} G$ be the associated principal $G$-bundle. Then
the space $\cA(P)$ of $G_\bR$-connections on $P$ is isomorphic
to the space $\cC(\xi_0)$ of $(0,1)$-connections ($\dbar$-operators)
on $\xi_0$ as infinite dimensional complex affine spaces.
More explicitly, $\cA(P)$ and $\cC(\xi_0)$ are affine spaces
whose associated vector spaces are
$\Omega^1(\Si,\fg_\bR)$ and $\Omega^{0,1}(\Si,\fg)$, respectively, where
$\fg_\bR$ and $\fg=\fg_\bR \otimes_\bR \bC$ are the Lie algebras of $G_\bR$ and $G$, respectively.
Choose a local orthonormal frame $(\theta^1,\theta^2)$ of cotangent bundle $T^*_\Si$
of $\Si$ such that $*\theta^1=\theta^2$. Define an isomorphism $
j:\Omega^1(\Si,\fg_\bR)\to \Omega^{0,1}(\Si,\fg)$ by
$$
j(X_1\otimes \theta^1+ X_2\otimes \theta^2)= (X_1+ \sqrt{-1} X_2)\otimes (\theta^1-\sqrt{-1}\theta^2)
$$
where $X_1,X_2\in\Omega^0(\Si,\fg_\bR)$. It is easily checked that the definition is independent of
choice of $(\theta^1,\theta^2)$.

Harder and Narasimhan \cite{HN} defined a stratification on $\cC(\xi_0)$
when $G=GL(n,\bC)$, and Ramanathan \cite{Ra} extended this to general reductive groups.
It was conjectured by Atiyah and Bott in \cite{ym}, and proved by
Daskalopoulos in \cite{Da} (see also \cite{Rad}), that under the isomorphism $\cA(P)\cong \cC(\xi_0)$, the
stratification on $\cC(\xi_0)$ coincides with the Morse stratification
of the Yang-Mills functional on $\cA(P)$.

In this section, we first review the description of the
stratification in terms of Atiyah-Bott points, following
\cite[Section 10]{ym} and \cite{FM}. Then we write down the
Atiyah-Bott points for classical groups explicitly, similar to the
description of the stratification in terms of slopes when
$G_\bR=U(n)$.

\subsection{Preliminaries on reductive Lie groups and Lie algebras}
We have
$$
\fg= \fz_G \oplus [\fg,\fg]
$$
where $\fz_G$ is the center of $\fg$ and $[\fg,\fg]$ is the maximal semisimple
subalgebra of $\fg$.
Let $H_\bR$ be a maximal torus of $G_\bR$, and let $\fh_\bR$ be the Lie algebra of $H_\bR$.
Then $\fh=\fh_\bR\otimes_{\bR}\bC$ is a Cartan subalgebra of $\fg$.
Recall that any two maximal tori of $G_\bR$ are conjugate to each other,
and any two Cartan subalgebras of $\fg$ are conjugate to each other.
We have $\fh= \fz_G\oplus \fh'$
where $\fh'=\fh\cap[\fg,\fg]$.  Here we fix a choice of $H_\bR$, or equivalently,
we fix a Cartan subalgebra $\fh$ of $\fg$. Let $R$ be the root system associated
to $\fh$. We have
$$
\fg=\fh\oplus\bigoplus_{\alpha\in R}\fg_\alpha
=\fz_G \oplus \fh' \oplus \bigoplus_{\alpha\in R}\fg_\alpha.
$$

We choose a system of simple roots $\Delta\subset R$, and let $R_+$ be
the set of positive roots. The {\em Borel subalgebra} associated to $\Delta$
is given by
$$
\fb=\fh \oplus \bigoplus_{\alpha\in R_+} \fg_\alpha.
$$
The Lie algebra of a Borel subgroup $B$ of $G$ is a Borel subalgebra of $\fg$.
We have $B\cap G_\bR =H_\bR$.

A {\em parabolic subgroup} $P$ of $G$ is a subgroup containing
a Borel subgroup, and a {\em parabolic subalgebra} $\fp$ of $\fg$ is a subalgebra
containing a Borel subalgebra. A parabolic subalgebra containing $\fb$ is of the form
$$
\fp=\fh\oplus \bigoplus_{\alpha\in \Gamma}\fg_\alpha
$$
where
\begin{equation}\label{gamma}
\Gamma=R_+\cup\{\alpha\in R \mid \alpha\in \mathrm{span}(\Delta- I)\}.
\end{equation}
for some subset $I$ of the set $\Delta$ of simple roots. There is a one-to-one
correspondence between any two of the following:
\begin{enumerate}
\item[(i)] Subsets $I\subseteq \Delta$.
\item[(ii)] Parabolic subalgebras containing a fixed Borel subalgebra $\fb$.
\item[(iii)] Parabolic subgroups containing a fixed Borel subgroup $B$.
\end{enumerate}
In particular, $I$ being the empty set corresponds to $G$ (or $\fg$),
and $I$ being the entire set $\Delta$ corresponds
to $B$ (or $\fb$).

Given a parabolic subalgebra
$$
\fp=\fh\oplus \bigoplus_{\alpha\in \Gamma}\fg_\alpha,
$$
with $\Gamma$ as in \eqref{gamma}, define $-\Gamma$ to be the set of negatives of the members
of $\Gamma$. In other words,
$-\Gamma=-R_+\cup\{\alpha\in R\mid \alpha\in \mathrm{span}(\Delta-I)\}$.
let
$$
\mathfrak{l}=\fh\oplus \bigoplus_{\alpha\in \Gamma\cap -\Gamma} \fg_\alpha,\quad
\fu=\bigoplus_{\alpha\in \Gamma, \alpha\notin -\Gamma} \fg_\alpha
$$
so that $\fp=\mathfrak{l}\oplus \fu$. Then $\mathfrak{l}$, $\fu$ are subalgebras of $\fp$
and $\fu$ is an ideal of $\fp$. The subalgebra $\fu$ is nilpotent, and is called the
{\em nilpotent radical} of $\fp$. The subalgebra $\mathfrak{l}$ is reductive, and is called
the {\em Levi factor} of $\fp$.  Let $P$ be the parabolic subgroup
with Lie algebra $\fp$. Let $P=L U$ be the semi-direct product associated
to the direct sum $\fp=\mathfrak{l}\oplus \fu$, so that
the Lie algebras of $L$ and $U$ are $\mathfrak{l}$ and $\fu$, respectively.
The reductive Lie group $L$ is
called the {\em Levi factor} of $P$, and $U$ is called the {\em unipotent radical} of $P$.
We have $P\cap G_\bR =L_\bR$, the maximal compact subgroup of $L$; $L$ is
the complexification of $L_\bR$.

For simple Lie groups, there is a one-to-one correspondence between
simple roots and nodes of the Dynkin diagram.
In particular, a (proper) maximal parabolic subgroup corresponds to omitting one
node of the Dynkin diagram. See for example \cite[Lecture 23]{FH}.

\begin{enumerate}
\item[(A)] $G_\bR=SU(n)$, $G=SL(n,\bC)$, $n\geq 2$.\\
The Dynkin diagram of $\fsl(n,\bC)$ is $A_{n-1}$.  Omitting a node of
$A_{n-1}$, we get the disjoint union of $A_{n_1-1}$ and $A_{n_2-1}$,
where $n_1+n_2=n$, $n_1, n_2\geq 1$ (with the convention that $A_0$ is
empty).  The corresponding parabolic subgroup $P$ of $SL(n,\bC)$
is the subgroup which leaves the subspace
$\bC^{n_1}\times \{0\}$ of $\bC^n$ invariant. We have
$$
P\cap SU(n)=\{ \diag(A,B)\mid A\in U(n_1), B\in U(n_2), \det(A)\det(B)=1\}.
$$
For a general parabolic subgroup $P$ of $SL(n,\bC)$, we have
$$
P\cap SU(n)=\{ A\in U(n_1)\times \cdots\times U(n_r)\mid \det(A)=1 \}
$$
corresponding to omitting $(r-1)$ nodes, where $n_1+\cdots+n_r=n$, $n_i\geq 1$.

\item[(B)] $G_\bR=SO(2n+1)$, $G=SO(2n+1,\bC)$, $n\geq 1$.\\
The Dynkin Diagram of $\fso(2n+1,\bC)$ is $B_n$ (with the convention
$B_1=A_1$). Omitting a node of $B_n$, we get the
disjoint union of $A_{n_1-1}$ and $B_{n_2}$, where $n_1+n_2=n$,
$n_1\geq 1$, $n_2\geq 0$ (with the convention that $B_0$ is empty).
The corresponding parabolic subgroup of $SO(2n+1,\bC)$ is the subgroup
which leaves the following $n_1$-dimensional subspace of $\bC^{2n+1}$ invariant:
$$
\{ (z_1,\sqrt{-1}z_1,\ldots, z_{n_1}, \sqrt{-1}z_{n_1}, 0,\ldots,0)
\mid z_1,\ldots,z_{n_1}\in \bC \}.
$$
We have
$$
P\cap SO(2n+1)\cong U(n_1)\times SO(2n_2+1).
$$
For a general parabolic subgroup $P$ of $SO(2n+1,\bC)$, we have
$$
P\cap SO(2n+1)\cong U(n_1)\times \cdots \times U(n_{r-1})\times SO(2n_r+1)
$$ corresponding to omitting $(r-1)$ nodes,
where $n_1+\cdots +n_r=n$, $n_i\geq 1$ for $i\neq r$, and $n_r\geq 0$
(with the convention that $SO(1)$ is the trivial group).

\item[(C)]  $G_\bR=Sp(n)$, $G=Sp(n,\bC)$, $n\geq 1$.  \\
The Dynkin diagram of $\fsp(n,\bC)$ is $C_n$ (with the convention $C_1=A_1$).
Omitting a node from $C_n$, we get the disjoint union
of $A_{n_1-1}$ and $C_{n_2}$, where $n_1+n_2=n$, $n_1\geq 1$,
$n_2\geq 0$ (with the convention that $C_0$ is the empty set).
The corresponding parabolic subgroup of $Sp(n,\bC)$ is the subgroup which
leaves the subspace $\bC^{n_1}\times \{0\}$ of $\bC^{2n}$ invariant.
We have
$$
P\cap Sp(n)\cong U(n_1)\times Sp(n_2).
$$
For a general parabolic subgroup $P$ of $Sp(n,\bC)$, we have
$$
P\cap Sp(n)\cong U(n_1)\times \cdots \times U(n_{r-1})\times Sp(n_r)
$$ corresponding to omitting $(r-1)$ nodes,
where $n_1+\cdots +n_r=n$, $n_i\geq 1$ for $i\neq r$, and $n_r\geq 0$ (with
the convention that $Sp(0)$ is the trivial group).

\item[(D)] $G_\bR=SO(2n)$, $G=SO(2n,\bC)$, $n\geq 1$.\\
The Dynkin diagram of $\fso(2n,\bC)$ is $D_n$  (with
the convention $D_1=A_1$, $D_2= A_1\times A_1$, $D_3=A_3$).
Omitting a node of  $D_n$, we get the disjoint union
of $A_{n_1-1}$ and $D_{n_2}$, where $n_1+n_2=n$, $n_1\geq 1$, $n_2\geq 0$
(with the convention that $D_0$ is empty).
The corresponding parabolic subgroup of $SO(2n,\bC)$ is the subgroup
which leaves the following $n_1$-dimensional subspace of $\bC^{2n}$ invariant:
$$
\{ (z_1,\sqrt{-1}z_1,\ldots, z_{n_1}, \sqrt{-1}z_{n_1}, 0,\ldots,0)\mid z_1,\ldots,z_{n_1}\in \bC \}.
$$
We have
$$
P\cap SO(2n)\cong U(n_1)\times SO(2n_2).
$$
For a general parabolic subgroup $P$ of $SO(2n,\bC)$, we have
$$
P\cap SO(2n)\cong U(n_1)\times \cdots \times U(n_{r-1})\times SO(2n_r)
$$ corresponding to omitting $(r-1)$ nodes,
where $n_1+\cdots+n_r=n$, $n_i\geq 1$ for $i\neq r$, and $n_r\geq 0$ (with
the convention that $SO(0)$ is the trivial group). Note that
$SO(2)=U(1)$.
\end{enumerate}

\subsection{Harder-Narasimhan filtrations of dual vector bundles}

Let $E$ be a holomorphic vector bundle over $\Si$, and let
$$
0=E_0\subset E_1 \subset \cdots\subset  E_r=E
$$
be the Harder-Narasimhan filtration, where
$D_j=E_j/E_{j-1}$ is semi-stable, and the slopes
$\mu_j=\deg(D_j)/\rank(D_j)$ satisfy
$\mu_1>\cdots >\mu_r$. The vector $\mu=(\mu_1,\ldots,\mu_r)$
is the type of $E$. Let $\bI$ denote the trivial
holomorphic line bundle over $\Si$, and let
$E^\vee=\Hom(E,\bI)$ be the dual vector bundle, so that
$$
E^\vee_x =\Hom(E_x, \bC).
$$
Define the subbundle $E^\vee_{-j}$ of $E^\vee$ by
$$
(E^\vee_{-j})_x =\{ \alpha\in E^\vee_x \mid \alpha(v)=0\ \forall v\in (E_j)_x \}.
$$
then $(E^\vee_{-j})_x=(E_x/(E_j)_x)^\vee$, and we have
$$
0=E^\vee_{-r}\subset E^\vee_{-r+1}\subset \cdots \subset E^\vee_{-1}\subset E^\vee_0=E^\vee
$$
Let $F_j= E^\vee_{-r+j}/E^\vee_{-r+j-1}$. Then
$\rank F_j = \rank D_{r+1-j}$, $\deg F_j =-\deg D_{r+1-j}$, so
$\mu(F_j)=-\mu(D_{r+1-j})=-\mu_{r+1-j}$. The type of $E^\vee$ is given by
$(-\mu_r,\ldots,-\mu_1)$, where $-\mu_r>\cdots>-\mu_1$.

\subsection{Atiyah-Bott points}\label{sec:ABpts}

Let $\xi$ be a holomorphic principal $G$-bundle over a Riemann surface, and let
$E=\ad \xi =\xi \times_G \fg$ be the associated adjoint bundle.

The Lie algebra $\fg$ has a nondegenerate
invariant quadratic form $\fg \to \bC$. Therefore, there is a nondegenerate
invariant quadratic form $I$ on $E$, which implies $E$ is self-dual $E^\vee =E$. So the
Harder-Narasimhan filtration of $E$ is of the form
$$
0\subset E_{-r} \subset E_{-r+1}\subset \cdots \subset E_{-1}\subset E_0
\subset E_1 \subset \cdots \subset E_{r-1}\subset E.
$$
where
$$
(E_{-j})_x=\{ v\in E_x \mid I(u,v)=0\ \forall u\in (E_{j-1})_x \}
$$
and $D_0=E_0/E_{-1}$ has slope zero. Then $E_0$ is a parabolic subbundle
of the Lie algebra bundle $E$. The structure group $G$ of
$\xi$ can then be reduced to a parabolic subgroup $Q$, such that
$\xi=\xi_Q \times _Q G$, where $\xi_Q$ is a holomorphic principal $Q$-bundle with
$\ad\xi_Q = E_0$. The parabolic group is unique up to conjugation, and there is
a canonical choice for a fixed Borel subgroup $B$. This choice gives the
{\em Harder-Narasimhan reduction} and $Q$ is called
the {\em Harder-Narasimhan parabolic} of $\xi$.

The stratification of the space of holomorphic structures
on a fixed topological principal $G$-bundle $\xi$  is determined by
the Harder-Narasimhan parabolic $Q$ together with the topological type
of the underlying principal $Q$-bundle which is an element in $\pi_1(Q)$.
To make this more explicit, we describe the stratification
in terms of {\em Atiyah-Bott points}, following
\cite[Section 2]{FM}.

Let $H$ be a Cartan subgroup of $G$. Then $\pi_1(H)$ can be
viewed as a lattice in $\sqrt{-1}\fh_\bR$ such that $\pi_1(H)\otimes_\bZ \bR = \sqrt{-1}\fh_\bR$.
$$
\pi_1(H)\cong \{ X\in \sqrt{-1}\fh_\bR \mid \exp(2\pi\sqrt{-1}X)=e\} \subset \sqrt{-1}\fh_\bR.
$$
For example, $G_\bR=U(n)$, $\fh_\bR =\{2\pi\sqrt{-1}\diag(t_1,\ldots,t_n)\mid
t_1,\ldots,t_n \in\bR \}$, and $\pi_1(H)$ can be identified
with the lattice $\{ \diag(k_1,\ldots,k_n)\mid k_1,\ldots,k_n\in \bZ\}\subset \sqrt{-1}\fh_\bR$.

The set $\Delta^\vee$ of simple coroots span a sublattice $\Lambda$ of $\pi_1(H)$,
and $\pi_1(G)=\pi_1(H)/\Lambda$. The lattice $\Lambda$ is called the {\em coroot lattice} of $G$.
Let $\hLambda$ be the saturation of $\Lambda$ in
$\pi_1(H)$. Then $\pi_1(G_{ss})\cong \widehat{\Lambda}/\Lambda$. Under the above identification,
the short exact sequence of abelian groups
$$
1\to \pi_1(G_{ss})\to \pi_1(G)\to \pi_1(G/G_{ss}) \to 1
$$
can be rewritten as
$$
0 \to \hLambda/\Lambda \to \pi_1(H)/\Lambda \to \pi_1(H)/\hLambda \to 0,
$$
where $\hLambda/\Lambda$ is a finite abelian group, and $\pi_1(H)/\hLambda$
is a lattice.
Let $Z_0$ denote the connected component of the center of $G$ containing identity. Then
$D=Z_0\cap G_{ss}$ is a finite abelian group, and $G/G_{ss}\cong Z_0/D$.
$\pi_1(G/G_{ss})=\pi_1(H)/\hLambda$ can be identified with a lattice in
$\sqrt{-1} \fz_{G_\bR}$, where $\fz_{G_\bR}=\fz_G\cap \fh_\bR$,
such that $\pi_1(G/G_{ss})\otimes_\bZ \bR =\sqrt{-1}\fz_{G_\bR}$.

Let $\xi_0$ be a principal $G$-bundle over a Riemann surface $\Sigma$. Its topological type is classified by
the second obstruction class $c_1(\xi_0)\in H^2(\Sigma;\pi_1(G))\cong \pi_1(G)$. Let
$$
\mu(\xi_0)\in \pi_1(G/G_{ss})\subset \sqrt{-1}\fz_{G_\bR}
$$
be the image of $c_1(\xi_0)$ under the projection
$$
\pi_1(G)=\pi_1(H)/\Lambda \to \pi_1(G/G_{ss})=\pi_1(H)/\widehat{\Lambda}.
$$
The group $\hat{G}=\Hom(G,\bC^*)=\Hom(G/G_{ss},\bC^*)$
can be identified with the dual lattice of $\pi_1(H)/\widehat{\Lambda}$.

Let $P^I$ be a parabolic subgroup determined by $I\subseteq \Delta$, and let $L^I$ be its Levi factor.
The topological type of a principal $L^I$ bundle $\eta_0$ is determined by $c_1(\eta_0)\in \pi_1(L)$.
Given $\xi_0\in \Prin_G(\Sigma)$, we want to enumerate
\begin{equation}\label{eqn:set}
\{ \eta_0 \in \Prin_{L^I}(\Sigma) \mid \eta_0\times_{L^I} G= \xi_0\}.
\end{equation}

Consider the commutative diagram
$$
\begin{CD}
 0     @.  0\\
 @VVV @VVV\\
0 \ \longrightarrow\  \pi_1(L_{ss})=\hat{\Lambda}_L/\Lambda_L\  @>{j_{ss}}>> \pi_1(G_{ss})=\hLambda/\Lambda
@>{\oplus_{\alpha\in I}\varpi_\alpha}>> \oplus_{ \alpha\in I}\bQ/\bZ \\
  @V{i_L}VV @V{i_G}VV @|\\
  \pi_1(L)=\pi_1(H)/\Lambda_L @>{j}>> \pi_1(G)=\pi_1(H)/\Lambda
@>{\oplus_{\alpha\in I}\varpi_\alpha}>> \oplus_{\alpha\in I}\bQ/\bZ \\
 @V{p_L}VV @V{p_G}VV\\
 \pi_1(L/L_{ss})=\pi_1(H)/\hLambda_L @>{p}>> \pi_1(G/G_{ss})=\pi_1(H)/\hLambda \\
 @VVV @VVV\\
 0 @. 0 \\
\end{CD}
$$
where $\varpi_\alpha$ are the fundamental weights. In the above diagram,
the columns and the first row are exact.

Given a principal $L$-bundle $\eta_0$, $c_1(\eta_0)\in \pi_1(L)$ is determined by
$$
j(c_1(\eta_0))=c_1(\eta_0\times_L G) \in\pi_1(G),\quad p_L(c_1(\eta_0))=\mu(\eta_0)\in \pi_1(L/L_{ss}).
$$
Given $\xi_0\in \Prin_G(\Sigma)$, we have $c_1(\xi_0)\in \pi_1(G)$ and $\mu(\xi_0)\in \pi_1(G/G_{ss})$.
The map $p_L$ restricts to a bijection $j^{-1}(c_1(\xi_0)) \to p^{-1}(\mu(\xi_0))$. Note that
the set in \eqref{eqn:set} can be identified with $j^{-1}(c_1(\xi_0))$.

\begin{lm}[{\cite[Lemma 2.1.2]{FM}}]
Suppose that $\eta_0$ is a reduction of $\xi_0$ to a standard parabolic group $P^I$ for some
$I\subseteq \Delta$, possibly empty. The Atiyah-Bott point $\mu(\eta_0)$ and the topological type
of $\xi_0$ as a $G$-bundle determine the topological type of $\eta_0/U^I$ as an $L^I$-bundle
(and hence of $\eta_0$ as a $P^I$ bundle). Given a point
$\mu\in \fh_\bR$, there is a reduction of $\xi$ to a $P^I$-bundle whose Atiyah-Bott point is
$\mu$ if and only if the following conditions hold:
\begin{enumerate}
\item[(i)] $\mu\in \sqrt{-1}\fz_{{L^I}_\bR}$, where $\fz_{{L^I}_\bR}$
is the center of the Lie algebra of $L^I_\bR = L^I\cap G_\bR$.
\item[(ii)] For every simple root $\alpha\in I$ we have
$\varpi_\alpha(\mu)\equiv \varpi_\alpha(c)$ (mod $\bZ$).
\item[(iii)] $\chi(\mu)=\chi(c)$ for all characters $\chi$ of $G$.
\end{enumerate}
\end{lm}

\begin{df}[{\cite[Definition 2.1.3]{FM}}]
A pair $(\mu,I)$ consisting of a point $\mu\in\sqrt{-1}\fh_\bR$ and
a subset $I\subseteq \Delta$ is said to be of {\em Atiyah-Bott type} for
$c\in \pi_1(G)$ (or $\xi_0$ where $c_1(\xi_0)=c$) if (i)-(iii) hold.
A point $\mu\in\sqrt{-1}\fh_\bR$ is said to be of {\em Atiyah-Bott type} for $c$
if there is $I\subseteq \Delta$ such that $(\mu,I)$ is a pair of Atiyah-Bott type for $c$.
\end{df}

One may assume $\mu\in \overline{C}_0$, where $\overline{C}_0$ is the closure
of the fundamental Weyl chamber
$$
C_0=\{ X\in \sqrt{-1}\fh_\bR\mid \alpha(X)>0 \  \forall \alpha \in \Delta\}.
$$
We may choose the minimal $I$ such
that $\alpha(\mu)>0$ for all $\alpha\in I$. Then the stratum $\cC_\mu$ of the
space of $(0,1)$-connections on $\xi_0$ are indexed by points $\mu$ of Atiyah-Bott type of
$c_1(\xi_0)$ such that $\mu\in \overline{C}_0$. We may incorporate this by adding
\begin{enumerate}
\item[(iv)] $\alpha(\mu)>0$ for all $\alpha\in I$.
\end{enumerate}

Let $\cC(\xi_0)$ be the space
of all $(0,1)$-connections defining holomorphic structures
on a principal $G$-bundle $\xi_0$ with $c_1(\xi_0)=c\in \pi_1(G)$.
As a  summary of the above discussion, we have following
description of the Harder-Narasimhan stratification of $\cC$.
\begin{df}
Given a point $\mu\in\overline{C}_0$ of Atiyah-Bott type for $c$, the stratum $\cC_\mu \subset \cC(\xi_0)$
is the set of all $(0,1)$-connections defining holomorphic structures on $\xi_0$ whose
Harder-Narasimhan reduction has Atiyah-Bott type equal to $\mu$. The strata are preserved
by the action of gauge group. The union of these strata over all $\mu\in\overline{C}_0$ of
Atiyah-Bott type for $\xi_0$ is $\cC(\xi_0)$.
\end{df}

\subsection{Atiyah-Bott points for classical groups} \label{sec:ABC}
In this subsection, we assume
$$
n_1,\ldots, n_r\in \bZ_{>0}, \quad n_1+\cdots +n_r=n.
$$
\subsubsection{$G_\bR=U(n)$}\label{sec:UC}
$G=GL(n,\bC)$, and
$$
\sqrt{-1}\fh_\bR=\{  \diag(t_1,\ldots,t_n)\mid t_i \in \bR \}.
$$
Let $e_i\in \sqrt{-1}\fh_\bR$ be defined by $t_j=\delta_{ij}$. Then
$\{e_1,\ldots,e_n\}$ is a basis of $\sqrt{-1}\fh_\bR$. Let $\{\theta_1,\ldots,\theta_n\}$ be
the dual basis of $(\sqrt{-1}\fh_\bR)^\vee=\Hom_\bR(\sqrt{-1}\fh_\bR,\bR)$.  Then
\begin{eqnarray*}
&& \pi_1(H)=\bZ e_1\oplus \cdots \oplus \bZ e_n \subset \sqrt{-1}\fh_\bR\\
&& \Delta =\{ \alpha_i=\theta_i-\theta_{i+1}\mid i=1,\ldots,n-1\} \subset (\sqrt{-1}\fh_\bR)^\vee\\
&& \Delta^\vee =\{ \alpha_i^\vee =e_i-e_{i+1}\mid i=1,\ldots,n-1\}\subset \sqrt{-1}\fh_\bR
\end{eqnarray*}
$\pi_1(U(n))\cong \pi_1(GL(n,\bC))\cong \pi_1(H)/\Lambda \cong \bZ$
is generated by $e_1$ (mod $\Lambda$).
Let $c=ke_1$ (mod $\Lambda$). Then
$\mu$ satisfies (i)-(iv) in Section \ref{sec:ABpts} if and only if
$$
\mu=\diag\Bigl(\frac{k_1}{n_1}I_{n_1},\ldots,\frac{k_r}{n_r}I_{n_r}\Bigr)
$$
where
$$
k_1,\ldots,k_r\in\bZ,\quad k_1+\cdots+k_r=k,\quad
\frac{k_1}{n_1} >\frac{k_1}{n_2}>\cdots> \frac{k_r}{n_r}.
$$

\subsubsection{$G_\bR=SO(2n+1)$}\label{sec:SOoddC} $G=SO(2n+1,\bC)$, and
$$
\sqrt{-1}\fh_\bR=\{\sqrt{-1}\diag(t_1 J,\ldots,t_n J,0I_1)\mid t_i \in \bR \}
$$
where
$$
J=\left(\begin{array}{cc}0 & -1\\ 1 & 0 \end{array}\right)
$$
Let $e_i\in \sqrt{-1}\fh_\bR$ be defined by $t_j=\delta_{ij}$. Then
$\{e_1,\ldots,e_n\}$ is a basis of $\sqrt{-1}\fh_\bR$. Let $\{\theta_1,\ldots,\theta_n\}$ be
the dual basis of $(\sqrt{-1}\fh_\bR)^\vee$.  Then
\begin{eqnarray*}
&& \pi_1(H)=\bZ e_1\oplus \cdots \oplus \bZ e_n\subset \sqrt{-1}\fh_\bR\\
&& \Delta=\{ \alpha_i=\theta_i-\theta_{i+1}\mid i=1,\ldots,n-1\}\cup \{\alpha_n=\theta_n\} \subset
(\sqrt{-1}\fh_\bR)^\vee\\
&& \Delta^\vee =\{ \alpha_i^\vee =e_i-e_{i+1}\mid i=1,\ldots,n-1\}\cup\{\alpha_n^\vee =2e_n\} \subset
\sqrt{-1}\fh_\bR
\end{eqnarray*}
$\pi_1(SO(2n+1))\cong \pi_1(SO(2n+1,\bC)\cong \bZ/2\bZ$
is generated by $e_n$ (mod $\Lambda$).
$c=ke_n$ (mod $\Lambda$) corresponds to $w_2=k$ where $k=0,1$.

\paragraph{Case 1} $\alpha_n\in I$. Then $\mu$ satisfies (i)-(iv) in Section \ref{sec:ABpts} if
and only if
$$
\mu=\sqrt{-1}\diag\Bigl(\frac{k_1}{n_1}J_{n_1},\ldots,\frac{k_r}{n_r}J_{n_r},0 I_1\Bigr)
$$
where
$$
k_1,\ldots,k_r\in \bZ,\quad
k_1+\cdots +k_r \equiv k\ (\mathrm{mod}\ 2\bZ),\quad
\frac{k_1}{n_1}>\frac{k_2}{n_2}>\cdots \frac{k_r}{n_r}>0.
$$

\paragraph{Case 2} $\alpha_n \notin I$.
Then $\mu$ satisfies (i)-(iv) in Section \ref{sec:ABpts} if and only if
$$
\mu=\sqrt{-1}\diag\Bigl(\frac{k_1}{n_1}J_{n_1},\ldots,\frac{k_{r-1}}{n_{r-1}}J_{n_{r-1}},0 I_{2n_r+1}\Bigr)
$$
where
$$
k_1,\ldots, k_{r-1}\in\bZ,\quad
\frac{k_1}{n_1}>\frac{k_2}{n_2}>\cdots >\frac{k_{r-1}}{n_{r-1}}>0.
$$

\subsubsection{$G_\bR=SO(2n)$}\label{sec:SOevenC} $G=SO(2n,\bC)$, and
$$
\sqrt{-1}\fh_\bR=\{ \sqrt{-1} \diag(t_1 J,\ldots,t_n J)\mid t_i \in \bR \}
$$
where
$$
J=\left(\begin{array}{cc}0 & -1\\ 1 & 0 \end{array}\right)
$$
Let $e_i\in \sqrt{-1}\fh_\bR$ be defined by $t_j=\delta_{ij}$. Then
$\{e_1,\ldots,e_n\}$ is a basis of $\sqrt{-1}\fh_\bR$. Let $\{\theta_1,\ldots,\theta_n\}$ be
the dual basis of $(\sqrt{-1}\fh_\bR)^\vee$.  Then
\begin{eqnarray*}
&& \pi_1(H)=\bZ e_1\oplus \cdots \oplus \bZ e_n\subset \sqrt{-1}\fh_\bR\\
&& \Delta=\{ \alpha_i=\theta_i-\theta_{i+1}\mid i=1,\ldots,n-1\}\cup \{\alpha_n=\theta_{n-1}+\theta_n\}
\subset (\sqrt{-1}\fh_\bR)^\vee\\
&& \Delta^\vee =\{ \alpha_i^\vee =e_i-e_{i+1}\mid i=1,\ldots,n-1\}\cup\{\alpha_n^\vee =e_{n-1}+e_n\}
\subset \sqrt{-1}\fh_\bR
\end{eqnarray*}
$\pi_1(SO(2n))\cong \pi_1(SO(2n,\bC))\cong \bZ/2\bZ$ is generated by $e_n$ (mod $\Lambda$).
$c=ke_n$ (mod $\Lambda$) corresponds to $w_2=k$ where $k=0,1$.

\paragraph{Case 1} $\alpha_{n-1},\alpha_n\in I$, $n_r=1$.
Then $\mu$ satisfies (i)-(iv) in Section \ref{sec:ABpts} if and only if
$$
\mu=\sqrt{-1}\diag\Bigl(\frac{k_1}{n_1}J_{n_1},\ldots,\frac{k_{r-1}}{n_{r-1}}J_{n_{r-1}},k_r J\Bigr)
$$
where
$$
k_1,\ldots,k_r \in \bZ,\quad
k_1+\cdots +k_r \equiv k\ (\mathrm{mod}\ 2\bZ),\quad
\frac{k_1}{n_1}>\frac{k_2}{n_2}>\cdots>\frac{k_{r-1}}{n_{r-1}}>|k_r|.
$$

\paragraph{Case 2} $\alpha_{n-1}\in I$, $\alpha_n\notin I$, $n_r>1$.
Then $\mu$ satisfies (i)-(iv) in Section \ref{sec:ABpts} if and only if
$$
\mu=\sqrt{-1}\diag\Bigl(\frac{k_1}{n_1}J_{n_1},\ldots,
\frac{k_{r-1}}{n_{r-1}} J_{n_{r-1}},
\frac{k_r}{n_r} J_{n_r-1},-\frac{k_r}{n_r} J\Bigr)
$$
where
$$
k_1,\ldots, k_r\in\bZ,\quad
k_1+\cdots+k_r=k\ (\mathrm{mod}\ 2\bZ),\quad
\frac{k_1}{n_1}>\frac{k_2}{n_2}>\cdots>\frac{k_r}{n_r}>0.
$$
\paragraph{Case 3} $\alpha_{n-1}\notin I$, $\alpha_n \in I$, $n_r>1$.
Then $\mu$ satisfies (i)-(iv) in Section \ref{sec:ABpts} if and only if
$$
\mu=\sqrt{-1}\diag\Bigl(\frac{k_1}{n_1}J_{n_1},\ldots,\frac{k_r}{n_r}J_{n_r}\Bigr)
$$
where
$$
k_1,\ldots, k_r\in\bZ,\quad
k_1+\cdots+k_r=k\ (\mathrm{mod}\ 2\bZ),\quad
\frac{k_1}{n_1}>\frac{k_2}{n_2}>\cdots>\frac{k_r}{n_r}>0.
$$

\paragraph{Case 4} $\alpha_{n-1}\notin I$, $\alpha_n\notin I$.
Then $\mu$ satisfies (i)-(iv) in Section \ref{sec:ABpts}  if and only if
$$
\mu=\sqrt{-1}\diag\Bigl(\frac{k_1}{n_1}J_{n_1},\ldots,\frac{k_{r-1}}{n_{r-1}} J_{n_r}, 0 J_{n_r}\Bigr)
$$
where
$$
k_1,\ldots, k_{r-1}\in\bZ,\quad
\frac{k_1}{n_1}>\frac{k_2}{n_2}>\cdots>\frac{k_{r-1}}{n_{r-1}}>0.
$$

\subsubsection{$G_\bR=Sp(n)$}\label{sec:SpC} $G=Sp(n,\bC)$, and
$$
\sqrt{-1}\fh_\bR=\{  \diag(t_1,\ldots,t_n,-t_1,\ldots,-t_n)\mid t_i \in \bR \}
$$
Let $e_i\in \sqrt{-1}\fh_\bR$ be defined by $t_j=\delta_{ij}$. Then
$\{e_1,\ldots,e_n\}$ is a basis of $\sqrt{-1}\fh_\bR$. Let $\{\theta_1,\ldots,\theta_n\}$ be
the dual basis of $(\sqrt{-1}\fh_\bR)^\vee$.  Then
\begin{eqnarray*}
&& \pi_1(H)=\bZ e_1\oplus \cdots \oplus \bZ e_n\subset \sqrt{-1}\fh_\bR\\
&& \Delta=\{ \alpha_i=\theta_i-\theta_{i+1}\mid i=1,\ldots,n-1\}\cup\{2\theta_n\}
\subset (\sqrt{-1}\fh_\bR)^\vee\\
&& \Delta^\vee =\{ \alpha_i^\vee =e_i-e_{i+1}\mid i=1,\ldots,n-1\} \cup\{ e_n\}\subset
\sqrt{-1}\fh_\bR
\end{eqnarray*}
$\pi_1(Sp(n))\cong \pi_1(Sp(n,\bC))$ is trivial.
\paragraph{Case 1.}$\alpha_n\in I$.
Then $\mu$ satisfies (i)-(iv) in Section \ref{sec:ABpts} if and only of
$$
\mu=\diag\Bigl(\frac{k_1}{n_1}I_{n_1},\ldots,\frac{k_r}{n_r}I_{n_r},
-\frac{k_1}{n_1}I_{n_1},\ldots,-\frac{k_r}{n_r}I_{n_r} \Bigr)
$$
where
$$
k_1,\ldots,k_r\in\bZ,\quad
\frac{k_1}{n_1}>\frac{k_2}{n_2}>\cdots >\frac{k_r}{n_r}>0.
$$

\paragraph{Case 2.}$\alpha_n\notin I$.
Then $\mu$ satisfies (i)-(iv) of Section \ref{sec:ABpts} if and only if
$$
\mu=\diag\Bigl(\frac{k_1}{n_1}I_{n_1},\ldots,\frac{k_{r-1}}{n_{r-1}}I_{n_{r-1}},0 I_{n_r}
-\frac{k_1}{n_1}I_{n_1},\ldots,-\frac{k_{r-1}}{n_{r-1}}I_{n_{r-1}}, 0 I_{n_r} \Bigr)
$$
where
$$
k_1,\ldots,k_{r-1}\in\bZ,\quad
\frac{k_1}{n_1}>\frac{k_2}{n_2}>\cdots >\frac{k_{r-1}}{n_{r-1}}>0.
$$

\section{Yang-Mills Connections and Representation Varieties}
\label{sec:YMrepresentation}

Let $G_\bR$ be a compact connected Lie group, and let
$P$ be a $C^\infty$ principal $G_\bR$-bundle
over a closed (orientable or nonorientable) surface.
In \cite[Section 3]{HL4}, we introduced Yang-Mills
functional and Yang-Mills connections on closed
nonorientable surfaces.

In this section, we study the connected components of the
representation variety of Yang-Mills connections. We recover the
description of the Morse stratification in terms of Atiyah-Bott
points for orientable $\Si$ (Section
\ref{sec:orientable-components}), and determine candidates of
Atiyah-Bott points for nonorientable $\Si$ (Section
\ref{sec:connected}). We also discuss and give a closed formula for
$G_\bR$-equivariant rational Poincar\'{e} series of the
representation variety of central Yang-Mills connections (Section
\ref{sec:poincare}). In Section \ref{sec:twisted} and Section
\ref{sec:twistedO}, we introduce certain twisted representation
varieties that will arise in Section \ref{sec:SOodd}, Section
\ref{sec:SOeven}, and Section \ref{sec:Sp}, and study their
connectedness.

\subsection{Representation varieties for Yang-Mills connections}
Let $\cA(P)$ be the space of $G_\bR$-connections
on $P$, and let $\cN(P)$ be the space of Yang-Mills connections
on $P$.
Let $\cG(P)=\Aut(P)$ be the gauge group, and let $\cG_0(P)$ be the
base gauge group. Let $\Gamma_\bR(\Si)$ be the super central
extension of $\pi_1(\Si)$ defined in \cite[Section 4.1]{HL4}.
\begin{thm}[{\cite[Theorem 6.7]{ym}, \cite[Theorem 4.6]{HL4}}]
There is a bijective correspondence between
conjugacy classes of homomorphisms $\Gamma_\bR(\Si)\to G_\bR$
and gauge equivalence classes of Yang-Mills $G_\bR$-connections
over $\Si$. In other words,
\begin{eqnarray*}
\bigcup_{P\in\Prin_{G_\bR}(\Si)} \cN(P)/\cG_0(P)
&\cong&\Hom(\Gamma_\bR(\Si),G_\bR)\\
\bigcup_{P\in\Prin_{G_\bR}(\Si)} \cN(P)/\cG(P)
&\cong&\Hom(\Gamma_\bR(\Si),G_\bR)/G_\bR
\end{eqnarray*}
\end{thm}

To describe $\Hom(\Gamma_\bR(\Si),G_\bR)$ more explicitly,
we introduce some notation.
Let $\Si^\ell_0$ be the closed, compact, connected, orientable
surface with $\ell\geq 0$ handles. Let $\Si^\ell_1$ be the connected
sum of $\Si^\ell_0$ and $\RP^2$, and let $\Si^\ell_2$ be the
connected sum of $\Si^\ell_0$ and a Klein bottle. Any closed,
compact, connected surface is of the form $\Si^\ell_i$, where
$\ell$ is a nonnegative integer and $i=0,1,2$. $\Si^\ell_i$ is
orientable if and only if $i=0$. Let $(G_\bR)_X$ denote
the stabilizer of $X$ of the adjoint action of $G_\bR$ on $\fg_\bR$.
With the above notation, $\Hom(\Gamma_\bR(\Si^\ell_i),G_\bR)$ can be
identified with the representation variety $\ym{\ell}{i}$, where

\begin{eqnarray*}
\ym{\ell}{0}&=&\bigl\{(\ab,X)\in {G_\bR}^{2\ell}\times \fg_\bR\mid \\
 && \quad a_i,b_i\in (G_\bR)_X,\ \pab=\exp(X)\bigr\}\\
\ym{\ell}{1}&=&\bigl\{(\ab,c,X)\in {G_\bR}^{2\ell+1}\times \fg_\bR\mid \\
&& \quad a_i,b_i\in (G_\bR)_X,\ \Ad(c)X=-X,\ \pab=\exp(X)c^2 \bigr\}\\
\ym{\ell}{2}&=&\bigl\{(\ab,d,c,X)\in {G_\bR}^{2\ell+2}\times \fg_\bR\mid \\
&& \quad a_i,b_i,d\in (G_\bR)_X,\ \Ad(c)X=-X,\ \pab=\exp(X)cdc^{-1}d \bigr\}
\end{eqnarray*}
The $G_\bR$-action on $\ym{\ell}{i}$ is given by
$$
g\cdot (c_1,\dots,c_{2\ell+i},X) = (g c_1 g^{-1},\dots, g c_{2\ell+i} g^{-1}, \Ad(g)X).
$$

\subsection[Connectedness of the representations for orientable surfaces]
{Connected components of the representation variety for orientable surfaces}\label{sec:orientable-components}

$G_\bR$ is connected, so the natural projection
$$
\ym{\ell}{0}
\to \ym{\ell}{0}/G_\bR
$$
induces a bijection
$$
\pi_0(\ym{\ell}{0})
\to \pi_0(\ym{\ell}{0}/G_\bR).
$$

Any point in $\ym{\ell}{0}/G_\bR$ can
be represented by
$$
(\ab,X)
$$
where $X\in \fh_\bR$. Such representative is unique if we
require that $\sqrt{-1}X$ is in the closure $\BC_0$ of the fundamental
Weyl chamber
$$
C_0= \{ Y \in \sqrt{-1}\fh_\bR \mid \alpha(Y) >0,\ \forall \alpha\in R_+ \}
=\{Y\in \sqrt{-1}\fh_\bR\mid \alpha(Y)>0,\ \forall \alpha\in\Delta\}.
$$

Given $X$ such that $\sqrt{-1}X\in \BC_0$, we want to find the
stabilizer $(G_\bR)_X$ of the adjoint action of $G_\bR$ on
$\fg_\bR$. Let $G$ be the complexification of $G_\bR$. We use the
notation in Section \ref{sec:algebraic_geometry}. Let
$$
I_X=\{ \alpha \in\Delta \mid \alpha(\sqrt{-1}X) > 0 \}.
$$
Then $I_X=\Delta$ if $\sqrt{-1}X\in C_0$, and $I_X$ is empty if and only if $X$ is
in the center $\fz_{G_\bR}$ of $\fg_\bR$.
Let
$$
\Gamma_X=R_+ \cup \{ \alpha\in R\mid \alpha\in \mathrm{span}(\Delta-I_X)\}.
$$
The stabilizer $\fg_X$ of the adjoint action of $\fg$ on itself is the Levi factor of the
parabolic subalgebra
$$
\fp_X= \fh \oplus\bigoplus_{\alpha\in \Gamma_X} \fg_\alpha.
$$
We have $\fp_X = \fg_X \oplus \fu_X$, where
$\fg_X$ and $\fu_X$ are the Levi factor and the nilpotent radical
of $\fp_X$, respectively. The Lie algebra of $G_X$ is $\fg_X$.
We conclude that
$$
(G_\bR)_X = L^{I_X}\cap G_\bR = L^{I_X}_\bR.
$$
Note that
$$
X\in \fz_{L^{I_X}_\bR}, \quad
\exp(X)=\pab \in (L^{I_X}_\bR)_{ss}.
$$
Let $\mu_X=\frac{\sqrt{-1}}{2\pi}X$. Then
$$
\mu_X \in \pi_1(H)/\widehat{\Lambda}_{L^{I_X}} \subset \sqrt{-1}\fz_{L^{I_X}_\bR}
\subset \sqrt{-1}\fh_\bR
$$
and $(\mu_X,I_X)$ is of Atiyah-Bott type
for some $c\in \pi_1(G)=\pi_1(G_\bR)$.

We now state the condition for $X \in \fh_\bR$ such that
$(\ab,X)\in \ym{\ell}{0}$ for some $(\ab)\in G^{2\ell}$.
Given $I\subseteq \Delta$, let $Z^I$ be the connected component
of the identity of the center of $L^I_\bR$, and let $D^I$ be the center of
$(L^I_\bR)_{ss}$. Then the Lie algebra for $Z^I$ is $\fz_{L^I_\bR}$. Denote
\begin{eqnarray*}
\Xi^I&=&\{ \mu \in\sqrt{-1} \fz_{L^I_\bR}\mid
\exp(-2\pi\sqrt{-1}\mu)\in D^I\} \cong \pi_1(Z^I/D^I)\cong \pi_1(L^I_\bR/(L^I_\bR)_{ss})\\
\Xi^I_+ &=&\{ \mu \in \Xi^I \cap\BC_0 \mid \alpha(\mu)>0
\textup{ iff }\alpha\in I\}.
\end{eqnarray*}

Given $\mu\in \Xi^I_+$, let $X_\mu=-2\pi\sqrt{-1}\mu \in \fh_\bR$.
Suppose that $(\ab,X)\in \ym{\ell}{0}$. Then there is a unique
pair $(\mu,I)$, where $I\subseteq\Delta$ and $\mu\in\Xi^I_+$, such
that $X$ is conjugate to $X_\mu$. Let $C_\mu\subset \fg_\bR$ denote the
conjugacy class of $X_\mu$, and define
$$
\ym{\ell}{0}_\mu=\{(\ab,X)\in  {G_\bR}^{2\ell}\times C_\mu\mid a_i,~b_i\in (G_\bR)_X, \pab=\exp(X)\}.
$$
Then $\ym{\ell}{0}$ is a disjoint union of
$$
\{\ym{\ell}{0}_\mu \mid \mu\in \Xi^I_+, I \subseteq \Delta \}.
$$
Each $\ym{\ell}{0}_\mu$ is a union of finitely many connected
components of $\ym{\ell}{0}$.

Note that $(G_\bR)_{X_\mu} = L^I_\bR$ for $\mu\in \Xi^I_+$. We
define {\em reduced representation varieties}
\begin{equation}\label{eqn:Vzero}
\Vym{\ell}{0}_\mu=\{ (\ab)\in (L^I_\bR)^{2\ell} \mid \pab =\exp(X_\mu)\}
\cong X_{\mathrm{YM}}^{\ell,0}(L_\bR^I)_\mu.
\end{equation}
They correspond to the reduction from $G_\bR$ to the subgroup $L^I_\bR$.
More precisely, we have a homeomorphism
$$
\ym{\ell}{0}_\mu/G_\bR \cong  \Vym{\ell}{0}_\mu/L^I_\bR
$$
and a homotopy equivalence
$$
{\ym{\ell}{0}_\mu}^{h G_\bR} \sim  {\Vym{\ell}{0}_\mu}^{h L^I_\bR}
$$
where $X^{hG}$ denote the homotopic orbit space $EG\times_G X$.

We now recall the formulation in \cite[Section 2.1]{HL3}.
Let $\rho_{ss}:\widetilde{(L^I_\bR)_{ss}} \to (L^I_\bR)_{ss}$ be the universal
cover. Then the universal cover of $L^I_\bR$ is given by
$$
\rho:\widetilde{L^I_\bR}=\fz_{L^I_\bR}\times \widetilde{(L^I_\bR)_{ss}}
\to L^I_\bR,\quad
(X,g)\mapsto \exp_{Z^I}(X)\rho_{ss}(g)
$$
where $\exp_{Z^I}:\fz_{L^I_\bR}\to Z^I$ is the exponential map.
We have
$$
\pi_1((L^I_\bR)_{ss})\cong \Ker (\rho_{ss}),\quad
\pi_1(L^I_\bR)\cong \Ker\rho \subset (-2\pi\sqrt{-1}\Xi^I)\times Z(\widetilde{(L^I_\bR)_{ss}})
\subset \fz_{L^I_\bR}\times \widetilde{(L^I_\bR)_{ss}}.
$$
The map
$$
p_{L^I_\bR}:\Ker\rho \to \Xi^I,\quad
(X,g)\mapsto \frac{\sqrt{-1}}{2\pi}X
$$
coincides with the surjective group homomorphism
$$
p_{L^I_\bR}: \pi_1(L^I_\bR)\to \pi_1(L^I_\bR/(L^I_\bR)_{ss})
$$
under the isomorphisms
$\Ker\rho\cong \pi_1(L^I_\bR)$ and
$\Xi^I\cong \pi_1(L^I_\bR/(L^I_\bR)_{ss})$.

Define the obstruction map $o:\Vym{\ell}{0}_\mu \to p_{L^I_\bR}^{-1}(\mu)$ as follows.
Given a point $(a_1, b_1, \ldots, a_\ell, b_\ell)\in \Vym{\ell}{0}_\mu$, choose $\ta_i \in \rho^{-1}(a_i)$,
$\tb_i \in \rho^{-1}(b_i)$. Define
$o(\ab)=\prod_{i=1}^\ell[\ta_i,\tb_i]$. Note that this definition
does not depend on the choice of $\ta_i,\tb_i$.
We have
$o(\ab)\in \{0\}\times \widetilde{(L^I_\bR)_{ss}}$, and
\[
\rho_{ss}(o(\ab))= \exp(X_\mu).\]
More geometrically,
given $(\ab)\in \Vym{\ell}{0}_\mu$,
let $P$ be the underlying topological $L^I_\bR$-bundle. Then
$o(\ab)= o_2(P)$ under the identification $\pi_1(L^I_\bR)\cong H^2(\Si^\ell_0;\pi_1(L^I_\bR))$.
It is shown in \cite{HL3} that for $\ell\geq 1$, $o^{-1}(k)$ is nonempty
and connected for all $k\in p_{L^I_\bR}^{-1}(\mu)$. We conclude that
\begin{pro}
For any $I\subseteq \Delta$ and $\mu\in \Xi^I_+$, there is
a bijection
$$
\pi_0\left(\Vym{\ell}{0}_\mu\right) \cong p_{L^I_\bR}^{-1}(\mu).
$$
\end{pro}
Consider the short exact sequence of abelian groups:
$$
\begin{CD}
0 @>>> \pi_1((L^I_\bR)_{ss}) @>{i}>> \pi_1(L^I_\bR) @>{p_{L^I_\bR}}>> \pi_1(L^I_\bR/(L^I_\bR)_{ss}) @>>> 0\\
@. @| @| @|\\
 @. \widehat{\Lambda}^{L^I}/\Lambda^{L^I} @. \pi_1(H)/\Lambda^{L^I} @. \pi_1(H)/\widehat{\Lambda}^{L^I}
\end{CD}
$$
There is a bijection
$$
\pi_0(\Vym{\ell}{0}_\mu/L^I_\bR) \to p_{L^I_\bR}^{-1}(\mu).
$$
Given any $\beta \in p_{L^I_\bR}^{-1}(\mu)$, there is a bijection
$$
\pi_1((L^I_\bR)_{ss})\to p_{L^I_\bR}^{-1}(\mu),
\quad \alpha \mapsto i(\alpha) +\beta.
$$

\subsection{Equivariant Poincar\'{e} series}
\label{sec:poincare}

Given a $C^\infty$ principal $G$-bundle $\xi_0$ over
$\Si^\ell_0$, let
$$
\Xi_{\xi_0}=\bigl\{\mu\in \bigcup_{I\subseteq \Delta}\Xi_+^I\ \bigl|
\ \mu \textup{ is of Atiyah-Bott type for }\xi_0 \bigr. \bigr\}.
$$
The Harder-Narasimhan stratification of the space $\cC(\xi_0)$ of
$(0,1)$-connections on $\xi_0$ is given by
$$
\cC(\xi_0)=\bigcup_{\mu\in \Xi_{\xi_0}}\cC_\mu(\xi_0).
$$
Recall that $\cC(\xi_0)$ is an infinite dimensional
complex affine space, and each strata
$\cC_\mu(\xi_0)$ is a complex submanifold
of complex codimension
\begin{equation}\label{eqn:codim}
d_\mu=\sum_{\alpha(\mu)>0,\alpha\in R^+}(\alpha(\mu)+\ell-1)
\end{equation}

Let $P$ be a $C^\infty$ principal $G_\bR$-bundle
over $\Si^\ell_0$ such that $P\times_{G_\bR} G=\xi_0$,
and let $\cA(P)$ be the space of $G_\bR$-connections
on $P$. Then $\cA(P)\cong \cC(\xi_0)$ as infinite
dimensional complex affine spaces. In \cite{ym},
Atiyah and Bott conjectured that the Morse stratification
of the Yang-Mills functional on $\cA(P)$ exists
and coincides with the Harder-Narasimhan stratification
on $\cC(\xi_0)$ under the isomorphism $\cA(P)\cong \cC(\xi_0)$. The
conjecture was proved by Daskalopoulos in \cite{Da}.
Atiyah and Bott showed that the Harder-Narasimhan stratification
is $\cG(\xi_0)$-perfect over $\bQ$, where $\cG(\xi_0)=\Aut(\xi_0)$
is the gauge group of $\xi_0$. Therefore,
\begin{equation}\label{eqn:PtC}
P_t^{\cG(\xi_0)}(\cC(\xi_0);\bQ)=
\sum_{\mu\in \Xi_{\xi_0}} t^{2d_\mu}P_t^{\cG(\xi_0)}(\cC_\mu(\xi_0);\bQ).
\end{equation}

Let $\cA_\mu(P)\subset \cA(P)$ be the Morse stratum corresponding to
$\cC_\mu(\xi_0)\subset \cC(\xi_0)$. It is the stable manifold of
a connected component $\cN_\mu(P)$ of $\cN(P)$.
Let $(G_\bR)_\mu=(G_\bR)_{X_\mu}$. Then
$\mu$ and $P$ uniquely determine a topological principal
$(G_\bR)_\mu$-bundle $P_\mu$. Let $\ym{\ell}{0}_\mu^P$ denote the
connected component of $\ym{\ell}{0}_\mu$
which corresponds to $P\in \Prin_{G_\bR}(\Si^\ell_0)$,
and let $\Vym{\ell}{0}_\mu^{P_\mu}$ denote the
connected component of $\Vym{\ell}{0}_\mu$
which corresponds to $P_\mu \in \Prin_{(G_\bR)_\mu}(\Si^\ell_0)$.
Then  $\Vym{\ell}{0}_\mu^{P_\mu}$ can be identified
with the representation variety $V_{ss}(P_\mu)$ of {\em central}
Yang-Mills connections on $P_\mu$. We have homeomorphisms
$$
\cN_\mu(P)/\cG(P)\cong \ym{\ell}{0}^P_\mu/G_\bR
\cong V_{ss}(P_\mu)/(G_\bR)_\mu
$$
and homotopy equivalences of homotopic orbit spaces:
$$
{\cN_\mu(P)}^{h \cG(P)}\sim \left(\ym{\ell}{0}^P_\mu\right)^{h G_\bR}
\sim  V_{ss}(P_\mu)^{h (G_\bR)_\mu}.
$$
Combined with the homotopy equivalence
${\cC_\mu(\xi_0)}^{h \cG(\xi_0)}\sim {\cN_\mu(P)}^{h \cG(P)}$, we conclude that
$$
P_t^{\cG(\xi_0)}(\cC_\mu(\xi_0);\bQ)
=P_t^{G_\bR}(\ym{\ell}{0}^P_\mu;\bQ)
=P_t^{(G_\bR)_\mu}(V_{ss}(P_\mu);\bQ).
$$

\begin{rem}
The connectedness of $\cN_\mu(P)$ implies the
connectedness of $V_{ss}(P_\mu)$, but not
vise versa, because $\cG_0(P)$ is not connected
in general. We know $\cC_\mu(\xi_0)$
is connected by results in \cite{ym},
and $\cN_\mu(P) =\cN(P)\cap \cA_\mu(P)$ is a deformation
retract of $\cA_\mu(P) \cong \cC_\mu(\xi_0)$
by results in \cite{Da,Rad}, so
$\cN_\mu(P)$ is connected.
\end{rem}

Suppose that $\ell\geq 2$. Then there is a unique $\mu_0\in \Xi_{\xi_0}$
such that $d_{\mu_0}=0$. Then $\cC_{\mu_0}(\xi_0)=\cC_{ss}(\xi_0)$, the semi-stable
stratum. Let
$$
\cA_{ss}(P)=\cA_{\mu_0}(P),\quad
\cN_{ss}(P)=\cN_{\mu_0}(P),\quad
\Xi'_{\xi_0}=\Xi_{\xi_0}\setminus\{\mu_0\}.
$$
Then
$$
\cN_{ss}(P)/\cG_0(P) \cong V_{ss}(P).
$$
The identity \eqref{eqn:PtC} can be rewritten as
\begin{equation}\label{eqn:PtV}
P_t(B\cG(P);\bQ)=P_t^{G_\bR}(V_{ss}(P);\bQ)
+\sum_{\mu\in \Xi_{\xi_0}'} t^{2d_\mu}
P_t^{(G_\bR)_\mu}(V_{ss}(P_\mu);\bQ)
\end{equation}
where $P_t(B\cG(P);\bQ)$ is given by
Theorem \ref{thm:BG_orientable}. This allows one to compute
$$
P_t^{G_\bR}(V_{ss}(P);\bQ)
$$
recursively.

When $G=GL(n,\bC)$, equivalent inductive procedure was derived
by Harder and Narasimhan by number theoretic method in \cite{HN}.
Zagier provided an explicit closed formula which solves
the recursion relation for $GL(n,\bC)$ \cite{Za}. Laumon and Rapoport found
an explicit closed formula which solves the recursion relation
for general compact $G$. When $G_{ss}$ is not simply connected,
the recursion relation \cite[Theorem 3.2]{LR} that they solved
is not exactly the recursion relation \eqref{eqn:PtV}.
The closed formula which solves \eqref{eqn:PtV} is the following
slightly modified version of  \cite[Theorem 3.4]{LR} (see Appendix
\ref{sec:LR} for details):
\begin{thm}\label{thm:Pt}
Suppose that $\xi_0=P\times_{G_\bR}G$ and
$$
c_1(\xi_0)=\mu \in \pi_1(G)=\pi_1(H)/\Lambda.
$$
Then
\begin{eqnarray*}
\lefteqn{P^{G_\bR}_t(V_{ss}(P))=}\\ &&
\sum_{I\subseteq \Delta}
(-1)^{\dim_\bC \fz_{L^I} -\dim_\bC \fz_G}
P_t(B\cG^{L^I_\bR};\bQ)\frac{t^{2\dim_\bC U^I(\ell-1)}}
{\prod_{\alpha\in I}
(1-t^{4\langle \rho_I,\alpha^\vee\rangle}) }\cdot
t^{4\sum_{\alpha\in I}\langle \rho_I,\alpha^\vee\rangle
\langle \varpi_\alpha(\mu)\rangle}
\end{eqnarray*}
where
$$
\rho^I=\frac{1}{2}\sum_{\tiny \begin{array}{c}\beta\in R_+\\
 \langle \beta,\alpha^\vee \rangle >0
\textup{ for some }\alpha\in I\end{array} } \beta,
$$
$\varpi_\alpha(\mu)\in \bQ/\bZ$,
and $\langle x\rangle \in\bQ$ is the unique representative
of the class $x\in \bQ/\bZ$ such
 that $0< \langle x\rangle\leq 1$.
\end{thm}

Theorem \ref{thm:Pt} coincides with \cite[Theorem 3.2]{LR}
when $G_{ss}$ is simply connected, for example, when
$G_\bR=U(n)$, $G=GL(n,\bC)$. When $G_\bR=U(n)$,
Theorem \ref{thm:Pt} specializes to the closed formula
derived by Zagier in \cite{Za} (see \cite[Section 4]{LR}
for details):

\begin{thm}[{\cite{Za}, \cite[Section 4]{LR}}] \label{thm:PtU}
\begin{eqnarray*}
&& P^{U(n)}_t(\ymU{\ell}{0}_{\frac{k}{n},\ldots,\frac{k}{n}})\\
&=& \sum_{r=1}^n\sum_{\tiny \begin{array}{c}n_1,\ldots, n_r\in\bZ_{>0}\\ \sum n_j=n \end{array}}
(-1)^{r-1} \prod_{i=1}^r \frac{\prod_{j=1}^{n_i} (1+t^{2j-1})^{2\ell} }
{ (1-t^{2n_i})\prod_{j=1}^{n_i-1}(1-t^{2j})^2 }\\
&& \quad \cdot \frac{t^{2(\ell-1)\sum_{i<j} n_i n_j} }
{\prod_{i=1}^{r-1}(1-t^{2(n_i+n_{i+1})} ) }
\cdot t^{2\sum_{i=1}^{r-1} (n_i+n_{i+1}) \langle (n_1+\cdots + n_i)(-\frac{k}{n})\rangle}
\end{eqnarray*}
\end{thm}

\begin{rem}
For $n\geq 2$,  we have
\begin{eqnarray*}
&& P^{U(n)}_t\bigl(\ymU{\ell}{0}_{0,\ldots,0}\bigr)
= P^{U(n)}_t\bigl(\flU{\ell}{0}\bigr)\\
&=& P^{U(1)}_t\bigl(\flS{\ell}{0}\bigr)
P^{SU(n)}_t\bigl(X_{\mathrm{flat}}^{\ell,0}(SU(n))\bigr)
=\frac{(1+t)^{2\ell}}{1-t^2}
P^{SU(n)}_t\bigl(X_{\mathrm{flat}}^{\ell,0}(SU(n))\bigr)
\end{eqnarray*}
So Theorem \ref{thm:PtU} also gives a formula
for $P^{SU(n)}_t\bigl(X_{\mathrm{flat}}^{\ell,0}(SU(n))\bigr)$.
\end{rem}

\begin{ex}\label{UtwoSUtwo}
\begin{eqnarray*}
&& P^{U(2)}_t(X_{\mathrm{YM}}^{\ell,0}(U(2))_{\frac{k}{2},\frac{k}{2}})\\
&=& \frac{(1+t)^{2\ell}(1+t^3)^{2\ell} }{(1-t^4)(1-t^2)^2}
+(-1)\left(\frac{(1+t)^{2\ell}}{1-t^2}\right)^2\cdot \frac{t^{2(\ell-1)}}{1-t^4}
t^{4\langle{-\frac{k}{2} }\rangle }\\
&=& \frac{(1+t)^{2\ell}}{(1-t^2)^2(1-t^4)}
\left( (1+t^3)^{2\ell} - t^{2\ell-2 + 4\langle -\frac{k}{2}\rangle }(1+t)^{2\ell}\right)
\end{eqnarray*}
where
$$
\langle -k/2\rangle=\left\{
\begin{array}{ll}
1 & k \textup{ even}\\
1/2 & k \textup{ odd}
\end{array}\right.
$$
So
\begin{eqnarray*}
P^{U(2)}_t(X_{\mathrm{YM}}^{\ell,0}(U(2))_{\frac{k}{2},\frac{k}{2}})
&=& \left\{ \begin{array}{ll}
 \frac{(1+t)^{2\ell}}{(1-t^2)^2(1-t^4)}
\left( (1+t^3)^{2\ell} - t^{2\ell+2  }(1+t)^{2\ell}\right)
& k \textup{ even}\\
 \frac{(1+t)^{2\ell}}{(1-t^2)^2(1-t^4)}
\left( (1+t^3)^{2\ell} - t^{2\ell  }(1+t)^{2\ell}\right)
& k \textup{ odd} \end{array}\right.
\end{eqnarray*}
and
$$
P^{SU(2)}_t(X_{\mathrm{flat}}^{\ell,0}(SU(2)))
= \frac{(1+t^3)^{2\ell}}{(1-t^2)(1-t^4)}
- \frac{t^{2\ell+2}(1+t)^{2\ell}}{(1-t^2)(1-t^4)}.
$$
\end{ex}

\begin{ex}\label{SUthreeSUfour}

\begin{eqnarray*}
&& P^{SU(3)}_t(X_{\mathrm{flat}}^{\ell,0}(SU(3)))\\
&=& \frac{(1+t^3)^{2\ell}(1+t^5)^{2\ell} }{(1-t^2)(1-t^4)^2(1-t^6)}
-2\frac{(1+t)^{2\ell}(1+t^3)^{2\ell} t^{4\ell+2}}
{(1-t^2)^2(1-t^4)(1-t^6)}
+\frac{(1+t)^{4\ell} t^{6\ell+2} }{(1-t^2)^2(1-t^4)^2}
\end{eqnarray*}
\begin{eqnarray*}
&& P^{SU(4)}_t(X_{\mathrm{flat}}^{\ell,0}(SU(4)))\\
&=& \frac{(1+t^3)^{2\ell}(1+t^5)^{2\ell}(1+t^7)^{2\ell}}
{(1-t^2)(1-t^4)^2(1-t^6)^2(1-t^8)}
-2 \frac{(1+t)^{2\ell}(1+t^3)^{2\ell}(1+t^5)^{2\ell} t^{6\ell+2}}
{(1-t^2)^2(1-t^4)^2(1-t^6)(1-t^8)}\\
&& -\frac{(1+t)^{2\ell}(1+t^3)^{4\ell}t^{8\ell}}
{(1-t^2)^3(1-t^4)^2(1-t^8)}
+2\frac{(1+t)^{4\ell}(1+t^3)^{2\ell}t^{10\ell}}
{(1-t^2)^3(1-t^4)^2(1-t^6)}\\
&& +\frac{(1+t)^{4\ell}(1+t^3)^{2\ell}t^{10\ell+2}}
{(1-t^2)^3(1-t^4)(1-t^6)^2}
-\frac{(1+t)^{6\ell}t^{12\ell}}{(1-t^2)^3(1-t^4)^3}
\end{eqnarray*}
\end{ex}

We will use Theorem \ref{thm:Pt} to write down explicit
closed formula for
$SO(2n+1)$, $SO(2n)$, and $Sp(n)$
in Section \ref{sec:SOodd_poincare},
Section \ref{sec:SOeven_poincare},
and Section \ref{sec:Sp_poincare},
respectively.

\subsection{Involution on the Weyl Chamber}
Let $\pi:\tSi\to\Si$ be the orientable double cover of
a closed, compact, connected, nonorientable surface $\Si$,
and let $\tau:\tSi\to \tSi$ be the deck transformation.
Let $P$ be a principal $G_\bR$-bundle over $\Sigma$, and
let $\tP=\pi^*P$.  Then $\tP$ and $\xi_0=\tP\times_{G_\bR} G$ are topologically trivial.
There is an involution $\tilde{\tau}_s:\tP\to\tP$
which covers the anti-holomorphic involution $\tau:\tSi\to \tSi$.
Under the trivialization $\tP\cong \tSi\times G_\bR$, $\tau_s$ is given
by $(x,h)\mapsto (\tau(x), s(x)h)$, where $s:\tSi\to G_\bR$ satisfies
$s(\tau(x))=s(x)^{-1}$ (see \cite[Section 3.2]{HL4} for details).

Let $\cA(P)$ and $\cA(\tP)$ denote the space of $G_\bR$-connections
on $P$ and on $\tP$ respectively, and let $\cC(\xi_0)$ be the space
of $(0,1)$-connections on the principal $G$-bundle $\xi_0$. Then
$\tilde{\tau}_s$ induces an involution $\tilde{\tau}_s^*:\cA(\tP)
\to \cA(\tP)$. Since $\tP$ and $\xi_0$ are topologically trivial, we
may identify $\cA(\tP)$ with $\Omega^1(\tSi,\fg_\bR)$ and identify
$\cC(\xi_0)$ with $\Omega^{0,1}(\tSi,\fg)$. Let
$j:\Omega^1(\tSi,\fg_\bR)\to \Omega^{0,1}(\tSi,\fg)$ be defined as
in the first paragraph of Section \ref{sec:algebraic_geometry}.
Given $X=X_1+\sqrt{-1} X_2\in \fg$, where $X_1,X_2\in\fg_\bR$,
define $\bar{X}= X_1-\sqrt{-1}X_2$; given $X:\tSi\to \fg$, define
$\bar{X}:\tSi\to \fg$ by $x\mapsto \overline{X(x)}$. Then $j\circ
\tilde{\tau}^*_s \circ j^{-1}: \cC(\xi_0)\to \cC(\xi_0)$ is given by
$$
 X\otimes \theta \mapsto \Ad(s)\overline{\tau^*X}\otimes \overline{\tau^*\theta}
 $$
 where $X\in\Omega^0(\tSi,\fg)$ and $\theta\in\Omega^{0,1}(\tSi)$.
 From now on, we denote $j\circ \tilde{\tau}_s^*\circ j^{-1}$ by $\tilde{\tau}_s^*$. We have
isomorphisms of real affine spaces
$\cA(P)\cong \cA(\tP)^{\tilde{\tau}_s^*}\cong \cC(\xi_0)^{\tilde{\tau}_s^*}$.

We will define an involution $\tau'$ on the positive Weyl chamber
$\BC_0$ such that $\tilde{\tau}_s^*\cC_\mu = \cC_{\tau'(\mu)}$,
where $\mu \in \BC_0$ is of Atiyah-Bott type for $\xi_0$ and $\cC_\mu$
is the associated stratum in $\cC(\xi_0)$.

The set
$$
-C_0=\{ -Y \mid Y\in C_0\} \subset \sqrt{-1}\fh_\bR.
$$
is another Weyl chamber. There is a unique element $w$ in the Weyl group
$W$ such that $w\cdot C_0=-C_0$. We have $w^2 \cdot C_0=C_0$, so
$w^2=\mathrm{id}_{\sqrt{-1}\fh_{\bR}}$.
Define $\tau': \sqrt{-1}\fh_\bR \to \sqrt{-1}\fh_\bR$ by
$X\mapsto  w\cdot (-X)$. Recall that $\tau$ induces an involution
on the symmetric representation variety which maps $X\in\fg_\bR$ to
$-\Ad(\bar{c})X \in\fg_\bR$ (see \cite[Section 4.5]{HL4}).
Given $Y\in \BC_0$, $\tau'(Y)$ is the unique vector in $\BC_0$ which
is in the orbit $G\cdot(-Y)=G\cdot(-\Ad(\bar{c})(Y))$ of
the adjoint action of $G$ on $\fg$.
Thus $\tau'$ is induced by the involution $\tau$ on the symmetric
representation variety. To simplify notation, from now on we will write $\tau$ instead of $\tau'$.
Obviously $\tau(\BC_0)=\BC_0$. Given $Y\in \BC_0$,
$\tau(Y)=Y$ if and only if $Y\in \BC_0$ is conjugate to $-Y$.
In this case, we have
$\Ad(\ep)Y=-Y$, where $\ep\in N(H_\bR)\subset G_\bR$ represents $w\in W=N(H_\bR)/H_\bR$.

To demonstrate the above discussion, we list some examples of classical Lie groups.
\begin{ex}
Let $G_\bR=U(n)$. then
\begin{eqnarray*}
\BC_0 &=&\{\diag(t_1,\ldots,t_n)\mid t_1,\ldots,t_n\in\bR, t_1\geq\cdots\geq t_n\}\\
-\BC_0&=&\{\diag(v_1,\ldots,v_n) \mid v_1,\ldots,v_n\in\bR, v_1 \leq \cdots \leq v_n\}
\end{eqnarray*}
There exists a unique $w$ in $W\cong S(n)$ , the symmetric group,
such that $w(\BC_0)=-\BC_0$. In fact,
$w\cdot \diag(t_1,\ldots,t_n)=\diag(t_n,\ldots,t_1)$ is the action of such $w$ on $\sqrt{-1}\fh_\bR$.
Thus, the involution $\tau(Y)$ defined as $w\cdot(-Y)$ gives us $\tau(\diag(t_1,\ldots,t_n))=\diag(-t_n,\ldots,-t_1)$,
and $Y$ is conjugate to $-Y$ (i.e. $\tau(Y)=Y$) if and only if $(t_1,\ldots,t_n)=(-t_n,\ldots,-t_1)$,
or equivalently, if and only if $Y$ is of the form $\diag(v_1,\ldots,v_k,0,\ldots,0,-v_k,\ldots,-v_1)$.
\end{ex}

\begin{ex}\label{ISO(2n+1)}
Let $G_\bR=SO(2n+1)$. then
\begin{eqnarray*}
\BC_0 &=&\{\sqrt{-1}\diag(t_1 J,\ldots,t_n J,0 I_1)\mid t_1 \geq\cdots\geq t_n\geq 0\},\\
-\BC_0&=&\{\sqrt{-1}\diag(v_1 J,\ldots,v_n J,0 I_1) \mid v_1 \leq \cdots \leq v_n \leq 0\},
\end{eqnarray*}
where
$$
J=\left(\begin{array}{cc} 0& -1\\1&0\end{array}\right).
$$
The unique $w$ in $W\cong G(n)$, the wreath product of $\bZ_2$ by $S(n)$,
that maps $\BC_0$ to $-\BC_0$, acts as
$w\cdot\sqrt{-1}\diag(t_1 J,\ldots,t_n J,0I_1)=\sqrt{-1}\diag(-t_1 J,\ldots,-t_n J,0I_1)$. Thus
$\tau:\sqrt{-1}\fh_\bR\to \sqrt{-1}\fh_\bR$ is the identity map.
Any $Y\in \BC_0$ is conjugate to the negative of itself.
Let
$$
H=\left(\begin{array}{rr} 1& 0\\0&-1 \end{array}\right)
$$
and let $H_n=\diag(\underbrace{H,\ldots,H}_n)$.
The element
$$
\ep=\diag(H_n, (-1)^n )\in SO(2n+1)
$$
satisfies $\Ad(\epsilon)Y=-Y$ for all $Y\in \sqrt{-1}\fh_\bR$ and $\ep^2=e$.
\end{ex}

\begin{ex}\label{ISO(2n)}
Let $G_\bR=SO(2n)$. Then
\begin{eqnarray*}
\BC_0 &=&\{\sqrt{-1}\diag(t_1 J,\ldots,t_n J)\mid t_1\geq\cdots\geq\mid t_n\mid\geq 0\},\\
-\BC_0 &=&\{\sqrt{-1}\diag(v_1 J,\ldots,v_n J) \mid v_1 \leq \cdots \leq -\mid v_n \mid \leq 0\}.
\end{eqnarray*}
The unique $w$ in $W\cong SG(n)$, the subgroup of $G(n)$ consisting of even permutations,
that maps $\BC_0$ to $-\BC_0$, belongs to the $\bZ_2$ part of $SG(n)$, and
$$
w\cdot \sqrt{-1}\diag(t_1 J,\ldots,t_n J)
= \sqrt{-1}\diag(-t_1 J,\ldots,-t_{n-1}J,(-1)^{n-1}t_n J).
$$
Thus
$$
\tau\left(\sqrt{-1}\diag(t_1 J,\ldots,t_n J)\right)
=\sqrt{-1} \diag(t_1 J,\ldots,t_{n-1}J,(-1)^n t_n J)
$$
If $n$ is even, then any $Y\in\BC_0$ in conjugate to the negative of
itself. If $n$ is odd, then $Y\in \BC_0$ is conjugate to $-Y$ if and
only if $Y$ is of the form $\sqrt{-1}\diag(t_1
J,\ldots,t_{n-1}J,0)$. Define
$$
\ep=\left\{\begin{array}{c}
H_n\\
\diag(H_{n-1},I_2) \end{array}\right.\in SO(2n)
\quad \textup{if $n$ is}
\begin{array}{l}\textup{even}\\ \textup{odd}\end{array}
$$
Then $\ep$ satisfies
$\Ad(\ep)Y=-Y$ for all $Y\in \sqrt{-1}\fh_\bR$ and $\epsilon^2=e$.
\end{ex}

\begin{ex}\label{ISp(n)}
Let $G_\bR=Sp(n)$. Then
\begin{eqnarray*}
\BC_0&=&\{\diag(t_1,\ldots,t_n,-t_1,\ldots,-t_n)\mid t_1\geq\cdots\geq t_n\geq 0\},\\
-\BC_0&=&\{\diag(v_1,\ldots,v_n,-v_1,\ldots,-v_n) \mid v_1 \leq \cdots \leq v_n \leq 0\}.
\end{eqnarray*}
the unique $w$ in $W\cong G(n)$, the wreath product of $\bZ_2$ by $S(n)$,
that maps $\BC_0$ to $-\BC_0$, acts as
$w\cdot\diag(t_1,\ldots,t_n,-t_1,\ldots,-t_n)=\diag(-t_1,\ldots,-t_n,t_1,\ldots,t_n)$.
Thus $\tau:\sqrt{-1}\fh_\bR \to \sqrt{-1}\fh_\bR$ is the identity map, and any
$Y\in\BC_0$ is conjugate to the negative of itself just as
in the $SO(2n+1)$ case. The element
$$
\epsilon=\left(\begin{array}{cc}
 0  &  -I_{n}   \\
 I_{n}  &  0
\end{array}\right) \in Sp(n)
$$
satisfies $Ad(\ep)Y=-Y$ for all $Y\in \sqrt{-1}\fh_\bR$  but $\ep^2\neq e$.
Indeed, let $\tilde{\ep}$ be any element that satisfies $\Ad(\tilde{\ep})Y=-Y$ for all $Y\in \sqrt{-1}\fh_\bR$.
Then we must have $\tilde{\ep}=\ep u$ for some $u$ in the maximal torus, and it is straightforward to
check that $\tilde{\ep}^2= -e$.
\end{ex}

\subsection[Connectedness of the representations for nonorientable sufaces]
{Connected components of the representation variety for nonorientable surfaces}
\label{sec:connected}

$G_\bR$ is connected, so the natural projection
$$
\ym{\ell}{i}\to \ym{\ell}{i}/G_\bR
$$
induces a bijection
$$
\pi_0(\ym{\ell}{i})
\to \pi_0(\ym{\ell}{i})/G_\bR).
$$

Any point in $\ym{\ell}{1}/G_\bR$ can
be represented uniquely by
$$
(\ab,c,X)
$$
where $X\in \BC_0$. Moreover, we must have
$X\in\BC_0^\tau$. Similarly, any point in $\ym{\ell}{2}/G_\bR$
can be represented uniquely by
$$
(\ab,d,c,X)
$$
where $X\in \BC_0^\tau$.

Recall that $\tau$ is an $\bR$-linear map from $\sqrt{-1}\fh_\bR$ to $\sqrt{-1}\fh_\bR$.
Its dual $\tau^*$ is an $\bR$-linear map
from $(\sqrt{-1}\fh_\bR)^\vee$ to $(\sqrt{-1}\fh_\bR)^\vee$. This
$\tau^*$ preserves $\Delta$, the set of simple roots,
and restricts to an involution on it. To simplify notation, we will also denote this involution
by $\tau$. Given $I\subseteq \Delta$ such that
$\tau(I)=I$, let
$$
(\Xi^I_+)^\tau =\{ \mu \in \Xi^I_+ \mid \tau(\mu)=\mu\}
$$

Suppose that $(\ab,c,X)\in \ym{\ell}{1}$. Then there
is a unique pair $(\mu,I)$, where $I\subseteq \Delta$, $\tau(I)=I$,
and $\mu\in (\Xi^I_+)^\tau$, such that $X$ is conjugate to $X_\mu =-2\pi\sqrt{-1}\mu$.
Given $\mu\in (\Xi^I_+)^\tau$, where $I\subseteq \Delta$ and $\tau(I)=I$, define
\begin{eqnarray*}
\ym{\ell}{1}_\mu &=&  \{(\ab,c,X)\in {G_\bR}^{2\ell+1}\times C_{\mu/2}\mid \\
&&\ab\in (G_\bR)_X, \Ad(c)X=-X, \pab= \exp(X)c^2 \}
\end{eqnarray*}
Where $C_{\mu/2}$ is the conjugacy class of $X_\mu/2$. We define $\ym{\ell}{2}_\mu$ similarly. For $i=1,2$,
$\ym{\ell}{i}$ is a disjoint union of
$$
\{ \ym{\ell}{i}_\mu \mid  \mu\in (\Xi^I_+)^\tau, I \subseteq \Delta, \tau(I)=I\}.
$$
When $G_\bR=U(n)$, $\ell\geq 1$,
each $\ym{\ell}{i}_\mu$ is nonempty and has one or two connected
components (see \cite[Section 7]{HL4}). We will see later that
$\ym{\ell}{i}_\mu$ can be empty for other classical groups (Section \ref{sec:SOodd_nonorientable},
Section \ref{sec:SO4m+2_nonorientable}, Section \ref{sec:SO4m_nonorientable} and
Section \ref{sec:Sp_nonorientable}).
When $\ym{\ell}{i}_\mu$ is nonempty, it is a union of finitely many connected
components of $\ym{\ell}{i}$.

The reduction of $\ym{\ell}{i}_\mu$ is more complicated
because $c$ is not in $G_X$. To do the reduction, we fix
some $\ep\in G_\bR$ such that the involution on $\BC_0$ is given
by $X\mapsto -\Ad(\ep)X$. Thus $\Ad(\ep)X=-X$ if $X$ is fixed by
the involution. For
any $\mu \in (\Xi^I_+)^\tau$, where $\tau(I)=I$,
we define {\em $\ep$-reduced representation varieties}
\begin{equation}\label{eqn:Vone}
\Vym{\ell}{1}_\mu=\{ (\ab,c')\in (L^I_\bR)^{2\ell+1} \mid
\pab=\exp(\frac{X_\mu}{2}) \ep c'  \ep  c' \}
\end{equation}
\begin{equation}\label{eqn:Vtwo}
\Vym{\ell}{2}_\mu=\{ (\ab,d,c')\in (L^I_\bR)^{2\ell+2} \mid
\pab=\exp(\frac{X_\mu}{2}) \ep c' d (\ep c')^{-1}  d \}
\end{equation}
For $i=1,2$, $L^I_\bR$ acts on $\Vym{\ell}{i}_\mu$ by
$$
g\cdot (c_1,\ldots, c_{2\ell+i})= (g c_1 g^{-1},\ldots, g c_{2\ell+i-1} g^{-1},
\ep^{-1}g \ep c_{2\ell+i} g^{-1}).
$$
Recall that $\Ad(\ep)(X_\mu)=-X_\mu$ and $L^I_\bR=(G_\bR)_{X_\mu}$. So
we have a homeomorphism
$$
\ym{\ell}{i}_\mu/G_\bR\cong \Vym{\ell}{i}_\mu/L^I_\bR
$$
and a homotopy equivalence between homotopic orbit spaces:
$$
{\ym{\ell}{i}_\mu}^{h G_\bR} \sim  {\Vym{\ell}{i}_\mu}^{h L^I_\bR}
$$

When $G_\bR=U(n)$, $\Vym{\ell}{i}_\mu$ can be viewed as a product of
representation varieties for $U(m)$ ($m<n$) of $\Si^\ell_i$ and of
its double cover $\Si^{2\ell+i-1}_0$ (see \cite[Section 7]{HL4}.
This is not the case for other classical groups.
We will see in Section \ref{sec:SOodd_nonorientable},
Section \ref{sec:SO4m+2_nonorientable}, Section \ref{sec:SO4m_nonorientable}, and
Section \ref{sec:Sp_nonorientable} that when $G_\bR=SO(n)$ or $Sp(n)$,
$\Vym{\ell}{i}_\mu$ is a product of {\em twisted representation varieties}
defined in Section \ref{sec:twisted} and Section \ref{sec:twistedO} below.

\subsection{Twisted representation varieties: $U(n)$}
\label{sec:twisted}

Given $n, k\in \bZ$, $n>0$, define {\em twisted representation varieties}
\begin{equation}\label{eqn:twistI}
\tV_{n,k}^{\ell,1}=\Bigl\{(\ab,c)
\in U(n)^{2\ell+1}\mid \pab=e^{-2\pi\sqrt{-1}k/n}I_n \bar{c}c
  \Bigr \}
\end{equation}
\begin{equation}\label{eqn:twistII}
\tV_{n,k}^{\ell,2}=\Bigl\{(\ab,d,c)\in U(n)^{2\ell+2}\mid \pab=
e^{-2\pi\sqrt{-1}k/n} I_n \bar{c}\bar{d}\bar{c}^{-1}d \}
\end{equation}
where $\bar{c}$ is the complex conjugate of $c$. In particular,
$$
\tV_{1,k}^{\ell,1}=U(1)^{2\ell+1},\quad
\tV_{1,k}^{\ell,2}=U(1)^{2\ell+2}.
$$

For $i=1,2$, $U(n)$ acts on $\tV_{n,k}^{\ell,i}$ by
\begin{equation}\label{eqn:actI}
g\cdot(\ab,c)= (g a_1 g^{-1}, gb_1 g^{-1},\ldots, g a_\ell g^{-1}, g b_\ell g^{-1},
\bar{g} c g^{-1})
\end{equation}
\begin{equation}\label{eqn:actII}
g\cdot(\ab,d, c)= (g a_1 g^{-1}, gb_1 g^{-1},\ldots, g a_\ell g^{-1}, g b_\ell g^{-1},
g d g^{-1}, \bar{g} c g^{-1})
\end{equation}

We will show that
\begin{pro}\label{thm:tV}
$\tV_{n,k}^{\ell,i}$ is nonempty and connected for $\ell\geq 2i$.
\end{pro}

\begin{proof}[Proof for $i=1$]
For any $(\ab,c)\in \tV_{n,k}^{\ell,1}$, we have
$$
\det(a_i)= e^{\sqrt{-1}\theta_i},\quad
\det(b_i)= e^{\sqrt{-1}\phi_i},\quad
\det(c)=e^{\sqrt{-1}\theta}.
$$
Define $\beta:[0,1] \to U(n)^{2\ell+1}$ by
$$
\beta(t)=(e^{-\sqrt{-1}t\theta_1/n}a_1,
e^{-\sqrt{-1}t\phi_1/n}b_1,\ldots,
e^{-\sqrt{-1}t\theta_\ell/n}a_\ell,
e^{-\sqrt{-1}t\phi_\ell/n}b_\ell,
e^{-\sqrt{-1}t\theta/n} c).
$$
Then the image of $\beta$ lies
in $\tV_{n,k}^{\ell,1}$, $\beta(0)=(\ab,c)$, and
\begin{eqnarray*}
\beta(1)\in W_{n,k}^{\ell,1} & \stackrel{\mathrm{def}}{=} &
\Bigl\{(\ab,c)\in SU(n)^{2\ell+1} \bigr|
\\ && \quad \pab=e^{-2\pi\sqrt{-1}k/n}I_n \bar{c}c\Bigr \} \subset \tV_{n,k}^{\ell,1}.
\end{eqnarray*}
So it suffices to show that $W_{n,k}^{\ell,1}$ is nonempty and connected.

Define $\pi:W_{n,k}^{\ell,1}\to SU(n)$ by
$(\ab,c)\mapsto c$. Then $\pi^{-1}(c)$ is nonempty and connected for any $c\in SU(n)$.
It remains to show that for any $c\in SU(n)$, there is a path
$\gamma:[0,1]\to W_{n,k}^{\ell,1}$ such
that $\gamma(0)\in \pi^{-1}(e)$ and $\gamma(1)\in \pi^{-1}(c)$.

Let $T$ be the maximal torus which consists of diagonal matrices in $SU(n)$.
For any $c\in SU(n)$, there exist
$g \in SU(n)$ such that $g^{-1} c g \in T$. We have
$$
c=g \exp \xi g^{-1}, \quad \bar{c} =\bar{g} \exp(-\xi) \bar{g}^{-1}
$$
for some $\xi\in\ft$.  Let
$$
\xi_0=-2\pi\sqrt{-1}\frac{k}{n}\diag(I_{n-1},(1-n)I_1) \in \ft.
$$
Then $\exp(\xi_0)=e^{-2\pi\sqrt{-1}k/n}I_n$.
Let $\omega$ be the coxeter element and $a$ be the corresponding element in $SU(n)$.
There are $\eta_0,\eta\in \ft$ such that
$$
\omega\cdot\eta_0-\eta_0=\xi_0,\quad
\omega\cdot\eta-\eta =\xi.
$$
Let $a\in N(T)$ represent $\omega\in W=N(T)/T$. Then
\begin{eqnarray*}
a\exp(\eta_0 - t \eta) a^{-1}\exp(-\eta_0 + t\eta)&=&
\exp(\omega\cdot(\eta_0 - t\eta)-(\eta_0 - t\eta))\\
&=&\exp(\xi_0-t\xi)\\
&=& e^{-2\pi\sqrt{-1}k/n}\exp(-t\xi)\\
a\exp(t\eta)a^{-1}\exp(-t\eta)&=&\exp(\omega\cdot(t\eta)-t\eta)=\exp(t\xi).
\end{eqnarray*}

Now since $SU(n)$ is connected, there are paths $\tg:[0,1]\to SU(n)$
such that $\tg(0)=e$ and $\tg(1)=g$. Now define $\gamma:[0,1]\to SU(n)^{2\ell+1}$
by $$\gamma(t)=(a_1(t),b_1(t),a_2(t),b_2(t),e,\ldots,e,c(t))$$
where
\begin{eqnarray*}
&& a_1(t)=\overline{\tg(t)}a\left(\overline{\tg(t)}\right)^{-1},\quad
b_1(t)=\overline{\tg(t)}\exp(\eta_0-t\eta)\left(\overline{\tg(t)}\right)^{-1}\\
&& a_2(t)=\tg(t) a \tg(t)^{-1},\quad
b_2(t)=\tg(t)\exp(t\eta)\tg(t)^{-1}\\
&& c(t)=\tg(t)\exp(t\xi)\tg(t)^{-1}.
\end{eqnarray*}
Then
$$
[a_1(t),b_1(t)]=e^{-2\pi\sqrt{-1}k/n} \overline{c(t)},\quad
[a_2(t),b_2(t)]=c(t),
$$
so the image of $\gamma$ lies in $W^{\ell,1}_{n,k}$.
We have
\begin{eqnarray*}
\gamma(0)&=&(a,\exp(\eta_0),a,e,e,\ldots,e,e,e)\in \pi^{-1}(e)\\
\gamma(1)&=&(\bar{g} a \bar{g} ^{-1},\bar{g} \exp(\eta_0-\eta)\bar{g}^{-1},
g a g^{-1},g \exp(\eta) \bar{g}^{-1}, e,\cdots,e,c)\in \pi^{-1}(c).
\end{eqnarray*}
\end{proof}

\begin{proof}[Proof for $i=2$]
For any $(\ab,d,c)\in \tV_{n,k}^{\ell,2}$, we have
$$
\det(a_i)= e^{\sqrt{-1}\theta_i},\quad
\det(b_i)= e^{\sqrt{-1}\phi_i},\quad
\det(c)=e^{\sqrt{-1}\theta},\quad
\det(d)=e^{\sqrt{-1}\phi}.
$$
Define $\beta:[0,1] \to U(n)^{2\ell+2}$ by
\begin{eqnarray*}
\lefteqn{
\beta(t)=
(e^{-\sqrt{-1}t\theta_1/n}a_1,
e^{-\sqrt{-1}t\phi_1/n}b_1,\ldots,}\\&&
e^{-\sqrt{-1}t\theta_\ell/n}a_\ell,
e^{-\sqrt{-1}t\phi_\ell/n}b_\ell,
e^{-\sqrt{-1}t\phi/n}d,
e^{-\sqrt{-1}t\theta/n} c).
\end{eqnarray*}
Then the image of $\beta$ lies
in $\tV_{n,k}^{\ell,2}$, $\beta(0)=(\ab,d,c)$, and
\begin{eqnarray*}
\beta(1)\in W_{n,k}^{\ell,2} &\stackrel{\mathrm{def}}{=}&
\Bigl\{(\ab,d,c)\in SU(n)^{2\ell+2} \bigr|\\
&& \quad \pab=e^{-2\pi\sqrt{-1}k/n}I_n \bar{c}\bar{d}\bar{c}^{-1}d\Bigr \}
\subset \tV_{n,k}^{\ell,2}.
\end{eqnarray*}
So it suffices to show that $W_{n,k}^{\ell,2}$ is nonempty and connected.

Define $\pi:W_{n,k}^{\ell,2}\to SU(n)^2$ by
$(\ab,d,c)\mapsto (d,c)$. Then $\pi^{-1}(d,c)$ is nonempty and connected for any $(d,c)\in SU(n)^2$.
It remains to show that for any $(d,c)\in SU(n)^2$, there is a path
$\gamma:[0,1]\to W_{n,k}^{\ell,2}$ such
that $\gamma(0)\in \pi^{-1}(e,e)$ and $\gamma(1)\in \pi^{-1}(d,c)$.

Let $T$ be the maximal torus which consists of diagonal matrices in $SU(n)$.
For any $c,d\in SU(n)$, there exist
$g_1,g_2 \in SU(n)$ such that $g_1^{-1} c g_1, g_2^{-1} d g_2 \in T$. We have
$$
c=g_1 \exp \xi_1 g_1^{-1}, \quad
\bar{c} =\bar{g}_1 \exp(-\xi_1) \bar{g}_1^{-1},\quad
d=g_2 \exp \xi_2 g_2^{-1}, \quad
\bar{d} =\bar{g}_2 \exp(-\xi_2) \bar{g}_2^{-1}
$$
for some $\xi_1,\xi_2\in\ft$.  Let
$$
\xi_0=-2\pi\sqrt{-1}\frac{k}{n}\diag(I_{n-1},(1-n)I_1) \in \ft.
$$
Then $\exp(\xi_0)=e^{-2\pi\sqrt{-1}k/n}I_n$.
Let $\omega$ be the coxeter element and $a$ be the corresponding element in $SU(n)$.
There are $\eta_0,\eta_1,\eta_2\in \ft$ such that
$$
\omega\cdot\eta_j-\eta_j=\xi_j,\quad j=0,1,2.
$$
Let $a\in N(T)$ represent $\omega\in W=N(T)/T$. Then
\begin{eqnarray*}
a\exp(\eta_0-t \eta_1) a^{-1}\exp(-\eta_0+t\eta_1)&=&
\exp(\omega\cdot(\eta_0-t\eta_1)-(\eta_0-t\eta_1))\\
&=&\exp(\xi_0 - t\xi_1)\\
&=& e^{-2\pi\sqrt{-1}k/n}\exp(-t\xi_1)\\
a\exp(-t\eta_1)a^{-1}\exp(t\eta_1)&=&\exp(\omega\cdot(-t\eta_1)+t\eta_1)=\exp(-t\xi_1)\\
a\exp(-t\eta_2)a^{-1}\exp(t\eta_2)&=&\exp(\omega\cdot(-t\eta_2)+t\eta_2)=\exp(-t\xi_2).
\end{eqnarray*}

Now since $SU(n)$ is connected, there are paths $\tg_1,\tg_2:[0,1]\to SU(n)$
such that $\tg_j(0)=e$ and $\tg_j(1)=g_j$ for $j=1,2$. Now define $\gamma:[0,1]\to SU(n)^{2\ell+2}$
by $$\gamma(t)=(a_1(t),b_1(t),a_2(t),b_2(t),a_3(t), b_3(t), a_4(t), b_4(t),e,\ldots,e,d(t),c(t))$$
where
\begin{eqnarray*}
&& a_1(t)=\overline{\tg_1(t)}a\left(\overline{\tg_1(t)}\right)^{-1},\quad
b_1(t)=\overline{\tg_1(t)}\exp(\eta_0-t\eta_1)\left(\overline{\tg_1(t)}\right)^{-1}\\
&& a_2(t)=\overline{\tg_2(t)} a \left(\overline{\tg_2(t)}\right)^{-1},\quad
b_2(t)=\overline{\tg_2(t)}\exp(-t\eta_2)\left(\overline{\tg_2(t)}\right)^{-1}\\
&& a_3(t)=\overline{\tg_1(t)}a\left(\overline{\tg_1(t)}\right)^{-1},\quad
b_3(t)=\overline{\tg_1(t)}\exp(t\eta_1)\left(\overline{\tg_1(t)}\right)^{-1}\\
&& a_4(t)=\tg_2(t) a \left(\tg_2(t)\right)^{-1},\quad
b_4(t)=\tg_2(t)\exp(t\eta_2)\left(\tg_2(t)\right)^{-1}\\
&& c(t)=\tg_1(t)\exp(t\xi_1)\tg_1(t)^{-1},\quad d(t)=\tg_2(t)\exp(t\xi_2)\tg_2(t)^{-1}
\end{eqnarray*}
Then
\begin{eqnarray*}
&& [a_1(t),b_1(t)]=e^{-2\pi\sqrt{-1}k/n} \overline{c(t)},\quad
   [a_2(t),b_2(t)]=\overline{d(t)},\\
&& [a_3(t),b_3(t)]=\overline{c(t)}^{-1},\quad
   [a_4(t),b_4(t)]= d(t).
\end{eqnarray*}
so the image of $\gamma$ lies in $W^{\ell,2}_{n,k}$.
We have
\begin{eqnarray*}
\gamma(0)&=&(a,\exp(\eta_0),a,e,a,e,a,e,e,\ldots,e,e,e)\in \pi^{-1}(e,e)\\
\gamma(1)&=&(\overline{g_1} a \overline{g_1} ^{-1}, \overline{g_1} \exp(\eta_0-\eta_1)\overline{g_1}^{-1},
\overline{g_2} a \overline{g_2}^{-1},\overline{g_2}\exp(-\eta_2) \overline{g_2}^{-1},\\
&& \overline{g_1} a \overline{g_1} ^{-1}, \overline{g_1} \exp(\eta_1)\overline{g_1}^{-1},
g_2 a g_2^{-1},g_2\exp(\eta_2) g_2^{-1},
e,\ldots,e,d,c)\in \pi^{-1}(d,c).
\end{eqnarray*}
\end{proof}

\subsection{Twisted representation varieties: $SO(n)$}
\label{sec:twistedO}

\newcommand{\tc}{\tilde{c}}
\newcommand{\td}{\tilde{d}}
\newcommand{\tab}{\ta_1,\tb_1,\ldots,\ta_\ell, \tb_\ell}
\newcommand{\ptab}{\prod_{i=1}^\ell[\ta_i,\tb_i]}

Let
$$
O(n)_\pm =\{ A\in O(n)\mid \det(A)=\pm 1\}.
$$
Then $O(n)_+$ and $O(n)_-$ are the two connected components
of $O(n)$, where $O(n)_+=SO(n)$. For $n\geq 2$, define
\begin{equation}\label{eqn:Oone}
V^{\ell,1}_{O(n),\pm 1} = \{(\ab,c)\in SO(n)^{2\ell}\times O(n)_\pm\mid
\pab=c^2\}
\end{equation}
\begin{equation}\label{eqn:Otwo}
V^{\ell,2}_{O(n),\pm 1} =\{ (\ab,d,c)\in SO(n)^{2\ell+1}\times O(n)_\pm \mid
\pab= cdc^{-1}d\}
\end{equation}
Note that $V^{\ell,i}_{O(n),+1}=X_{\mathrm{flat}}^{\ell,i}(SO(n))$.
Recall that $X_{\mathrm{flat}}^{\ell,i}(SO(n))$ has two connected
components $X_{\mathrm{flat}}^{\ell,1}(SO(n))^{+1}$ and
$X_{\mathrm{flat}}^{\ell,2}(SO(n))^{-1}$.

For $i=1,2$, $SO(n)$ acts on $V^{\ell,i}_{O(n),\pm 1}$ by
\begin{equation}\label{eqn:OactI}
g\cdot(\ab,c)=(ga_1g^{-1},gb_1g^{-1},\ldots,ga_\ell g^{-1}, gb_\ell g^{-1}, gcg^{-1})
\end{equation}
\begin{equation}\label{eqn:OactII}
g\cdot(\ab,d,c)=(ga_1g^{-1},gb_1g^{-1},\ldots,ga_\ell g^{-1}, gb_\ell g^{-1},
gdg^{-1},gcg^{-1})
\end{equation}

When $n=2$, we have diffeomorphisms
$O(2)_+\cong O(2)_-\cong U(1)$, and
diffeomorphisms
$$
V^{\ell,i}_{O(2),+1}\cong X_{\mathrm{flat}}^{\ell,i}(U(1))
\cong U(1)^{2\ell+i-1}\times\{\pm 1\}
$$
where $i=1,2$. For any $d\in SO(2)$ and $c\in O(2)_-$, we have
$$
c^2=I_2,\quad cdc^{-1}d=I_2,
$$
so
\begin{eqnarray*}
V^{\ell,1}_{O(2),-1}&=&
\{(\ab,c)\in SO(2)^{2\ell}\times O(2)_-\mid I_2=c^2\}
\\&=& SO(2)^{2\ell}\times O(2)_-, \\
V^{\ell,2}_{O(2),-1}&=&
\{ (\ab,d,c)\in SO(2)^{2\ell+1}\times O(2)_-\mid
I_2 = cdc^{-1} d\}
\\&=& SO(2)^{2\ell+1}\times O(2)_-.
\end{eqnarray*}
For $i=1,2$, $V^{\ell,i}_{O(2), -1}$ is diffeomorphic to $U(1)^{2\ell+i}$,
thus nonempty and connected.

From now on, we assume that $n\geq 3$ so that $SO(n)$ is semisimple.
Let $\rho:Pin(n)\to O(n)$ be the double cover defined in \cite[Chapter I, Section 6]{BD},
and let $Pin(n)_\pm =\rho^{-1}(O(n)_\pm)$.
Then $Pin(n)_+$ and $Pin(n)_-$ are the two connected
components of $Pin(n)$, where  $Pin(n)_+ = Spin(n)$. Note that
$Pin(n)_-$ is not a group because if $x,y \in Pin(n)_-$ then $xy \in Pin(n)_+$.

Recall that there is an obstruction map
$$
o_2: V^{\ell,1}_{O(n),+1}=X_{\mathrm{flat}}^{\ell,1}(SO(n))\to \Ker(\rho)=\{ 1, -1\}\subset Spin(n)
$$
given by
$$
(\ab,c)\mapsto \ptab \tc^{-2}
$$
where $(\tab,\tc)$ is
the preimage of $(\ab,c)$ under $\rho^{2\ell+1}:Spin(n)^{2\ell+1}\to SO(n)^{2\ell+1}$.
It is easy to check that $o_2$ does not depend on the choice of the liftings $(\tab,\tc)$
because $2\Ker(\rho)=\{1\}$. Similarly,
 there is an obstruction map $o_2: V^{\ell,2}_{O(n),+1}=X_{\mathrm{flat}}^{\ell,2}(SO(n))
\to \{ 1, -1\}$ given by
$$
(\ab,d,c)\mapsto \ptab (\tc \td \tc^{-1}\td)^{-1}
$$
where $(\tab,\td,\tc)$ is
the preimage of $(\ab,d,c)$ under $\rho^{2\ell+2}:Spin(n)^{2\ell+2}\to SO(n)^{2\ell+2}$.
Again, $o_2$ does not depend on the choice of $\tilde{a}_i,~\tilde{b}_i,~\td,~\tc$.

For $i=1,2$, define
$V^{\ell,i,\pm 1}_{O(n),+1}=X_{\mathrm{flat}}^{\ell,i}(SO(n))^{\pm 1} =o_2^{-1}(\pm 1)$. Then
$V^{\ell,i,+1}_{O(n),+1}=X_{\mathrm{flat}}^{\ell,i}(SO(n))^{+1}$
corresponds to flat connections on the trivial $SO(n)$-bundle
($w_2=0 \in H^2(\Si^\ell_i;\bZ/2\bZ)\cong \bZ/2\bZ$),
while
$V^{\ell,i,-1}_{O(n),+1}=X_{\mathrm{flat}}^{\ell,i}(SO(n))^{-1}$
corresponds to flat connections on the nontrivial $SO(n)$-bundle
($w_2=1 \in H^2(\Si^\ell_i;\bZ/2\bZ)\cong \bZ/2\bZ$).
It was proved in \cite{HL2} that $X_{\mathrm{flat}}^{\ell,i}(SO(n))^{+1}$ and
$X_{\mathrm{flat}}^{\ell,i}(SO(n))^{-1}$ are nonempty and connected if $\ell\geq i$, i.e.,
$(\ell,i)\neq (0,1), (0,2), (1,2)$. The result is extended to the case
$(1,2)$ in \cite{HL4}.

We now extend the definition of $o_2$ to $V^{\ell,i}_{O(n),-1}$. Define
$o_2: V^{\ell,1}_{O(n),-1}\to \{ 1, -1\} \subset Spin(n)$ by
$$
(\ab,c)\mapsto \ptab \tc^{-2}
$$
where $(\tab,\tc)$ is
the preimage of $(\ab,c)$ under $\rho^{2\ell+1}:Spin(n)^{2\ell}\times Pin(n)_-\to
SO(n)^{2\ell}\times O(n)_-$.
It is easy to check that $o_2$ does not depend on the choice of $(\tab,\tc)$. Similarly,
define $o_2: V^{\ell,2}_{O(n),-1}\to \{ 1, -1\} \subset Spin(n)$ by
$$
(\ab,d,c)\mapsto \ptab (\tc \td \tc^{-1}\td)^{-1}
$$
where $(\tab,\td,\tc)$ is
the preimage of $(\ab,d,c)$ under $\rho^{2\ell+2}:Spin(n)^{2\ell+1}\times Pin(n)_-
\to SO(n)^{2\ell+1}\times O(n)_-$. Again, $o_2$ does not depend on the choice of $(\tab,\td,\tc)$.
Define $V^{\ell,i,\pm 1}_{O(n),-1}=o_2^{-1}(\pm 1)$. We will show that

\begin{pro}\label{thm:tO}
Suppose that $\ell\geq 2i$, where $i=1,2$, and $n\geq 3$. Then
$V^{\ell,i,+1}_{O(n),-1}$ and
$V^{\ell,i,-1}_{O(n),-1}$ are nonempty and connected.
\end{pro}
\begin{proof}
Define
\begin{eqnarray*}
\tV^{\ell,1,\pm 1}_{Pin(n)_-}&=&\{(\tab,\tc)\in Spin(n)^{2\ell}\times Pin(n)_-\mid \ptab\tc^{-2}=\pm 1\}\\
 \tV^{\ell,2,\pm 1}_{Pin(n)_-}&=&\{(\tab,\td,\tc)\in Spin(n)^{2\ell+1}\times Pin(n)_-
\mid \\&&\quad \ptab(\tc\td \tc^{-1}\td)^{-1} =\pm 1\}
\end{eqnarray*}
Then $\rho^{2\ell+i}: Spin(n)^{2\ell+i-1}\times Pin(n)_- \to SO(n)^{2\ell+i-1}\times O(n)_-$
restricts to a covering map $\tV^{\ell,i,\pm 1}_{Pin(n)_-}\to V^{\ell,i,\pm 1}_{O(n),-1}$.
It suffices to prove that $\tV^{\ell,i,+ 1}_{Pin(n)_-}$
and $\tV^{\ell,i,-1}_{Pin(n)_-}$ are nonempty and connected for $\ell\geq 2i$.

\paragraph{$i=1$} Define $\pi_\pm:\tV^{\ell,1,\pm 1}_{Pin(n)_-}\to Pin(n)_-$ by
$(\tab,\tc)\mapsto \tc$. Note that $Spin(n)$ is simply connected and $\tc^2, -\tc^2 \in Spin(n)$, so
$\pi_\pm^{-1}(\tc)$ is nonempty and connected for any $\tc\in Pin(n)_-$.
Let
$\ep_+=e_1e_2e_3$, and let $\ep_-=e_1$. Then $\ep_+,\ep_-\in Pin(n)_-$, and $(\ep_\pm)^2=\pm 1$.
It suffices to show that for any $\tc\in Pin(n)_-$, there is a path $\gamma_\pm:[0,1]\to
\tV^{\ell,1,\pm 1}_{Pin(n)_-}$ such that $\gamma_\pm(0)\in\pi_\pm^{-1}(\ep_\pm)$ and
$\gamma_\pm(1)\in \pi_\pm^{-1}(\tc)$.

Let $T$ be the maximal torus of $Spin(n)$, and let $\ft$ be the Lie algebra of $T$.
For any $\tc\in Pin(n)_-$, we have
$(\ep_\pm)^{-1}\tc\in Spin(n)$, so there exists $g_\pm \in Spin(n)$ such that
$(g_\pm)^{-1} (\ep_\pm)^{-1}\tc g_\pm \in T$. We have
$$
\tc=\ep_\pm g_\pm \exp(\xi_\pm) (g_\pm)^{-1}
$$
for some $\xi_+,\xi_-\in \ft$. Let $\omega$ be the coxeter element.
There are $\eta_+,\eta_-\in \ft$ such that
$$
\omega\cdot \eta_\pm -\eta_\pm=\xi_\pm.
$$
Let $a\in N(T) \subset Spin(n)$  be the corresponding
element which represents $\omega\in W=N(T)/T$. Then
$$
a \exp(t\eta_\pm) a^{-1}\exp(-t\eta_\pm)=\exp(\omega\cdot t\eta_\pm -t\eta_\pm)= \exp(t\xi_\pm).
$$

Now since $Spin(n)$ is connected, there are paths $\tg_\pm:[0,1]\to Spin(n)$ such that
$\tg_\pm(0)=1$ and $\tg_\pm(1)=g_\pm$. Now define $\gamma:[0,1]\to Spin(n)^{2\ell}\times Pin(n)_-$ by
$$
\gamma_\pm(t)=(a_1^\pm(t),b_1^\pm(t), a_2^\pm(t),b_2^\pm(t),1,\ldots,1,c^\pm(t))
$$
where
\begin{eqnarray*}
&& a_1^{\pm}(t)= \ep_\pm \tg_\pm(t) a (\ep_\pm \tg_\pm(t))^{-1},\quad
b_1^{\pm}(t)=\ep_\pm \tg_\pm(t) \exp(t\eta_\pm) (\ep_\pm \tg_\pm(t))^{-1},\\
&& a_2^{\pm}(t)= \tg_\pm(t) a (\tg_\pm(t))^{-1},\quad
b_2^{\pm}(t)=\tg_\pm(t) \exp(t\eta_\pm) (\tg_\pm(t))^{-1},\\
&& c^\pm(t)= \ep_\pm \tg_\pm(t) \exp(t\xi_\pm) (\tg_\pm(t))^{-1}
\end{eqnarray*}
Then
\begin{eqnarray*}
{[a_1^{\pm}(t),b_1^{\pm}(t)]}&= &\ep_\pm \tg_\pm(t)[a,\exp(t\eta_\pm)](\ep_\pm \tg_\pm(1))^{-1}\\
&=&\ep_\pm\tg_\pm(t) \exp(t\xi_\pm)(\ep_\pm \tg_\pm(t))^{-1}=c(t)(\ep^{-1}_\pm)=c(t)(\pm \ep_\pm),\\
{[a_2^{\pm}(t),b_2^{\pm}(t)]}&=& \tg_\pm(t)[a,\exp(t\eta_\pm)](\tg_\pm(1))^{-1}
=\tg_\pm(t)\exp(t\xi_\pm)(\tg_\pm(t))^{-1}=\ep_\pm^{-1}c(t),
\end{eqnarray*}
so the image of $\gamma_\pm$ lies in $\tV^{\ell,1,\pm}_{Pin(n)_-}$. We have
\begin{eqnarray*}
\gamma_\pm(0)&=& (\ep_\pm a \ep_\pm^{-1}, 1, a,1,1,\ldots, 1,\ep_\pm)\in \pi_\pm^{-1}(\ep_\pm)\\
\gamma_\pm(1)&=& (\ep_\pm g_\pm a (\ep_\pm g_\pm)^{-1},
\ep_\pm g_\pm \exp(\eta_\pm) (\ep_\pm g_\pm)^{-1},
 g_\pm a (g_\pm)^{-1}, g_\pm \exp(\eta_\pm) (g_\pm)^{-1}, \\
&& 1,\ldots,1,\tc)\in \pi_\pm^{-1}(\tc).
\end{eqnarray*}

\paragraph{$i=2$} Define $\pi_\pm:\tV^{\ell,2,\pm 1}_{Pin(n)_-}\to Spin(n)\times Pin(n)_-$ by
$(\tab,\td,\tc)\mapsto (\td,\tc)$. Note that $Spin(n)$ is simply connected and
$\tc\td \tc^{-1}\td,-\tc\td\tc^{-1}\td \in Spin(n)$, so
$\pi_\pm^{-1}(\td,\tc)$ is nonempty and connected for any $(\td,\tc)\in Spin(n)\times Pin(n)_-$.
Let $\ep_+=1$, and let $\ep_-=e_2e_3$. Then
$e_1 \ep_\pm e_1^{-1}\ep_\pm=e_1^{-1}\ep_\pm e_1\ep_\pm =\pm 1$.
It suffices to show that for any $(\td,\tc)\in Spin(n)\times Pin(n)_-$, there is a path
$\gamma_\pm:[0,1]\to \tV^{\ell,2,\pm 1}_{Pin(n)_-}$ such that
$\gamma_\pm(0)\in\pi_\pm^{-1}(\ep_\pm,e_1)$ and $\gamma(1)\in \pi_\pm^{-1}(\td,\tc)$.

Let $T$ be the maximal torus of $Spin(n)$, and let $\ft$ be the Lie algebra of $T$.
Given $\td\in Spin(n)$ and $\tc\in Pin(n)_-$, there exist $g_+,g_-,g \in Spin(n)$ such that
and $\xi, \xi_+,\xi_- \in \ft$ such that
$$
\tc=e_1 g \exp(\xi) g^{-1}, \td=\ep_\pm g_{\pm}\exp(\xi_\pm) (g_{\pm})^{-1}.
$$
Let $\omega$ be the coxeter element. There are $\eta,\eta_+,\eta_-\in \ft$ such that
$$
\omega\cdot \eta -\eta=\xi,\quad
\omega\cdot \eta_\pm -\eta_\pm=\xi_\pm.
$$
Let $a\in N(T) \subset Spin(n)$ be the corresponding element which represents
$\omega\in W=N(T)/T$. Then
\begin{eqnarray*}
a\exp(t\eta)a^{-1}\exp(-t\eta)&=&\exp(\omega\cdot t\eta -t\eta)= \exp(t\xi),\\
a\exp(t\eta_\pm)a^{-1}\exp(-t\eta_\pm)&=&\exp(\omega\cdot t\eta_\pm -t\eta_\pm)= \exp(t\xi_\pm).
\end{eqnarray*}

Now since $Spin(n)$ is connected, there are paths $\tg, \tg_+,
\tg_-:[0,1]\to Spin(n)$ such that
$$
\tg(0)=\tg_+(0)=\tg_-(0)=1, \quad \tg(1)=g,\quad \tg_\pm(1)=g_\pm.
$$
Now define $\gamma:[0,1]\to Spin(n)^{2\ell+1}\times Pin(n)_-$ by
$$
\gamma_\pm(t)=(a_1(t),b_1(t), a_2^\pm(t),b_2^\pm(t),
a_3^\pm(t),b_3^\pm(t), a_4^\pm(t), b_4^\pm(t),1,\ldots,1,d^\pm(t),c(t))
$$
where
\begin{eqnarray*}
&& a_1(t)= e_1 \tg(t) a (e_1\tg(t))^{-1},\quad
b_1(t)=e_1\tg(t) \exp(t\eta) (e_1\tg(t))^{-1},\\
&& a_2^{\pm}(t)= e_1\ep_\pm\tg_\pm(t) a (e_1\ep_\pm\tg_\pm(t))^{-1},\quad
b_2^{\pm}(t)=e_1\ep_\pm\tg_\pm(t)\exp(t\eta_\pm) (e_1\ep_\pm\tg_\pm(t))^{-1},\\
&& a_3^{\pm}(t)= e_1\ep_\pm\tg(t) a (e_1\ep_\pm\tg(t))^{-1},\quad
b_3^{\pm}(t)=e_1\ep_\pm \tg(t) \exp(-t\eta) (e_1\ep_\pm \tg(t))^{-1},\\
&& a_4^{\pm}(t)= \tg_\pm(t) a \tg_\pm(t)^{-1},\quad
b_4^{\pm}(t)=\tg_\pm(t)\exp(t\eta_\pm) \tg_\pm(t)^{-1},\\
&& c(t)= e_1 \tg(t) \exp(t\xi) \tg(t)^{-1},\quad
d^{\pm}(t)=\ep_\pm \tg_\pm(t) \exp(t\xi_\pm)\tg_\pm(t)^{-1}.
\end{eqnarray*}
Then
\begin{eqnarray*}
&& [a_1(t),b_1(t)]= e_1 \tg(t)[a,\exp(t\eta)](e_1\tg(t))^{-1}
=c(t)e_1^{-1},\\
&& [a_2^{\pm}(t),b_2^{\pm}(t)]= e_1\ep_\pm\tg_\pm(t)[a,\exp(t\eta_\pm)](e_1\ep_\pm\tg_\pm(t))^{-1}
=e_1d(t)(e_1\ep_\pm)^{-1}\\
 && [a_3^{\pm}(t),b_3^{\pm}(t)]= e_1\ep_\pm \tg(t)[a,\exp(-t\eta)]\tg(t)^{-1}(e_1 \ep_\pm)^{-1}
=e_1 \ep_\pm c(t)^{-1} (\pm \ep_\pm),\\
&& [a_4^{\pm}(t),b_4^{\pm}(t)]= \tg_\pm(t)[a,\exp(t\eta_\pm)](\tg_\pm(t))^{-1}=\ep_\pm^{-1} d(t),
\end{eqnarray*}
so the image of $\gamma_\pm$ lies in $\tV^{\ell,2,\pm}_{Pin(n)_-}$. We have
\begin{eqnarray*}
\gamma_\pm(0)&=& (e_1a e_1^{-1},1,~e_1\ep_\pm a (e_1\ep_\pm)^{-1},1,
~e_1\ep_\pm a (e_1\ep_\pm)^{-1},1,a,1,1,\ldots,1,\ep_\pm,e_1)\\&&\in \pi_\pm^{-1}(\ep_\pm,e_1)\\
\gamma_\pm(1)&=& (e_1 g a (e_1 g)^{-1},e_1 g \exp(\eta) (e_1 g)^{-1},
e_1\ep_\pm g_\pm a(e_1\ep_\pm g_\pm)^{-1},\\&&
e_1\ep_\pm g_\pm \exp(\eta_\pm)(e_1\ep_\pm g_\pm)^{-1},
e_1\ep_\pm g a(e_1\ep_\pm g)^{-1},
e_1\ep_\pm g \exp(-\eta)(e_1\ep_\pm g)^{-1},\\&&
g_\pm a  g_\pm^{-1},~g_\pm \exp(\eta_\pm) g_\pm^{-1},
~1,\ldots,1,\td,\tc) \in\pi_\pm^{-1}(\td,\tc)
\end{eqnarray*}

\end{proof}

\section{Yang-Mills $SO(2n+1)$-Connections}
\label{sec:SOodd}

The maximal torus of $SO(2n+1)$ consists of block diagonal matrices of the form
$$
\diag(A_1,\ldots,A_n, I_1),
$$
where $A_1,\ldots, A_n\in SO(2)$, and $I_1$ is the $1\times 1$
identity matrix. The Lie algebra of the maximal torus consists
of matrices of the form
$$
2\pi  \diag(t_1 J,\ldots,t_n J, 0I_1)=
2\pi \left(\begin{array}{ccccccc}
0    & -t_1 &       &       &     & 0    & 0 \\
t_1  &  0   &       &       &     &      & 0 \\
     &      & \cdot &       &     &      &   \\
     &      &       & \cdot &     &      &   \\
     &      &       &       & 0   & -t_n &   \\
0    &      &       &       & t_n & 0    & 0 \\
0    &  0   &       &       &     & 0    & 0
\end{array}\right),
$$
where
\begin{equation}\label{eqn:J}
J=\left(\begin{array}{cc} 0 & -1\\ 1 & 0 \end{array}\right).
\end{equation}

The fundamental Weyl chamber is
$$
\BC_0=\{\sqrt{-1}\diag(t_1 J,\ldots, t_n J,0I_1 )\mid
t_1\geq t_2\geq\cdots\geq t_n\geq 0\}.
$$

In this section, we assume
$$
n_1,\ldots,n_r\in \bZ_{>0},\quad n_1+\cdots+n_r=n.
$$

\subsection{$SO(2n+1)$-connections on orientable surfaces}
\label{sec:SOodd_orientable}

Let $J_m$ denote the $2m\times 2m$ matrix $\diag(\underbrace{J,\ldots,J}_m)$.
Any $\mu\in \BC_0$ is of the form
$$
\mu = \sqrt{-1}\diag(\lambda_1 J_{n_1},\ldots, \lambda_r J_{n_r},0I_1),
$$
where $\lambda_1>\cdots >\lambda_r \geq 0$.

Let $X_\mu=-2\pi\sqrt{-1}\mu$. Then
$$
SO(2n+1)_{X_\mu}\cong \left\{
\begin{array}{ll}
\Phi(U(n_1))\times \cdots\times \Phi(U(n_r)),& \lambda_r>0,\\
\Phi(U(n_1))\times \cdots\times \Phi(U(n_{r-1}))\times SO(2n_r+1), &\lambda_r=0,
\end{array}\right.
$$
where $\Phi:U(m)\hookrightarrow SO(2m)$ is the standard
embedding defined as follows. Consider the $\bR$-linear map
$L:\bR^{2m}\to \bC^m$ given by
$$
\left(\begin{array}{c}x_1\\y_1\\ \vdots\\x_m\\y_m \end{array}\right)\mapsto
\left(\begin{array}{c}x_1+\sqrt{-1}y_1 \\ \vdots\\ x_m+\sqrt{-1}y_m\end{array}\right).
$$
We have $L^{-1}\circ (\sqrt{-1}I_m)\circ L(v)= J_m v$ for $v\in \bR^{2m}$.
If $A$ is a $m\times m$ matrix, let
$\Phi(A)$ be the $2m\times 2m$ matrix defined by
\begin{equation}\label{eqn:PhiA}
L^{-1}\circ A\circ  L(v)= \Phi(A)(v),\quad v\in \bR^{2m}.
\end{equation}
Note that
$A (\sqrt{-1}I_m) =(\sqrt{-1}I_m) A \Rightarrow
J_m \Phi(A)=\Phi(A)J_m$.

Suppose that $(\ab,X_\mu)\in X_{\mathrm{YM}}^{\ell,0}(SO(2n+1))$. Then
$$
\exp(X_\mu)=\pab
$$
where $a_i,~b_i\in SO(2n+1)_{X_\mu}$.
This implies that $\exp(X_\mu)\in (SO(2n+1)_{X_\mu})_{ss}$, the semisimple
part of $SO(2n+1)_{X_\mu}$:
$$
(SO(2n+1)_{X_\mu})_{ss}
=\left\{\begin{array}{ll}
\Phi(SU(n_1))\times\cdots\times\Phi(SU(n_r)), & \lambda_r>0,\\
\Phi(SU(n_1))\times\cdots\times\Phi(SU(n_{r-1}))\times SO(2n_r+1), & \lambda_r=0.
\end{array}\right.
$$
Thus
\begin{eqnarray*}
X_\mu &=& 2\pi\diag\Bigl(\frac{k_1}{n_1} J_{n_1},\ldots,
\frac{k_r}{n_r} J_{n_r}, 0 I_{1}\Bigr),\\
\mu&=& \sqrt{-1}\diag\Bigl(\frac{k_1}{n_1} J_{n_1},\ldots,
\frac{k_r}{n_r} J_{n_r},0 I_{1}\Bigr),
\end{eqnarray*}
where
$$k_1,\ldots,k_r \in \bZ,
\quad\frac{k_1}{n_1}>\cdots > \frac{k_r}{n_r} \geq 0.
$$
This agrees with Section \ref{sec:SOoddC}.

Recall that for each $\mu$, the representation variety is
\[
\VymSOO{\ell}{0}_\mu=\{ (\ab)\in (SO(2n+1)_{X_\mu})^{2\ell}\mid
\pab=\exp(X_\mu) \}.
\]
For $i= 1,\cdots,\ell$, write
$$
\begin{array}{cl}
a_i=\diag(A^i_1,\ldots, A^i_r,I_1),\ b_i=\diag(B^i_1,\ldots, B^i_r,I_1), & \textup{when }k_r>0,\\
a_i=\diag(A^i_1,\ldots, A^i_r),\ b_i=\diag(B^i_1,\ldots, B^i_r), & \textup{when }k_r=0,
\end{array}
$$
where $A^i_j, ~B^i_j \in \Phi(U(n_{j}))$ for $j=1,\ldots,r-1$, and
$$
A^i_r, ~B^i_r \in \left\{\begin{array}{ll} \Phi(U(n_r)), & \textup{when }k_r>0,\\
SO(2n_r+1), & \textup{when }k_r=0.\end{array}\right.
$$
Let
$$
\hat{J}_t =\exp(2\pi t J)
=\left(\begin{array}{rr}\cos (2\pi t) & -\sin (2\pi t) \\
\sin (2\pi t) & \cos (2\pi t)\end{array}\right),
$$
and let
\begin{equation}\label{eqn:Tnk}
T_{n,k}=\Phi(e^{2\pi\sqrt{-1}k/n}I_n)=
\diag(\underbrace{\hat{J}_{k/n},\ldots,\hat{J}_{k/n} }_n) \in SO(2n).
\end{equation}
For $j=1,\ldots, r-1$, define
\begin{equation}\label{eqn:Vj}
\begin{aligned}
V_j&=\Bigl\{(A^1_j,B^1_j,\ldots, A^\ell_j,B^\ell_j)
\in \Phi(U(n_j))^{2\ell}\mid \prod_{i=1}^\ell[A^i_j,B^i_j]= T_{n_j,k_j}
  \Bigr \} \\
&\stackrel{\Phi}{\cong}
\Bigl\{(A^1_j,B^1_j,\ldots, A^\ell_j,B^\ell_j)
\in U(n_j)^{2\ell}\mid \prod_{i=1}^\ell[A^i_j,B^i_j]=
e^{2\pi\sqrt{-1}k_j/n_j}I_{n_j})
  \Bigr \}\\
&\cong X_{\mathrm{YM}}^{\ell,0}(U(n_j))_{-\frac{k_j}{n_j},\dots,-\frac{k_j}{n_j}}.
\end{aligned}
\end{equation}
If $k_r>0$, define $V_r$ by \eqref{eqn:Vj}. If $k_r=0$, define
\begin{eqnarray*}
V_r &=&\{(A^1_r,B^1_r,\ldots, A^\ell_r,B^\ell_r)
\in SO(2n_r+1)^{2\ell}\mid \prod_{i=1}^\ell[A^i_r,B^i_r]=I_{2n_r+1} \Bigr \}\\
&\cong& X^{\ell,0}_{\mathrm{flat}}(SO(2n_r+1)).
\end{eqnarray*}
Then $\VymSOO{\ell}{0}_\mu=\prod_{j=1}^r V_j$. We have
a homeomorphism
\[
\VymSOO{\ell}{0}_\mu/SO(2n+1)_{X_\mu}=
\begin{cases}
\displaystyle{\prod_{j=1}^r (V_j/U(n_j))},  &k_r>0, \\
\displaystyle{\prod_{j=1}^{r-1} (V_j/U(n_j)) \times V_r/SO(2n_r+1)},  &k_r=0,
\end{cases}
\]
and a homotopy equivalence
\[
{\VymSOO{\ell}{0}_\mu}^{h SO(2n+1)_{X_\mu} } \sim
\begin{cases}
\prod_{j=1}^r {V_j}^{h U(n_j)},  &k_r>0, \\
\prod_{j=1}^{r-1} {V_j}^{h U(n_j)} \times {V_r}^{h SO(2n_r+1)},  &k_r=0.
\end{cases}
\]

\begin{no}\label{notation}
Suppose that $m\geq 3$. Let $\Si$ be a closed, orientable
or nonorientable surface.
Let $P_{SO(m)}^{+ 1}$ and $P_{SO(m)}^{-1}$
denote the principal $SO(m)$-bundle on $\Si$
with $w_2(P_{SO(m)}^{+1})=0$ and $w_2(P_{SO(m)}^{-1})=1$ respectively in
$H^2(\Si;\bZ/2\bZ)\cong \bZ/2\bZ$. Let
$\cN(\Si)_{SO(m)}^{\pm 1}$ denote the
space of Yang-Mills connections on $P_{SO(m)}^{\pm 1}$,
and let $\cN_0(\Si)_{SO(m)}^{\pm}$ denote
the space of flat connections on $P_{SO(m)}^{\pm 1}$.

For $i=0,1,2$, we have
$$
X_{\mathrm{YM}}^{\ell,i}(SO(m))=
X_{\mathrm{YM}}^{\ell,i}(SO(m))^{+1}\cup X_{\mathrm{YM}}^{\ell,i}(SO(m))^{-1}
$$
where
$$
X_{\mathrm{YM}}^{\ell,i}(SO(m))^{\pm 1}\cong
\cN(\Si^\ell_i)_{SO(m)}^{\pm 1}/\cG_0(P_{SO(m)}^{\pm 1}),
$$
and
$$
X_{\mathrm{flat}}^{\ell,i}(SO(m))=
X_{\mathrm{flat}}^{\ell,i}(SO(m))^{+1}\cup X_{\mathrm{flat}}^{\ell,i}(SO(m))^{-1}
$$
where
$$
X_{\mathrm{flat}}^{\ell,i}(SO(m))^{\pm 1}
=\cN_0(\Si^\ell_i)_{SO(m)}^{\pm 1}/\cG_0(P_{SO(m)}^{\pm 1})
$$
is nonempty and connected for $\ell\geq 1$.
Let
$$
X_{\mathrm{YM}}^{\ell,i}(SO(m))^{\pm 1}_\mu
=X_{\mathrm{YM}}^{\ell,i}(SO(m))_\mu \cap X_{\mathrm{YM}}^{\ell,i}(SO(m))^{\pm 1}
$$
be the representation varieties for Yang-Mills connections of type $\mu$
on $P_{SO(m)}^{\pm 1}$.
Let
$$
\cM(\Si,P_{SO(m)}^{\pm 1})
=X_{\mathrm{flat}}^{\ell,i}(SO(m))^{\pm 1}/SO(m)
$$
be the moduli space of gauge equivalence classes of flat connections
on $P^{\pm 1}_{SO(m)}$ over $\Si$.
Let
$$
\cM(\Si^\ell_0, P^{n,k})=\ymU{\ell}{0}_{\frac{k}{n},\ldots,\frac{k}{n}}/U(n)
$$
be the moduli space of gauge equivalence classes of central
Yang-Mills connections on a degree $k$ principal $U(n)$-bundle
over $\Si^\ell_0$. Recall that there is no flat connection on a degree $k\neq 0$ principal
$U(n)$-bundle over $\Si^\ell_0$.
\end{no}

We have seen that $\VymSOO{\ell}{0}_\mu=\prod_{j=1}^r V_j$ is connected for $k_r>0$ and
disconnected with two connected components for $k_r=0$. To determine the underlying topological type
of the $SO(2n+1)$-bundle, let us consider the group homomorphism
$$
\phi_\mu:\pi_1(SO(2n+1)_{X_\mu})\to \pi_1(SO(2n+1))\cong \bZ/2\bZ
$$
induced by the inclusion $SO(2n+1)_{X_\mu} \hookrightarrow SO(2n+1)$. We have
$$
\pi_1(SO(2n+1)_{X_\mu})\cong\left\{ \begin{array}{ll}
\displaystyle{ \prod_{j=1}^r \pi_1(U(n_j)) } \cong \bZ^r, & \lambda_r >0,\\
\displaystyle{ \prod_{j=1}^{r-1} \pi_1(U(n_j)) } \times \pi_1(SO(2n_r+1))\cong \bZ^{r-1}\times \bZ/2\bZ, & \lambda_r =0,
\end{array}\right.
$$
and
$$
\phi_\mu(k_1,\ldots, k_r)= k_1+\cdots +k_r \quad (\mathrm{mod}\ 2).
$$
Thus, for $k_r>0$, $\VymSOO{\ell}{0}_\mu$ is from the trivial $SO(2n+1)$-bundle if and only
if $k_1+\cdots+k_r=0 \quad (\mathrm{mod}\ 2)$; and for $k_r=0$,  $\VymSOO{\ell}{0}_\mu$ has
two connected components $\VymSOO{\ell}{0}^+_\mu$ and $\VymSOO{\ell}{0}^-_\mu$, where
\begin{eqnarray*}
\VymSOO{\ell}{0}^+_\mu&=&\prod_{j=1}^{r-1}V_j \times
X^{\ell,0}_{\mathrm{flat}}(SO(2n_r+1))^{(-1)^{k_1+\cdots+k_{r-1}}},\\
\VymSOO{\ell}{0}^-_\mu&=&\prod_{j=1}^{r-1}V_j \times
X^{\ell,0}_{\mathrm{flat}}(SO(2n_r+1))^{(-1)^{k_1+\cdots+k_{r-1}+1}}.
\end{eqnarray*}

To simplify the notation, we write
$$
\mu=(\mu_1,\ldots,\mu_n)=\Bigl(\underbrace{\frac{k_1}{n_1},\ldots,\frac{k_1}{n_1} }_{n_1}, \ldots,
\underbrace{\frac{k_r}{n_r},\ldots,\frac{k_r}{n_r} }_{n_r}\Bigr)
$$
instead of
$$
\sqrt{-1}\diag\Bigl(\frac{k_1}{n_1}J_{n_1},\ldots, \frac{k_r}{n_r} J_{n_r}, 0I_1\Bigr).
$$

Let
$$
I_{SO(2n+1)}=\Bigl \{
\mu=
\Bigl(\underbrace{\frac{k_1}{n_1},\ldots,\frac{k_1}{n_1}}_{n_1}, \ldots,
\underbrace{\frac{k_r}{n_r},\ldots,\frac{k_r}{n_r}}_{n_r}\Bigr)\Bigr|
\begin{array}{c} n_j\in\bZ_{>0},\ n_1+\cdots+ n_r=n\\
k_j\in \bZ, \ \frac{k_1}{n_1}>\cdots >\frac{k_r}{n_r}\geq 0
\end{array}\Bigr\},
$$
\begin{eqnarray*}
I_{SO(2n+1)}^{\pm1}&=& \{\mu\in I_{SO(2n+1)}\mid \mu_n>0, (-1)^{k_1+\cdots+k_r}=\pm 1 \},\\
I_{SO(2n+1)}^0&=&\{ \mu\in I_{SO(2n+1)}\mid \mu_n=0\}.
\end{eqnarray*}

From the discussion above, we conclude:
\begin{pro}\label{thm:muSOodd_orientable}
Suppose that $\ell\geq 1$.
Let
\begin{equation}\label{eqn:muSOodd_orientable}
\mu=\Bigl(\underbrace{\frac{k_1}{n_1},\ldots,\frac{k_1}{n_1} }_{n_1}, \ldots,
\underbrace{\frac{k_r}{n_r},\ldots,\frac{k_r}{n_r} }_{n_r}\Bigr)\in I_{SO(2n+1)}.
\end{equation}
\begin{enumerate}
\item[(i)] If $\mu\in I_{SO(2n+1)}^{\pm 1}$, then
$\ymSOO{\ell}{0}_\mu =\ymSOO{\ell}{0}_\mu^{\pm 1}$ is nonempty and connected.
We have a homeomorphism
$$
\ymSOO{\ell}{0}_\mu /SO(2n+1) \cong
\prod_{j=1}^r \cM(\Si^\ell_0,P^{n_j,-k_j})
$$
and a homotopy equivalence
$$
{\ymSOO{\ell}{0}_\mu }^{h SO(2n+1)} \sim
\prod_{j=1}^r
\Bigl(X_{\mathrm{YM}}^{\ell,0}(U(n_j))_{-\frac{k_j}{n_j},...,-\frac{k_j}{n_j}} \Bigr)^{h U(n_j)}.
$$

\item[(ii)] If $\mu\in I_{SO(2n+1)}^0$, then $\ymSOO{\ell}{0}_\mu$ has two connected
components (from both bundles over $\Si^\ell_0$)
$$
\ymSOO{\ell}{0}_\mu^{+1} \quad \textup{and} \quad  \ymSOO{\ell}{0}_\mu^{-1}.
$$
We have a homeomorphism
$$
\ymSOO{\ell}{0}^{\pm 1}_\mu /SO(2n+1)
\cong \prod_{j=1}^{r-1} \cM(\Si^\ell_0,P^{n_j,-k_j}) \times
\cM\Bigl(\Si^\ell_0, P_{SO(2n_r+1)}^{\pm(-1)^{k_1+\cdots+k_{r-1}}}\Bigr)
$$
and a homotopy equivalence
\begin{eqnarray*}
 &&\Bigl(\ymSOO{\ell}{0}_\mu^\pm \Bigr)^{h SO(2n+1)}
 \sim\prod_{j=1}^{r-1}
\Bigl(X_{\mathrm{YM}}^{\ell,0}(U(n_j))_{-\frac{k_j}{n_j},...,-\frac{k_j}{n_j}} \Bigr)^{h U(n_j)} \times\\
&& \hspace{1in}
\Bigl(X_{\mathrm{flat}}^{\ell,0}(SO(2n_r+1))^{\pm(-1)^{k_1+\cdots+ k_{r-1}}}\Bigr)^{h SO(2n_r+1)}.
\end{eqnarray*}
\end{enumerate}
\end{pro}

\begin{pro}Suppose that $\ell\geq 1$.
The connected components of the representation variety
$\ymSOO{\ell}{0}^{\pm 1}$ are
$$
\{\ymSOO{\ell}{0}_\mu \mid \mu\in I^{\pm 1}_{SO(2n+1)}\}\cup
\{\ymSOO{\ell}{0}_\mu^{\pm 1}\mid \mu\in I^0_{SO(2n+1)}\}.
$$
\end{pro}

The following is an immediate consequence of
Proposition \ref{thm:muSOodd_orientable}.
\begin{thm}
  Suppose that $\ell\geq 1$, and let $\mu$ be as in \eqref{eqn:muSOodd_orientable}.
\begin{enumerate}
\item[(i)] If $\mu\in I_{SO(2n+1)}^{\pm 1} $, then
$$
P_t^{SO(2n+1)}\left(\ymSOO{\ell}{0}_\mu \right)
= \prod_{j=1}^r P_t^{U(n_j)}
\Bigl(X_{\mathrm{YM}}^{\ell,0}(U(n_i))_{-\frac{k_j}{n_j},\ldots,-\frac{k_j}{n_j}}\Bigr).
$$
\item[(ii)] If $\mu\in I_{SO(2n+1)}^0$, then
\begin{eqnarray*}
&&
P_t^{SO(2n+1)}\left(\ymSOO{\ell}{0}^{\pm 1}_\mu\right)=
\prod_{j=1}^{r-1} P_t^{U(n_j)}
\Bigl(X_{\mathrm{YM}}^{\ell,0}(U(n_j))_{-\frac{k_j}{n_j},\ldots,-\frac{k_j}{n_j}}\Bigr)
 \times\\&&\hspace{1in}
P_t^{SO(2n_r+1)}\left(X_{\mathrm{flat}}^{\ell,0}
(SO(2n_r+1))^{\pm (-1)^{k_1+\cdots+k_{r-1}}}\right).
\end{eqnarray*}
\end{enumerate}
\end{thm}

\subsection{Equivariant Poincar\'{e} series}
\label{sec:SOodd_poincare}

Recall from Section \ref{sec:SOoddC}:
\begin{eqnarray*}
&& \Delta=\{\alpha_i=\theta_i-\theta_{i+1}\mid i=1,\ldots, n-1\}\cup \{ \alpha_n=\theta_n\}\\
&&  \Delta^\vee =\{ \alpha_i^\vee = e_i - e_{i+1}\mid i=1,\ldots,n-1\}\cup
\{ \alpha_n^\vee =2e_n\}\\
&& \pi_1(H)=\bigoplus_{i=1}^n \bZ e_i,
\quad\Lambda=\bigoplus_{i=1}^{n-1} \bZ(e_i-e_{i+1}) \oplus \bZ(2e_n),\\
&&\pi_1(SO(2n+1))= \langle e_n \rangle \cong \bZ/2\bZ
\end{eqnarray*}

We now apply Theorem \ref{thm:Pt} to the case $G_\bR=SO(2n+1)$.
\begin{eqnarray*}
&& \varpi_{\alpha_i}= \theta_1+\cdots + \theta_i,
\quad i=1,\ldots, n-1,\quad \varpi_{\alpha_n}=\frac{1}{2}(\theta_1+\cdots+\theta_n) \\
&&
\varpi_{\alpha_i}(ke_n)=\left\{
\begin{array}{ll}
0 & i<n\\
k/2 & i=n
\end{array} \right.
\end{eqnarray*}
Case 1. $\alpha_n\in I$:
\begin{eqnarray*}
&& I= \{ \alpha_{n_1}, \alpha_{n_1+n_2},\ldots, \alpha_{n_1+\cdots + n_{r-1}},\alpha_n  \}\\
&& L^I= GL(n_1,\bC)\times \cdots \times  GL(n_r,\bC), \quad n_1+\cdots +n_r=n \\
&& \dim_\bC \fz_{L^I}-\dim_\bC \fz_{SO(2n+1,\bC)}=r,\quad
\dim_\bC U^I = \sum_{1\leq i<j\leq r}n_i n_j +\frac{n(n+1)}{2},\\
&&\rho^I =\frac{1}{2}\sum_{i=1}^r
\biggl(n-2\sum_{j=1}^{i} n_j +n_i \biggr)
\biggl(\sum_{j=1}^{n_i}\theta_{n_1+\cdots+ n_{i-1}+j}\biggr)
+\frac{n}{2}(\theta_1+\cdots +\theta_n)
\end{eqnarray*}
$$
\langle \rho^I, \alpha_{n_1+\cdots+n_i}^\vee \rangle = \frac{n_i+n_{i+1}}{2}
\textup{ for } i=1,\ldots,r-1,
\quad \langle \rho^I, \alpha_n^\vee \rangle = n_r
$$
Case 2. $\alpha_n\notin I$:
\begin{eqnarray*}
&& I= \{ \alpha_{n_1}, \alpha_{n_1+n_2},\ldots, \alpha_{n_1+\cdots + n_{r-1}}  \}\\
&& L^I= GL(n_1,\bC)\times \cdots \times  GL(n_{r-1},\bC)\times SO(2n_r+1,\bC),
\quad n_1+\cdots +n_r=n \\
&& \dim_\bC \fz_{L^I}-\dim_\bC \fz_{SO(2n+1,\bC)}=r-1,\quad\\
&&\dim_\bC U^I = \sum_{1\leq i<j\leq r}n_i n_j
+\frac{n(n+1)-n_r(n_r+1)}{2},\\
&& \rho^I =\frac{1}{2}\sum_{i=1}^r
\biggl(n-2\sum_{j=1}^{i} n_j +n_i \biggr)
\biggl(\sum_{j=1}^{n_i}\theta_{n_1+\cdots + n_{i-1}+j}\biggr) \\
&& \quad\quad +
\frac{n}{2}(\theta_1+\cdots +\theta_{n_1+\cdots+n_{r-1}})
+ \frac{n-n_r}{2}(\theta_{n_1+\cdots+n_{r-1}+1}+\cdots +\theta_n)
\end{eqnarray*}
$$
\langle \rho^I,  \alpha_{n_1+\cdots+n_i}^\vee \rangle = \frac{n_i+n_{i+1}}{2}
\textup{ for }i=1,\ldots, r-2,\quad
\langle \rho^I, \alpha_{n_1+\cdots+n_{r-1}}^\vee \rangle =\frac{n_{r-1}}{2}+ n_r
$$

We have the following closed formula for the $SO(2n+1)$-equivariant poincar\'{e} series
of the representation of flat $SO(2n+1)$-connections:

\begin{thm}\label{thm:PtSOodd}
\begin{eqnarray*}
&& P_t^{SO(2n+1)}(X_{\mathrm{flat}}^{\ell,0}(SO(2n+1))^{(-1)^k} ) \\
&=&\sum_{r=1}^n\sum_{\tiny \begin{array}{c}n_1,\ldots, n_r\in\bZ_{>0}\\\sum n_j=n \end{array}}
\left( (-1)^r \prod_{i=1}^r \frac{\prod_{j=1}^{n_i} (1+t^{2j-1})^{2\ell} }
{ (1-t^{2n_i})\prod_{j=1}^{n_i-1}(1-t^{2j})^2 }\right. \\
&& \cdot \frac{t^{(\ell-1)(2\sum_{i<j} n_i n_j+n(n+1))} }
{\left[\prod_{i=1}^{r-1}(1-t^{2(n_i+n_{i+1} )} )\right](1-t^{4n_r}) }
\cdot t^{2\sum_{i=1}^{r-1}(n_i +n_{i+1})+4 n_r\langle k/2\rangle}\\
&& + (-1)^{r-1} \prod_{i=1}^{r-1} \frac{\prod_{j=1}^{n_i} (1+t^{2j-1})^{2\ell} }
{ (1-t^{2n_i})\prod_{j=1}^{n_i-1}(1-t^{2j})^2 }
\cdot\frac{\prod_{j=1}^{n_r}(1+t^{4j-1})^{2\ell}}{\prod_{j=1}^{2n_r}(1-t^{2j})}\\
&& \left.\cdot
 \frac{t^{(\ell-1)(2\sum_{i<j} n_i n_j+n(n+1)-n_r(n_r+1))} }
{\left[\prod_{i=1}^{r-2}(1-t^{2(n_i+n_{i+1})} ) \right](1-\epsilon(r)t^{2n_{r-1}+4n_r}) }
t^{2\sum_{i=1}^{r-1}(n_i+ n_{i+1})+ 2\epsilon(r)n_r }
\right)
\end{eqnarray*}
where
$$
\epsilon(r)=\left\{\begin{array}{ll}0 & r=1\\ 1 & r>1 \end{array}\right.
$$
\end{thm}

\begin{rem}
We have
$$
 P_t^{SO(2n+1)}(X_{\mathrm{flat}}^{\ell,0}(SO(2n+1))^{+1} )
= P_t^{Spin(2n+1)}(X_{\mathrm{flat}}^{\ell,0}(Spin(2n+1)) ),
$$
so Theorem \ref{thm:PtSOodd} also gives a formula
for $X_{\mathrm{flat}}^{\ell,0}(Spin(2n+1))$.
\end{rem}

\begin{ex}\label{SOthree}
\begin{eqnarray*}
&& P^{SO(3)}_t(X_{\mathrm{flat}}^{\ell,0}(SO(3))^{+1})
=P^{Spin(3)}_t(X_{\mathrm{flat}}^{\ell,0}(Spin(3)))\\
&=& -\frac{(1+t)^{2\ell} t^{2\ell+2} }{(1-t^2)(1-t^4)} +\frac{(1+t^3)^{2\ell}}{(1-t^2)(1-t^4)}\\
&& P^{SO(3)}_t(X_{\mathrm{flat}}^{\ell,0}(SO(3))^{-1})\\
&=& -\frac{(1+t)^{2\ell} t^{2\ell} }{(1-t^2)(1-t^4)} +\frac{(1+t^3)^{2\ell}}{(1-t^2)(1-t^4)}
\end{eqnarray*}

Note that $Spin(3)=SU(2)$, so
$$
P^{Spin(3)}_t(X_{\mathrm{flat}}^{\ell,0}(Spin(3)))
=P^{SU(2)}_t(X_{\mathrm{flat}}^{\ell,0}(SU(2)))
$$
as expected, where $P^{SU(2)}_t(X_{\mathrm{flat}}^{\ell,0}(SU(2))$
is calculated in Example \ref{UtwoSUtwo}.
\end{ex}

\begin{ex}\label{SOfive}
\begin{eqnarray*}
&& P^{SO(5)}_t(X_{\mathrm{flat}}^{\ell,0}(SO(5))^{+1})
=P^{Spin(5)}_t(X_{\mathrm{flat}}^{\ell,0}(Spin(5)))\\
&=& -\frac{(1+t)^{2\ell}(1+t^3)^{2\ell} t^{6\ell+2} }
{(1-t^2)^2(1-t^4)(1-t^8)} +\frac{(1+t^3)^{2\ell}(1+t^7)^{2\ell}}{(1-t^2)(1-t^4)(1-t^6)(1-t^8)}\\
&& +\frac{(1+t)^{4\ell} t^{8\ell}}{(1-t^2)^2(1-t^4)^2}
 - \frac{(1+t)^{2\ell}(1+t^3)^{2\ell} t^{6\ell}}{(1-t^2)^2 (1-t^4)(1-t^6)}\\
&& P^{SO(5)}_t(X_{\mathrm{flat}}^{\ell,0}(SO(5))^{-1})\\
&=&  -\frac{(1+t)^{2\ell}(1+t^3)^{2\ell} t^{6\ell-2} }
{(1-t^2)^2(1-t^4)(1-t^8)} +\frac{(1+t^3)^{2\ell}(1+t^7)^{2\ell}}{(1-t^2)(1-t^4)(1-t^6)(1-t^8)}\\
&& +\frac{(1+t)^{4\ell} t^{8\ell-2}}{(1-t^2)^2(1-t^4)^2}
 - \frac{(1+t)^{2\ell}(1+t^3)^{2\ell} t^{6\ell}}{(1-t^2)^2 (1-t^4)(1-t^6)}
\end{eqnarray*}
\end{ex}

\subsection{$SO(2n+1)$-connections on nonorientable surfaces}
\label{sec:SOodd_nonorientable}
We have $\BC_0^\tau=\BC_0$ (any $\mu\in \BC_0$ is conjugate to $-\mu$).
Any $\mu\in \BC_0^\tau$ is of the form
$$
\mu=\sqrt{-1}\diag(\lambda_1 J_{n_1},\ldots, \lambda_r J_{n_r},0I_1)
$$
where $\lambda_1>\cdots >\lambda_r \geq 0$. We have
$$
SO(2n+1)_{X_\mu}\cong \left\{\begin{array}{ll}
\Phi(U(n_1))\times \cdots\times\Phi(U(n_r)), & \lambda_r>0,\\
\Phi(U(n_1))\times\cdots\times\Phi(U(n_{r-1}))\times SO(2n_r+1), & \lambda_r=0,
\end{array}\right.
$$
where $X_\mu=-2\pi\sqrt{-1}\mu$, and $\Phi:U(m)\hookrightarrow SO(2m)$ is the standard embedding.

Given $\mu\in \BC_0$, define
$$
\ep_\mu=\left\{
\begin{array}{ll}
\diag\left(H_n, (-1)^n I_1\right),& \lambda_r>0,\\
\diag\left(H_{n-n_r},(-1)^{n-n_r} I_1, I_{2n_r}\right), & \lambda_r=0.
\end{array}\right.
$$
Then $\Ad(\ep_\mu) X_\mu=-X_\mu$. Suppose that
$$
(\ab,\ep_\mu c',X_\mu/2)\in X_{\mathrm{YM}}^{\ell,1}(SO(2n+1)).
$$
Then
\[
\exp(X_\mu/2)\ep_\mu c' \ep_\mu c'=\pab
\]
where
$$
a_i,~b_i,~c'\in\left\{
\begin{array}{ll}
\Phi(U(n_1))\times\cdots \times  \Phi(U(n_r)), & \lambda_r>0, \\
\Phi(U(n_1))\times\cdots \times \Phi(U(n_{r-1}))\times SO(2n_r+1), & \lambda_r=0.
\end{array}\right.
$$

We first assume that $\lambda_r>0$.
Let $L:\bR^{2n}\to \bC^n$ defined as in Section \ref{sec:SOodd_orientable}. Define
\begin{eqnarray*}
X'_\mu &=& L \circ 2\pi\diag(\lambda_1J_{n_1},\ldots,\lambda_rJ_{n_r}) \circ L^{-1}\\
&=& 2\pi\sqrt{-1}\diag(\lambda_1 I_{n_1},\ldots,\lambda_r I_{n_r})\in
\mathfrak{u}(n_1)\times\cdots\times\mathfrak{u}(n_r).
\end{eqnarray*}

We have $L\circ H_n \circ L^{-1}(v)=\bar{v}$ for $v\in\bC^n$,
where $\bar{v}$ is the complex conjugate of $v$. So
\begin{eqnarray*}
 L\circ H_n \Phi(c') H_n \circ L^{-1}(v)
&=&(L\circ H_n \circ L^{-1}) (L\circ \Phi(c')\circ L^{-1})(L\circ H_n \circ L^{-1})(v)\\
&=& (L\circ H_n \circ L^{-1})c'\bar{v}=\overline{c'\bar{v}} = \bar{c'}v.
\end{eqnarray*}
So the condition on $X'_\mu$ is
$$
\exp(X'_\mu/2)\bar{c'}c'=\pab\in SU(n_1)\times\cdots\times SU(n_r),
$$
where $a_i,~b_i,~c'\in U(n_1)\times\cdots \times  U(n_r)$, and $\bar{c'}$ is the complex conjugate
of $c'$. In order that this is nonempty, we need
$1=\det(e^{\pi\sqrt{-1}\lambda_j}I_{n_j})$, or equivalently
\begin{equation}\label{eqn:twok}
\lambda_j=\frac{2k_j}{n_j}\quad  k_j,~n_j \in \bZ_{>0}
\end{equation}
for $j=1,\ldots,r$.

When $\lambda_r=0$, the above argument gives the condition \eqref{eqn:twok}
for $j=1,\ldots,r-1$.

Similarly, suppose that
$$
(\ab,d, \ep_\mu c',X_\mu/2)\in X_{\mathrm{YM}}^{\ell,1}(SO(2n+1)).
$$
Then
\[
\exp(X_\mu/2)(\ep_\mu c') d(\ep_\mu c')^{-1}d =\pab,
\]
where
$$
a_i,~b_i,~d,~c'\in\left\{
\begin{array}{ll}
\Phi(U(n_1))\times\cdots \times  \Phi(U(n_r)), & \lambda_r>0, \\
\Phi(U(n_1))\times\cdots \times \Phi(U(n_{r-1}))\times SO(2n_r+1), & \lambda_r=0.
\end{array}\right.
$$
When $\lambda_r>0$, the condition on $X'_\mu$ is
$$
\exp(X'_\mu/2)\bar{c'}\bar{d}\bar{c'}^{-1}d =\pab\in SU(n_1)\times \cdots \times SU(n_r).
$$
Again, we need $1=\det(e^{\pi\sqrt{-1}\lambda_j}I_{n_j})$, or equivalently \eqref{eqn:twok}.
When $\lambda_r=0$ we get condition \eqref{eqn:twok} for $j=1,\ldots,r-1$.

We conclude that for nonorientable surfaces,
$$
\mu= \sqrt{-1}\diag\Bigl(\frac{2k_1}{n_1} J_{n_1},\ldots, \frac{2k_r}{n_r} J_{n_r},0 I_1\Bigr),
~\textup{where }
k_1,\ldots,k_r \in \bZ, ~\frac{k_1}{n_1}>\cdots> \frac{k_r}{n_r}\geq 0.
$$

Recall that for {\em orientable} surfaces we have
$$
\mu=\sqrt{-1}\diag\Bigl(\frac{k_1}{n_1} J_{n_1},\ldots, \frac{k_r}{n_r} J_{n_r},0 I_1\Bigr),
~\textup{where }k_1,\ldots,k_r\in \bZ,~\frac{k_1}{n_1}>\cdots >\frac{k_r}{n_r}\geq 0.
$$

For each $\mu$, define $\ep_\mu$-reduced representation varieties
\begin{eqnarray*}
\VymSOO{\ell}{1}_\mu&=&
\{ (\ab,c')\in (SO(2n+1)_{X_\mu})^{2\ell+1}\mid \\
&&\quad
\pab=\exp(X_\mu/2) \epsilon_\mu c' \epsilon_\mu c'\},\\
\VymSOO{\ell}{2}_\mu&=&
\{ (\ab,d,c')\in (SO(2n+1)_{X_\mu})^{2\ell+2}\mid \\
&&\quad
\pab=\exp(X_\mu/2)\epsilon_\mu c' d (\epsilon_\mu c')^{-1} d \}.
\end{eqnarray*}

For $i=1,\ldots,\ell$, write
\begin{eqnarray*}
 && a_i=\diag(A^i_1,\ldots, A^i_r, I_1), \quad b_i=\diag(B^i_1,\ldots, B^i_r,I_1),\\
&& c' =\diag(C_1,\ldots,C_r,I_1),\quad  d  =\diag(D_1,\ldots,D_r,I_1),
\end{eqnarray*}
when $k_r>0$, and write
\begin{eqnarray*}
&& a_i=\diag(A^i_1,\ldots, A^i_r),\quad b_i=\diag(B^i_1,\ldots, B^i_r),\\
&& c' =\diag(C_1,\ldots,C_r),\quad  d  =\diag(D_1,\ldots,D_r),
\end{eqnarray*}
when $k_r=0$, where $A^i_j, ~B^i_j, ~D_j, ~C_j \in \Phi(U(n_{j}))$ for $j=1,\ldots,r-1$, and
$$
A^i_r, ~B^i_r, ~D_r, ~C_r\in
\left\{
\begin{array}{ll}
\Phi(U(n_r)) & \textup{when }k_r>0,\\
SO(2n_r+1) &\textup{when }k_r=0.
\end{array}
\right.
$$

$i=1$. Let $T_{n,k}$ be defined as in \eqref{eqn:Tnk}, and let $\ep_j=\diag(H_{n_j})$.
For $j=1,\ldots,r-1$, define
\begin{equation}\label{eqn:VjRP}
\begin{aligned}
V_j&=
\Bigl\{(A^1_j,B^1_j,\ldots, A^\ell_j,B^\ell_j,C_j)
\in \Phi(U(n_j))^{2\ell+1}\mid \prod_{i=1}^\ell[A^i_j,B^i_j]= T_{n_j,k_j}\ep_j C_j \ep_j C_j
  \Bigr \} \\
&\stackrel{\Phi}{\cong}
\Bigl\{(A^1_j,B^1_j,\ldots, A^\ell_j,B^\ell_j,C_j)
\in U(n_j)^{2\ell+1}\mid \prod_{i=1}^\ell[A^i_j,B^i_j]=
e^{2\pi\sqrt{-1}k_j/n_j} \bar{C_j} C_j
  \Bigr \}\\
&\cong \tV^{\ell,1}_{n_j,-k_j}
\end{aligned}
\end{equation}
where $\tV^{\ell,1}_{n_j,-k_j}$ is the twisted representation
variety defined in \eqref{eqn:twistI} of Section \ref{sec:twisted}.
$\tV^{\ell,1}_{n_j,-k_j}$ is nonempty if $\ell\geq 1$. We have
shown that $\tV^{\ell,1}_{n_j,-k_j}$ is connected if $\ell\geq 2$
(Proposition \ref{thm:tV}).

When $k_r>0$, define $V_r$ by \eqref{eqn:VjRP}. When $k_r=0$, define
\begin{equation}
V_r=
\Bigl\{(A^1_r,B^1_r,\ldots, A^\ell_r,B^\ell_r,C_r)
\in SO(2n_r+1)^{2\ell+1}\mid \prod_{i=1}^\ell[A^i_r,B^i_r]= (\epsilon C_r)^2
  \Bigr \},
\end{equation} where $\ep=\diag((-1)^{n-n_r}I_1, I_{2n_r}),\det(\ep)=(-1)^{n-n_r}$. Let
$C_r'=\epsilon C_r$. We see that
\begin{eqnarray*}
V_r&\cong&
\Bigl\{(A^1_r,B^1_r,\ldots, A^\ell_r,B^\ell_r,C_r')
\in SO(2n_r+1)^{2\ell}\times O(2n_r+1)\mid \\
&&\prod_{i=1}^\ell[A^i_r,B^i_r]= (C_r')^2,
\det(C'_r)=(-1)^{n-n_r}  \Bigr \}\\
&\cong & V^{\ell,1}_{O(2n_r+1),(-1)^{n-n_r}}
\end{eqnarray*}
where
$V^{\ell,1}_{O(n),\pm 1}$ is the twisted representation variety defined in (\ref{eqn:Oone})
of Section \ref{sec:twistedO}. $V^{\ell,1}_{O(n),\pm 1}$ is nonempty if $\ell\geq 2$. We have shown
that $V^{\ell,1}_{O(n),\pm 1}$ is disconnected with two components
$V^{\ell,1,+ 1}_{O(n),\pm 1}$ and $V^{\ell,1,-1}_{O(n),\pm 1}$
if $\ell\geq 2$ and
$n>2$ (Proposition \ref{thm:tO}).

We have
$$
\VymSOO{\ell}{1}_\mu=\prod_{j=1}^r V_j.
$$
We define a $U(n_j)$-action on $V_j=\tV^{\ell,1}_{n_j,-k_j}$ by \eqref{eqn:actI}
of Section \ref{sec:twisted}; when $k_r=0$, we define an $SO(2n_r+1)$-action on
$V_r=V^{\ell,1}_{O(2n_r+1),(-1)^{n-n_r}}$ by \eqref{eqn:OactI} of Section \ref{sec:twistedO}
Then we have a homeomorphism
$$
\VymSOO{\ell}{1}_\mu/SO(2n+1)_{X_\mu}\cong
\begin{cases}
\displaystyle{ \prod_{j=1}^r (V_j/U(n_j)) }, & k_r >0,\\
\displaystyle{ \prod_{j=1}^{r-1}(V_j/U(n_j)) } \times V_r/SO(2n_r+1), & k_r=0,
\end{cases}
$$
and a homotopy equivalence
$$
{\VymSOO{\ell}{1}_\mu}^{h SO(2n+1)_{X_\mu}}\sim
\begin{cases}
\prod_{j=1}^r {V_j}^{h U(n_j)}, & k_r >0,\\
\prod_{j=1}^{r-1}{V_j}^{ hU(n_j)} \times {V_r}^{h SO(2n_r+1)}. & k_r=0,
\end{cases}
$$

$i=2$.  Let $\ep_j=\diag(H_{n_j})$. Define
\begin{equation}\label{eqn:VjKlein}
\begin{aligned}
V_j&=
\Bigl\{(A^1_j,B^1_j,\ldots, A^\ell_j,B^\ell_j,D_j,C_j)
\in \Phi(U(n_j))^{2\ell+2}\mid \\&\quad\quad\quad
 \prod_{i=1}^\ell[A^i_j,B^i_j]=
T_{n_j,k_j}\ep_j C_j D_j (\ep_j C_j)^{-1} D_j
  \Bigr \} \\
&=
\Bigl\{(A^1_j,B^1_j,\ldots, A^\ell_j,B^\ell_j,D_j,C_j)
\in \Phi(U(n_j))^{2\ell+2}\mid \\
&\quad\quad\quad \prod_{i=1}^\ell[A^i_j,B^i_j]=
T_{n_j,k_j}\ep_j C_j \ep_j^{-1} \ep_j D_j \epsilon^{-1}
\ep_j C_j^{-1} \ep_j^{-1} D_j
  \Bigr \} \\
&\stackrel{\Phi}{\cong}
\Bigl\{(A^1_j,B^1_j,\ldots, A^\ell_j,B^\ell_j,D_j,C_j)
\in U(n_j)^{2\ell+2}\mid \\&\quad\quad\quad
\prod_{i=1}^\ell[A^i_j,B^i_j]=
e^{2\pi\sqrt{-1}k_j/n_j} \bar{C_j} \bar{D_j}\bar{C}^{-1}_jD_j
  \Bigr \} \cong \tV^{\ell,2}_{n_j,-k_j}
\end{aligned}
\end{equation}
where $\tV^{\ell,2}_{n_j,-k_j}$ is the twisted representation
variety defined in \eqref{eqn:twistII} of Section \ref{sec:twisted}.
$\tV^{\ell,2}_{n_j,-k_j}$ is nonempty if $\ell\geq 1$. We have
shown that $\tV^{\ell,2}_{n_j,-k_j}$ is connected if $\ell\geq 4$
(Proposition \ref{thm:tV}).

When $k_r>0$, define $V_r$ by \eqref{eqn:VjKlein}. When $k_r=0$, define
\begin{equation}\begin{aligned}
V_r=&
\Bigl\{(A^1_r,B^1_r,\ldots, A^\ell_r,B^\ell_r,D_r,C_r)
\in SO(2n_r+1)^{2\ell+2}\mid \\&\quad
\prod_{i=1}^\ell[A^i_r,B^i_r]=\epsilon C_r D_r
(\epsilon C_r)^{-1}D_r
  \Bigr \},
\end{aligned}\end{equation}
where $\ep=\diag((-1)^{n-n_r}I_1, I_{2n_r}),\det(\ep)=(-1)^{n-n_r}$. Let
$C_r'=\epsilon C_r$. We see that
\begin{eqnarray*}
V_r&\cong&
\Bigl\{(A^1_r,B^1_r,\ldots, A^\ell_r,B^\ell_r,D_r,C_r')
\in SO(2n_r+1)^{2\ell+1}\times O(2n_r+1)\mid \\
&&\prod_{i=1}^\ell[A^i_r,B^i_r]=C_r' D_r
C_r'^{-1}D_r,\det(C'_r)=(-1)^{n-n_r} \Bigr \}\\
&\cong &V^{\ell,2}_{O(2n_r+1),(-1)^{n-n_r}}
\end{eqnarray*}
where
$V^{\ell,2}_{O(n),\pm 1}$ is the twisted representation variety defined in (\ref{eqn:Otwo})
of Section \ref{sec:twistedO}. $V^{\ell,2}_{O(n),\pm 1}$ is nonempty if $\ell\geq 4$. We have shown
that $V^{\ell,2}_{O(n),\pm 1}$ is disconnected with two components
$V^{\ell,2,+1}_{O(n),\pm 1}$ and $V^{\ell,2,-1}_{O(n),\pm 1}$ if $\ell\geq 4$ and
$n>2$ (Proposition \ref{thm:tO}).

We have
$$
\VymSOO{\ell}{2}_\mu=\prod_{j=1}^r V_j.
$$
We define a $U(n_j)$-action on $V_j=\tV^{\ell,2}_{n_j,-k_j}$ by \eqref{eqn:actII}
of Section \ref{sec:twisted}; when $k_r=0$, we define an $SO(2n_r+1)$-action on
$V_r=V^{\ell,2}_{O(2n_r+1),(-1)^{n-n_r}}$ by \eqref{eqn:OactII}
of Section \ref{sec:twistedO}. Then we have a homeomorphism
$$
\VymSOO{\ell}{2}_\mu/SO(2n+1)_{X_\mu}\cong
\begin{cases}
\displaystyle{ \prod_{j=1}^r (V_j/U(n_j)) }, & k_r>0,\\
\displaystyle{ \prod_{j=1}^{r-1}(V_j/U(n_j)) }\times V_r/SO(2n_r+1), & k_r=0,
\end{cases}
$$
and a homotopy equivalence
$$
{\VymSOO{\ell}{2}_\mu}^{h SO(2n+1)_{X_\mu} }\sim
\begin{cases}
\prod_{j=1}^r {V_j}^{h U(n_j)}, & k_r>0,\\
\prod_{j=1}^{r-1} {V_j}^{h U(n_j)}\times {V_r}^{h SO(2n_r+1)}, & k_r=0,
\end{cases}
$$

We have seen that for $i=1,2$, $\VymSOO{\ell}{i}_\mu$ is connected when
$k_r>0$. In this case, to determine the topological type of the underlying
$SO(2n+1)$-bundle $P$ over $\Si^{\ell}_{i}$,
we can just look at a special point in $\VymSOO{\ell}{i}_\mu$
where $c',d$ are the identity element $I_{2n+1}$. Then
$$
\pab =\exp(X_\mu/2),
$$
so $\ab$ can be viewed as the holonomies of a Yang-Mills connection on an $SO(2n+1)$-bundle
$Q_0\to \Si^\ell_0$. Also, $c=\ep=\diag(H_{n},(-1)^{n}I_1)$ can be viewed as the holonomy of a flat
connection on an $SO(2n+1)$-bundle $Q_1$ over $\Si^0_1= \RP^2$, and
$c=\ep$, $d=I_{2n+1}$ can be viewed as the holonomies of
a flat connection on an $SO(2n+1)$-bundle $Q_2$ over $\Si^0_2$ (a Klein bottle).
Let $\Si'$ be obtained by gluing $\Si^\ell_0$ and $\Si^0_i$ at a point, and
let $P'\to \Si'$ be the (topological) principal $SO(2n+1)$-bundle over $\Si'$ such that
$P'|_{\Si^{\ell}_0}=Q_0$ and $P'|_{\Si^0_i}=Q_i$.
Then $P=p^*P'$ where $p:\Si^\ell_i\rightarrow \Si'=\Si^\ell_0\cup\Si^0_i$ is the collapsing map.
Then $w_2(P')= (w_2(Q_0), w_2(Q_i))$ under the isomorphism
$$
H^2(\Si';\bZ/2\bZ)\cong H^2(\Si^\ell_0;\bZ/2\bZ)\oplus H^2(\Si^0_i;\bZ/2\bZ),
$$
and $w_2(P)=p^* w_2(P') = w_2(Q_0) + w_2(Q_i)$, if we identify
$H^2(\Si^\ell_i;\bZ/2\bZ)$, $H^2(\Si^\ell_0;\bZ/2\bZ)$, and
$H^2(\Si^0_i;\bZ/2\bZ)$ with $\bZ/2\bZ$. So it remains to compute
$w_2(Q_0)$, $w_2(Q_1)$, and $w_2(Q_2)$.
We have $Q_0\cong P^{n,-(k_1+\cdots+k_r)}\times_{U(n)} SO(2n+1)$, so
$w_2(Q_0)= k_1+\cdots+k_r$ (mod 2).
To compute $w_2(Q_1)$ and $w_2(Q_2)$, we lift $c=\ep$ to $\tilde{c}\in Spin(2n+1)$
and lift $d=I_{2n+1}$ to $\tilde{d}\in Spin(2n+1)$.
Since $2\pi_1(SO(2n+1))$ is the trivial group, we may choose any lifting for $c$ and $d$.
We choose $\tilde{d}=1\in Spin(2n+1)$ and
$$
\tilde{c}=\left\{\begin{array}{ll}
e_2 e_4 \cdots e_{2n}, & n\textup{ even},\\
e_2 e_4 \cdots e_{2n} e_{2n+1}, & n\textup{ odd}.
\end{array} \right.
$$
Then $\tilde{c}^2 = (-1)^{n(n+1)/2}$ and $\tilde{c}\tilde{d}\tilde{c}^{-1}\tilde{d} =1$.
We conclude that
$$
w_2(Q_1)= \frac{n(n+1)}{2}\ (\textup{mod }2),\quad
w_2(Q_2)=0\ (\textup{mod }2),
$$
so
$$
w_2(P)= k_1+\cdots +k_r +i\frac{n(n+1)}{2}\ (\textup{mod }2).
$$

When $k_r=0$, we have seen that $\VymSOO{\ell}{i}_\mu$ is disconnected with two components.
To determine the corresponding underlying topological types, we consider two special cases.

\paragraph{\em Case 1.} We consider special points
$$
(\ab,c)\in \VymSOO{\ell}{1}_\mu, ~(\ab,d,c)\in \VymSOO{\ell}{2}_\mu,
$$
where
\begin{eqnarray*}
&& a_i=\diag(A^i_1,\ldots, A^i_{r-1},I_{2n_r+1}),
\quad b_i=\diag(B^i_1,\ldots, B^i_{r-1},I_{2n_r+1}),\\
&&c=\ep_\mu=\diag(H_{n-n_r},(-1)^{n-n_r}I_1, I_{2n_r}),\quad
d=I_{2n+1}.
\end{eqnarray*}
Let $\ep_1=\diag((-1)^{n-n_r}I_1, I_{2n_r})$.
Then
\begin{eqnarray*}
&& (A^i_j,B^i_j,\ldots,A^i_j, B^i_j)\in X^{\ell,0}_{\mathrm{YM}}
(U(n_j))_{-\frac{k_j}{n_j},\ldots,-\frac{k_j}{n_j}},\quad j=1,\ldots, r-1,\\
&& (I_{2n_r+1},\ldots, I_{2n_r+1},\ep_1) \in V^{\ell,1,(-1)^{n-n_r}}_{O(2n_r+1), (-1)^{n-n_r}},\\
&&(I_{2n_r+1},\ldots,I_{2n_r+1},I_{2n_r+1},\ep_1) \in V^{\ell,2,1}_{O(2n_r+1),(-1)^{n-n_r}}.
\end{eqnarray*}
We have $P=P_1\times P_2$, where $P_1$ is an $SO(2(n-n_r)+1)$-bundle, and
$P_2$ is an $SO(2n_r)$-bundle with trivial holonomies $I_{2n_r}$. We have
\begin{eqnarray*}
w_2(P)=w_2(P_1)= k_1+\cdots+k_{r-1} + i\frac{(n-n_r)(n-n_r+1)}{2}.
\end{eqnarray*}

\paragraph{\em Case 2.} We consider special points
$$
(\ab,c)\in \VymSOO{\ell}{1}_\mu, ~(\ab,d,c)\in \VymSOO{\ell}{2}_\mu,
$$
where
\begin{eqnarray*}
&& a_i=\diag(A^i_1,\ldots, A^i_{r-1},I_{2n_r+1}),
\quad b_i=\diag(B^i_1,\ldots, B^i_{r-1},I_{2n_r+1}),\\
&&c=\diag(H_{n-n_r},(-1)^{(n-n_r)}I_1,-I_2, I_{2n_r-2}),
\quad d=\diag(I_{2(n-n_r)+1}, -I_2, I_{2n_r-2}).
\end{eqnarray*}
Let $\ep_1=\diag((-1)^{n-n_r}I_1,-I_2, I_{2n_r-2})$, $\ep_2=\diag(I_1,-I_2, I_{2n_r-2})$. Then
\begin{eqnarray*}
&& (A^i_j,B^i_j,\ldots,A^i_j, B^i_j)\in X^{\ell,0}_{\mathrm{YM}}
(U(n_j))_{-\frac{k_j}{n_j},\ldots,-\frac{k_j}{n_j}},\quad j=1,\ldots, r-1,\\
&& (I_{2n_r+1},\ldots,I_{2n_r+1},\ep_1) \in V^{\ell,1,-(-1)^{n-n_r}}_{O(2n_r+1), (-1)^{n-n_r}},\\
&& (I_{2n_r+1},\ldots,I_{2n_r+1},\ep_2,\ep_1) \in V^{\ell,2,-1}_{O(2n_r+1),(-1)^{n-n_r}}.
\end{eqnarray*}
We have $P=P_1\times P_2$, where $P_1$ is an $SO(2(n-n_r)+1)$-bundle, and
$P_2$ is an $SO(2n_r)$-bundle with holonomies $a_i=b_i=I_{2n_r}$, $c=d=\ep=\diag(-I_2, I_{2n_r-2})$.
Similarly, we can choose the lifting of $d$ and $c$ as $\tilde{d}=\tilde{c}=e_1e_2$. Then
$\tilde{c}^2=\tilde{c}\tilde{d}\tilde{c}^{-1}\tilde{d}=-1$.
We have
\begin{eqnarray*}
w_2(P_1)= k_1+\cdots+k_{r-1} + i\frac{(n-n_r)(n-n_r+1)}{2}\ (\mathrm{mod}\ 2),\quad
w_2(P_2)=1\ (\mathrm{mod}\ 2),
\end{eqnarray*}
so
$$
w_2(P)=w_2(P_1)+w_2(P_2)= k_1+\cdots+k_{r-1} + i\frac{(n-n_r)(n-n_r+1)}{2} + 1\
(\mathrm{mod}\ 2).
$$

To summarize, when $k_r=0$ we have
$$
\VymSOO{\ell}{i}_\mu^\pm  =\prod_{j=1}^{r-1} V_j \times
V_{O(2n_r+1),(-1)^{n-n_r}  }^{\ell, i,
\pm (-1)^{k_1+\cdots + k_{r-1} + i\frac{(n-n_r)(n-n_r-1)}{2}}  },
$$
where $\VymSOO{\ell}{i}_\mu^\pm$ is the $\ep_\mu$-reduced version of
$X_{\mathrm{YM}}^{\ell,i}(SO(2n+1))^{\pm 1}_\mu$.


To simplify the notation, we write
$$
\mu=(\mu_1,\ldots,\mu_n)=\Bigl(\underbrace{\frac{2k_1}{n_1},\ldots,\frac{2k_1}{n_1} }_{n_1}, \ldots,
\underbrace{\frac{2k_r}{n_r},\ldots,\frac{2k_r}{n_r} }_{n_r}\Bigr)
$$
instead of
$$
\sqrt{-1}\diag\Bigl(\frac{2k_1}{n_1}J_{n_1},\ldots, \frac{2k_r}{n_r} J_{n_r}, 0I_1\Bigr).
$$
Let
\begin{eqnarray*}
\hat{I}_{SO(2n+1)}&=&\Bigl \{
\mu=\Bigl(\underbrace{\frac{2k_1}{n_1},\ldots,\frac{2k_1}{n_1}}_{n_1}, \ldots,
\underbrace{\frac{2k_r}{n_r},\ldots,\frac{2k_r}{n_r}}_{n_r}
\Bigr) \bigl|\ n_j\in \bZ_{> 0},\\&&\quad n_1+\cdots+n_r=n,\ k_j\in\bZ,\
\frac{k_1}{n_1}>\cdots >\frac{k_r}{n_r}\geq 0\bigr. \Bigr\},\\
\hat{I}_{SO(2n+1)}^{\pm 1}&=& \{\mu\in \hat{I}_{SO(2n+1)}\mid \mu_n>0,
(-1)^{k_1+\cdots+k_r +\frac{i n(n+1)}{2}}=\pm 1 \},\\
\hat{I}_{SO(2n+1)}^0&=& \{\mu\in \hat{I}_{SO(2n+1)}\mid \mu_n=0\}.
\end{eqnarray*}
For $i=1,2$, define twisted moduli spaces
$$
\tilde{\cM}^{\ell,i}_{n,k}=\tV^{\ell,i}_{n,k}/U(n),\quad
\cM^{\ell,i,\pm 1}_{O(n),\pm 1} = V^{\ell, i, \pm 1}_{O(n),\pm
1}/SO(n).
$$

\begin{pro}\label{thm:muSOodd_nonorientable}
Suppose that $\ell\geq  2i$, where $i=1,2$.
Let
\begin{equation}\label{eqn:muSOodd_nonorientable}
\mu=\Bigl(\underbrace{\frac{2k_1}{n_1},\ldots,\frac{2k_1}{n_1} }_{n_1}, \ldots,
\underbrace{\frac{2k_r}{n_r},\ldots,\frac{2k_r}{n_r} }_{n_r}\Bigr)\in \hat{I}_{SO(2n+1)}.
\end{equation}
\begin{enumerate}
\item[(i)] If $\mu\in \hat{I}_{SO(2n+1)}^{\pm 1}$, then
$\ymSOO{\ell}{i}_\mu =\ymSOO{\ell}{i}_\mu^{\pm 1}$  is nonempty and
connected (coming from either the trivial bundle or the nontrivial
bundle). We have a homeomorphism
$$
\ymSOO{\ell}{i}_\mu /SO(2n+1) \cong
\prod_{j=1}^r \tilde{\cM}^{\ell,i}_{n_j,-k_j}
$$
and a homotopy equivalence
$$
{\ymSOO{\ell}{i}_\mu }^{h SO(2n+1)} \sim
\prod_{j=1}^r \bigl(\tV^{\ell,i}_{n_j,-k_j}\bigr)^{h U(n_j)}.
$$

\item[(ii)] If $\mu\in \hat{I}_{SO(2n+1)}^0$, then $\ymSOO{\ell}{i}_\mu$ has two connected
components (coming from both bundles)
$$
\ymSOO{\ell}{i}_\mu^{+1} \quad\textup{and}\quad \ymSOO{\ell}{i}_\mu^{-1}.
$$
We have homeomorphisms
\begin{eqnarray*}
 & &\ymSOO{\ell}{i}^{\pm 1}_\mu /SO(2n+1)\\
&\cong&\prod_{j=1}^{r-1} \tilde{\cM}^{\ell,i}_{n_j,-k_j}
\times \cM_{O(2n_r+1),(-1)^{n-n_r}}^{\ell,i,
\pm (-1)^{ k_1+\cdots+k_{r-1} +i \frac{(n-n_r)(n-n_r-1)}{2} } }
\end{eqnarray*}
and homotopy equivalences
\begin{eqnarray*}
&    & \Bigl(\ymSOO{\ell}{i}^{\pm 1}_\mu\Bigr)^{h SO(2n+1)} \\
&\sim& \prod_{j=1}^{r-1} \bigl(\tV^{\ell,i}_{n_j,-k_j} \bigr)^{h U(n_j)}
\times \Bigl(V_{O(2n_r+1),(-1)^{n-n_r}}^{\ell,i,
\pm (-1)^{ k_1+\cdots+k_{r-1} +i \frac{(n-n_r)(n-n_r-1)}{2} } }\Bigr)^{h SO(2n_r+1)}.
\end{eqnarray*}
\end{enumerate}
\end{pro}

\begin{pro}Suppose that $\ell\geq 2i$, where $i=0,1$.
The connected components of $\ymSOO{\ell}{i}^{\pm 1}$ are
$$
\{\ymSOO{\ell}{i}_\mu \mid \mu\in \hat{I}^{\pm 1}_{SO(2n+1)}\}
\cup \{\ymSOO{\ell}{i}_\mu^{\pm 1}\mid \mu\in \hat{I}^0_{SO(2n+1)}\}.
$$
\end{pro}

Notice that, the set
$\{\mu=\sqrt{-1}\diag(\mu_1 J,\ldots,\mu_n J,0 I_1)\mid
(\mu_1,\ldots,\mu_n)\in\hat{I}_{SO(2n+1)}\}
$ is a {\em proper} subset of $\{\mu\in(\Xi^I_+)^\tau|I\subseteq\Delta,\tau(I)=I\}$
as mentioned in Section \ref{sec:connected}.

The following is an immediate consequence of
Proposition \ref{thm:muSOodd_nonorientable}.
\begin{thm}\label{thm:PtmuSOodd_nonorientable}
  Suppose that $\ell\geq 2i$, where $i=1,2$,
and let $\mu$ be as in \eqref{eqn:muSOodd_nonorientable}.
\begin{enumerate}
\item[(i)] If $\mu\in \hat{I}_{SO(2n+1)}^{\pm 1} $, then
$$
P_t^{SO(2n+1)}\left(\ymSOO{\ell}{0}_\mu \right)
= \prod_{j=1}^r P_t^{U(n_j)}(\tV^{\ell,i}_{n_j,-k_j}).
$$
\item[(ii)] If $\mu\in \hat{I}_{SO(2n+1)}^0$, then
\begin{eqnarray*}
&& P_t^{SO(2n+1)}\left(\ymSOO{\ell}{i}^{\pm 1}_\mu\right)\\
&=&\prod_{j=1}^{r-1} P_t^{U(n_j)}(\tV^{\ell,i}_{n_j,-k_j})
\cdot P_t^{SO(2n_r+1)}\left(V^{\ell,i,\pm (-1)^{k_1+\cdots+k_{r-1}+i\frac{(n-n_r)(n-n_r-1)}{2}}}
_{O(2n_r+1), (-1)^{n-n_r} }\right).
\end{eqnarray*}
\end{enumerate}
\end{thm}

\section{Yang-Mills $SO(2n)$-Connections}
\label{sec:SOeven}

The maximal torus of $SO(2n)$ consists of block diagonal matrices
of the form
$$
\diag(A_1,\ldots,A_n)
$$
where $A_1,\ldots, A_n\in SO(2)$. The Lie algebra
of the maximal torus consists of matrices of the form
$$
2\pi\diag(t_1  J,\ldots,t_n J)
$$
where
$$
J=\left(\begin{array}{cc} 0 &-1 \\ 1&0\end{array}\right).
$$
The fundamental Weyl chamber is
$$
\BC_0=\{\sqrt{-1}\diag(t_1 J,\ldots,t_n J )\mid t_1\geq t_2\geq\cdots
\geq\mid t_n\mid\geq 0\}.
$$

As in Section \ref{sec:SOodd}, in this section we continue to assume
$$
n_1,\ldots,n_r\in\bZ_{>0},\quad n_1+\cdots+n_r=n.
$$

\subsection{$SO(2n)$-connections on orientable surfaces}
\label{sec:SOeven_orientable}

There are four cases.

\paragraph{Case 1. $t_{n-1}>|t_n|$, $n_r=1$}
$$
\mu= \sqrt{-1}\diag(\lambda_1 J_{n_1},\ldots, \lambda_{r-1} J_{n_{r-1}}, \lambda_r J),
$$
where $\lambda_1 >\cdots >\lambda_{r-1} > |\lambda_r|\geq 0$.
Thus
$$
SO(2n)_{X_\mu}\cong \Phi(U(n_1))\times\cdots\times
\Phi(U(n_{r-1}))\times \Phi(U(n_r)).
$$

Suppose that $(\ab,X_\mu)\in X_{\mathrm{YM}}^{\ell,0}(SO(2n))$. Then
$$
\exp(X_\mu)=\pab
$$
where $\ab\in SO(2n)_{X_\mu}$.  Then we have
$$
\exp(X_\mu)\in (SO(2n)_{X_\mu})_{ss}=\Phi(SU(n_1))\times\cdots\times\Phi(SU(n_{r-1}))\times\{I_2\}.
$$
Thus
\begin{eqnarray*}
X_\mu &=& 2\pi \diag\Bigl(\frac{k_1}{n_1} J_{n_1},\ldots, \frac{k_{r-1}}{n_{r-1}}
J_{n_{r-1}},k_r  J \Bigr),\\
\mu&=& \sqrt{-1}\diag\Bigl(\frac{k_1}{n_1} J_{n_1},\ldots, \frac{k_{r-1}}{n_{r-1}} J_{n_{r-1}},
k_r  J \Bigr),
\end{eqnarray*}
where
$$
k_1,\ldots,k_r\in \bZ,\quad
\frac{k_1}{n_1}>\cdots >\frac{k_{r-1}}{n_{r-1}}> |k_r|\geq 0.
$$

Recall that for each $\mu$, the representation variety is
$$
\VymSOE{\ell}{0}_\mu = \{ (\ab)\in (SO(2n)_{X_\mu})^{2\ell}\mid
\pab=\exp(X_\mu) \}.
$$

For $i= 1,\ldots,\ell$, write
$$
a_i=\diag(A^i_1,\ldots, A^i_r),\quad
b_i=\diag(B^i_1,\ldots, B^i_r),
$$
where $A^i_j, ~B^i_j \in \Phi(U(n_{j}))$.
Define $V_j$ as in \eqref{eqn:Vj}. Then
\begin{equation}\label{eqn:orient-SOeven}
\VymSOE{\ell}{0}_\mu =
\prod_{j=1}^r V_j.
\end{equation}
We have a homeomorphism
\begin{equation}\label{eqn:orient-SOeven-homeo}
\VymSOE{\ell}{0}_\mu/SO(2n)_{X_\mu}\cong \prod_{j=1}^r (V_j/U(n_j))
\end{equation}
and a homotopy equivalence
\begin{equation}\label{eqn:orient-SOeven-homotopy}
{\VymSOE{\ell}{0}_\mu}^{h SO(2n)_{X_\mu} }\sim \prod_{j=1}^r {V_j}^{h U(n_j)}.
\end{equation}

\paragraph{Case 2. $t_{n-1} = -t_n>0$, $n_r>1$}
$$
\mu =\sqrt{-1}\diag(\lambda_1 J_{n_1},\ldots,\lambda_{r-1}J_{n_{r-1}},
\lambda_{r} J_{n_r-1}, -\lambda_r J),
$$
where $\lambda_1 >\cdots >\lambda_r >0$. Thus
$$
SO(2n)_{X_\mu}\cong \Phi(U(n_1))\times\cdots\times
\Phi(U(n_{r-1}) )\times  \Phi'(U(n_r)),
$$
where
$\Phi':U(m)\hookrightarrow SO(2m)$ is the embedding
defined as follows.
Consider the $\mathbb{R}$-linear map $L':\mathbb{R}^{2m}\to \mathbb{C}^{m}$
given by
$$
\left(\begin{array}{c}x_1\\ y_1\\ \vdots\\ x_{m-1}\\y_{m-1}\\x_m\\ y_m\end{array}\right)
\mapsto \left(\begin{array}{c} x_1+\sqrt{-1}y_1 \\ \vdots \\ x_{m-1}+\sqrt{-1} y_{m-1}\\
x_m -\sqrt{-1}y_m
\end{array}\right).
$$
We have $(L')^{-1}\circ \sqrt{-1}I_m\circ L'(v)=
\diag(J_{m-1}, -J)(v)$ for $v\in \bR^{2m}$.
If $A$ is a $m\times m$ matrix, and let
$\Phi'(A)$ be the $(2m)\times (2m)$ matrix given by
 $$
(L')^{-1}\circ A \circ L'(v)=\Phi'(A)(v).
$$
Note that
$A(\sqrt{-1}I_{m}) =(\sqrt{-1}I_{m})A \Rightarrow \Phi'(A)
\diag(J_{m-1},-J)
=\diag(J_{m-1},-J)\Phi'(A)$.

Suppose that $(\ab,X_\mu)\in X_{\mathrm{YM}}^{\ell,0}(SO(2n))$. Then
$$
\exp(X_\mu)=\pab
$$
where $\ab\in SO(2n)_{X_\mu}$.  Then we have
$$
\exp(X_\mu)\in (SO(2n)_{X_\mu})_{ss}=\Phi(SU(n_1))\times\cdots\times
\Phi(SU(n_{r-1}))\times\Phi'(SU(n_r)).
$$
Thus
\begin{eqnarray*}
X_\mu &=& 2\pi\diag
\Bigl(\frac{k_1}{n_1} J_{n_1},\ldots, \frac{k_{r-1}}{n_{r-1}}J_{n_{r-1}},
\frac{k_r}{n_r} J_{n_r-1},\frac{-k_r}{n_r}J \Bigr),\\
\mu&=& \sqrt{-1}\diag
\Bigl(\frac{k_1}{n_1} J_{n_1},\ldots, \frac{k_{r-1}}{n_{r-1}}J_{n_{r-1}},
\frac{k_r}{n_r} J_{n_r-1},\frac{-k_r}{n_r}J \Bigr),
\end{eqnarray*}
where
$$
k_1,\ldots,k_r \in \bZ,\quad
\frac{k_1}{n_1}>\cdots >\frac{k_r}{n_r}>0.
$$
Recall that for each $\mu$, the representation variety is
$$
\VymSOE{\ell}{0}_\mu =\{ (\ab)\in (SO(2n)_{X_\mu})^{2\ell}\mid
\pab=\exp(X_\mu) \}.
$$

For $i= 1,\ldots,\ell$, write
$$
a_i=\diag(A^i_1,\ldots, A^i_r),\quad
b_i=\diag(B^i_1,\ldots, B^i_r),
$$
where $A^i_j, ~B^i_j \in \Phi(U(n_{j}))$ for $j=1,\ldots,r-1$, and
$A^i_r, ~B^i_r \in \Phi'(U(n_r))$.

For $j=1,\ldots,r-1$, define $V_j$ as in \eqref{eqn:Vj}.
Recall that
$$
\hat{J}_t=\left(\begin{array}{rr}
\cos(2\pi t) & -\sin(2\pi t)\\
\sin(2\pi t) & \cos(2\pi t) \end{array}\right).
$$
Define
\begin{eqnarray*}
V_r&=&
\Bigl\{(A^1_r,B^1_r,\ldots, A^\ell_r,B^\ell_r)
\in \Phi'(U(n_r))^{2\ell}\mid \\&&
\prod_{i=1}^\ell[A^i_r,B^i_r]=
\diag(\hat{J}_{k_r/n_r},\ldots,\hat{J}_{k_r/n_r},
\hat{J}_{-k_r/n_r})
  \Bigr \}\\&\stackrel{\Phi'}{\cong}&
\Bigl\{(A^1_r,B^1_r,\ldots, A^\ell_r,B^\ell_r)
\in U(n_r)^{2\ell}\mid
\prod_{i=1}^\ell[A^i_r,B^i_r]=
e^{2\pi\sqrt{-1}k_r/n_r}I_{n_r}
  \Bigr \} \\
&\cong&  X_{\mathrm{YM}}^{\ell,0}(U(n_r))_{-\frac{k_r}{n_r},\dots,-\frac{k_r}{n_r}}.
\end{eqnarray*}
Then we have \eqref{eqn:orient-SOeven}, \eqref{eqn:orient-SOeven-homeo},
and \eqref{eqn:orient-SOeven-homotopy}.

\paragraph{Case 3. $t_{n-1}=t_n> 0$, $n_r>1$}
$$
\mu= \sqrt{-1}\diag(\lambda_1 J_{n_1},\ldots, \lambda_{r} J_{n_r}),
$$
where $\lambda_1 >\cdots >\lambda_r >0$. Thus
$$
SO(2n)_{X_\mu}\cong \Phi(U(n_1))\times\cdots\times
\Phi(U(n_r)).
$$

Let $X_\mu=-2\pi\sqrt{-1}\mu$.
Suppose that $(\ab,X_\mu)\in X_{\mathrm{YM}}^{\ell,0}(SO(2n))$. Then
$$
\exp(X_\mu)=\pab
$$
where $\ab\in SO(2n)_{X_\mu}$.  Then we have
$$
\exp(X_\mu)\in (SO(2n)_{X_\mu})_{ss}
=\Phi(SU(n_1))\times\cdots\times\Phi(SU(n_r)).
$$
Thus
\begin{eqnarray*}
X_\mu &=& 2\pi\diag
\Bigl(\frac{k_1}{n_1} J_{n_1},\ldots, \frac{k_r}{n_r} J_{n_r}\Bigr),\\
\mu &=& \sqrt{-1}\diag\Bigl(\frac{k_1}{n_1} J_{n_1},\ldots, \frac{k_r}{n_r} J_{n_r} \Bigr),
\end{eqnarray*}
where
$$
k_1,\ldots,k_r\in\bZ,\quad
\frac{k_1}{n_1}>\cdots >\frac{k_r}{n_r}>0.
$$

Recall that for each $\mu$, the representation variety is
$$
\VymSOE{\ell}{0}_\mu = \{ (\ab)\in (SO(2n)_{X_\mu})^{2\ell}\mid
\pab=\exp(X_\mu) \}.
$$
For $i= 1,\ldots,\ell$, write
$$
a_i=\diag(A^i_1,\ldots, A^i_r),\quad
b_i=\diag(B^i_1,\ldots, B^i_r),
$$
where $A^i_j, ~B^i_j \in \Phi(U(n_{j}))$.

Define $V_j$ as in \eqref{eqn:Vj}.
Then we have \eqref{eqn:orient-SOeven}, \eqref{eqn:orient-SOeven-homeo},
and \eqref{eqn:orient-SOeven-homotopy}.

\paragraph{Case 4. $t_{n-1}=t_n= 0$, $n_r>1$}
$$
\mu= \sqrt{-1}\diag(\lambda_1 J_{n_1},\ldots, \lambda_{r-1} J_{n_{r-1}}, 0J_{n_r}),
$$
where $\lambda_1 >\cdots >\lambda_{r-1} >0$. Thus
$$
SO(2n)_{X_\mu}\cong \Phi(U(n_1))\times\cdots\times
\Phi(U(n_{r-1}))\times SO(2n_r).
$$

Let $X_\mu=-2\pi\sqrt{-1}\mu$.
Suppose that $(\ab,X_\mu)\in X_{\mathrm{YM}}^{\ell,0}(SO(2n))$. Then
$$
\exp(X_\mu)=\pab
$$
where $\ab\in SO(2n)_{X_\mu}$.  Then we have
$$
\exp(X_\mu)\in (SO(2n)_{X_\mu})_{ss}
= \Phi(SU(n_1))\times\cdots\times\Phi(SU(n_{r-1}))\times SO(2n_r).
$$
Thus
\begin{eqnarray*}
X_\mu &=& 2\pi\diag
\Bigl(\frac{k_1}{n_1} J_{n_1},\ldots, \frac{k_{r-1}}{n_{r-1}} J_{n_{r-1}},
0 J_{n_r}\Bigr),\\
\mu &=& \sqrt{-1}\diag\Bigl(\frac{k_1}{n_1} J_{n_1},\ldots, \frac{k_{r-1}}{n_{r-1}} J_{n_{r-1}},
0 J_{n_r}\Bigr),
\end{eqnarray*}
where
$$
k_1,\ldots,k_{r-1}\in\bZ,\quad
\frac{k_1}{n_1}>\cdots >\frac{k_{r-1}}{n_{r-1}}>0.
$$

Recall that for each $\mu$, the representation variety is
$$
\VymSOE{\ell}{0}_\mu = \{ (\ab)\in (SO(2n)_{X_\mu})^{2\ell}\mid
\pab=\exp(X_\mu) \}.
$$
For $i= 1,\ldots,\ell$, write
$$
a_i=\diag(A^i_1,\ldots, A^i_r),\quad
b_i=\diag(B^i_1,\ldots, B^i_r),
$$
where $A^i_j, ~B^i_j \in \Phi(U(n_{j}))$ for $j=1,\ldots,r-1$,
and $A^i_r, ~B^i_r \in SO(2n_r)$.

For $j=1,\ldots,r-1$, define $V_j$ as in \eqref{eqn:Vj}. Define
$$
V_r=
\Bigl\{(A^1_r,B^1_r,\ldots, A^\ell_r,B^\ell_r)
\in SO(2n_r)^{2\ell}\mid \prod_{i=1}^\ell[A^i_r,B^i_r]= I_{2n_r}
  \Bigr \}\cong X_{\mathrm{flat}}^{\ell,0}(SO(2n_r)).
$$
Then $\VymSOE{\ell}{0}_\mu =\prod_{j=1}^r V_j$. We have a
homeomorphism
$$
\VymSOE{\ell}{0}_\mu/SO(2n)_{X_\mu}\cong
\prod_{j=1}^{r-1} (V_j/U(n_j)) \times V_r/SO(2n_r)
$$
and a homotopy equivalence
$$
{\VymSOE{\ell}{0}_\mu}^{h SO(2n)_{X_\mu}}\sim
\prod_{j=1}^{r-1} {V_j}^{h U(n_j)} \times {V_r}^{h SO(2n_r)}.
$$

We can decide the topological type of the underlying $SO(2n)$ as in
Section \ref{sec:SOodd_orientable}. Then
Case 1, Case 2, Case 3 and Case 4 give exactly the same Atiyah-Bott points
as in Section \ref{sec:SOevenC}.

To simplify the notation, we write
$$
\mu=(\mu_1,\ldots,\mu_n)=\Bigl(\underbrace{\frac{k_1}{n_1},\ldots,\frac{k_1}{n_1} }_{n_1}, \ldots,
\underbrace{\frac{k_r}{n_r},\ldots,\frac{k_r}{n_r} }_{n_r}\Bigr)
$$
instead of
$$
\sqrt{-1}\diag\Bigl(\frac{k_1}{n_1}J_{n_1},\ldots, \frac{k_r}{n_r} J_{n_r}\Bigr),
$$
and write
$$
\mu=(\mu_1,\ldots,\mu_n)=\Bigl(\underbrace{\frac{k_1}{n_1},\ldots,\frac{k_1}{n_1} }_{n_1}, \ldots,
\underbrace{\frac{k_{r-1}}{n_{r-1}},\ldots,\frac{k_{r-1}}{n_{r-1}} }_{n_{r-1}},
\underbrace{\frac{k_r}{n_r},\ldots,\frac{k_r}{n_r} }_{n_r-1},\frac{-k_r}{n_r}\Bigr)
$$
instead of
$$
\sqrt{-1}\diag\Bigl(\frac{k_1}{n_1}J_{n_1},\ldots, \frac{k_{r-1}}{n_{r-1}} J_{n_{r-1}},
\frac{k_r}{n_r} J_{n_r-1},-\frac{k_r}{n_r}J\Bigr).
$$

Let
\begin{eqnarray*}
&&I_{SO(2n)}^{\pm 1}=\Bigl \{
\mu=\Bigl(\underbrace{\frac{k_1}{n_1},\ldots,\frac{k_1}{n_1}}_{n_1}, \ldots,
\underbrace{\frac{k_{r-1}}{n_{r-1}},\ldots,\frac{k_{r-1}}{n_{r-1}}}_{n_{r-1}}, k_r\Bigr)
\Bigl|\ n_j\in \bZ_{>0},\ k_j\in\bZ \\
&&\quad~
n_1+\cdots+n_{r-1}+1=n,\ \frac{k_1}{n_1}>\cdots >\frac{k_{r-1}}{n_{r-1}}>|k_r|\geq 0,\
 (-1)^{k_1+\cdots+k_r}=\pm 1 \Bigr\}\\
&&\bigcup \Bigl \{
\mu=\Bigl(\underbrace{\frac{k_1}{n_1},\ldots,\frac{k_1}{n_1}}_{n_1}, \ldots,
\underbrace{\frac{k_{r-1}}{n_{r-1}},\ldots,\frac{k_{r-1}}{n_{r-1}}}_{n_{r-1}},
\underbrace{\frac{k_r}{n_r},\ldots,\frac{k_r}{n_r}}_{n_r-1},\pm \frac{k_r}{n_r} \Bigr)
\Bigl|\ n_j \in \bZ_{>0}, \\
&& \quad~ k_j\in\bZ,
n_r\in\bZ_{>1}, n_1+\cdots+n_r=n,
\frac{k_1}{n_1}>\cdots >\frac{k_r}{n_r}>0, (-1)^{k_1+\cdots+k_r}=\pm 1 \Bigr\}
\end{eqnarray*}
\begin{eqnarray*}
&&I_{SO(2n)}^{0}=\Bigl \{
\mu=\Bigl(\underbrace{\frac{k_1}{n_1},\ldots,\frac{k_1}{n_1}}_{n_1}, \ldots,
\underbrace{\frac{k_{r-1}}{n_{r-1}},\ldots,\frac{k_{r-1}}{n_{r-1}}}_{n_{r-1}},
\underbrace{0,\ldots,0}_{n_r}\Bigr)\Bigl|\ n_j\in \bZ_{>0},\\
&&\quad\quad n_r\in\bZ_{>1},\ n_1+\cdots+n_r=n,\  k_j\in\bZ,\ \frac{k_1}{n_1}>\cdots >\frac{k_{r-1}}{n_{r-1}}>0  \Bigr\}.
\end{eqnarray*}

From the above discussion, we conclude that
\begin{pro}\label{thm:muSOeven_orientable}
Suppose that $\ell\geq 1$.
\begin{enumerate}
\item[(i)] If
$\displaystyle{
\mu=\Bigl(\underbrace{\frac{k_1}{n_1},\ldots,\frac{k_1}{n_1}}_{n_1}, \ldots,
\underbrace{\frac{k_{r-1}}{n_{r-1}},\ldots,\frac{k_{r-1}}{n_{r-1}}}_{n_{r-1}},k_r \Bigr)
\in I_{SO(2n)}^{\pm 1}
}$, or \\
$
\displaystyle{
\mu=\Bigl(\underbrace{\frac{k_1}{n_1},\ldots,\frac{k_1}{n_1}}_{n_1}, \ldots,
\underbrace{\frac{k_{r-1}}{n_{r-1}},\ldots, \frac{k_{r-1}}{n_{r-1}} }_{n_{r-1}},
\underbrace{\frac{k_r}{n_r},\ldots,\frac{k_r}{n_r}}_{n_r-1},\pm \frac{k_r}{n_r} \Bigr)
\in I_{SO(2n)}^{\pm 1}
}$,\\
then $\ymSOE{\ell}{0}_\mu =\ymSOO{\ell}{0}_\mu^{\pm 1}$  is nonempty and connected.
We have a homeomorphism
$$
\ymSOE{\ell}{0}_\mu /SO(2n) \cong
\prod_{j=1}^r \cM(\Si^\ell_0,P^{n_j,-k_j})
$$
and a homotopy equivalence
$$
{\ymSOE{\ell}{0}_\mu}^{h SO(2n)} \sim
\prod_{j=1}^r \Bigl(X_{\mathrm{YM}}^{\ell.0}(U(n_j))
_{-\frac{k_j}{n_j},\ldots,-\frac{k_j}{n_j}} \Bigr)^{h U(n_j)}.
$$

\item[(ii)] If
$\displaystyle{
\mu=\Bigl(\underbrace{\frac{k_1}{n_1},\ldots,\frac{k_1}{n_1} }_{n_1}, \ldots,
\underbrace{\frac{k_{r-1}}{n_{r-1}},\ldots,\frac{k_{r-1} }{n_{r-1}} }_{n_{r-1}},
\underbrace{0,\ldots,0}_{n_r}\Bigr)\in I_{SO(2n)}^0
}$,\\
then $\ymSOE{\ell}{0}_\mu$ has two connected components (from both bundles over $\Si^\ell_0$)
$$
\ymSOE{\ell}{0}_\mu^{+1}\quad\textup{and}\quad \ymSOE{\ell}{0}_\mu^{-1}.
$$
We have a homeomorphism
$$
\ymSOE{\ell}{0}^{\pm 1}_\mu /SO(2n)
\cong \prod_{j=1}^{r-1} \cM(\Si^\ell_0,P^{n_j,-k_j})
\times \cM\Bigl(\Si^\ell_0, P_{SO(2n_r)}^{\pm(-1)^{k_1+\cdots+k_{r-1}}}\Bigr)
$$
and a homotopy equivalence
\begin{eqnarray*}
&& {\ymSOE{\ell}{0}^{\pm 1}_\mu}^{h SO(2n)}\sim
\prod_{j=1}^{r-1} \Bigl(X_{\mathrm{YM}}^{\ell,0}(U(n_j))
     _{-\frac{k_j}{n_j},\ldots,-\frac{k_j}{n_j} }\Bigr)^{h U(n_j)}\times\\
&& \hspace{1in} \Bigl(X_{\mathrm{flat}}(SO(2n_r))^{\pm(-1)^{k_1+\cdots+k_{r-1}}}\Bigr)
^{h SO(2n_r)}.
\end{eqnarray*}
\end{enumerate}
\end{pro}

\begin{pro}
Suppose that $\ell\geq 1$. The connected components of the representation variety
$X^{\ell,0}_{\mathrm{YM}}(SO(2n))^{\pm 1}$
are
$$
\{X^{\ell,0}_{\mathrm{YM}}(SO(2n))_\mu\mid \mu\in I^{\pm 1}_{SO(2n)} \}\cup
\{X^{\ell,0}_{\mathrm{YM}}(SO(2n))_\mu^{\pm 1} \mid \mu\in I^0_{SO(2n)} \}.
$$
\end{pro}

The following is an immediate consequence of Proposition
\ref{thm:muSOeven_orientable}.
\begin{thm}
Suppose that $\ell\geq 1$.
\begin{enumerate}
\item[(i)] If
$\displaystyle{
\mu=\Bigl(\underbrace{\frac{k_1}{n_1},\ldots,\frac{k_1}{n_1}}_{n_1}, \ldots,
\underbrace{\frac{k_{r-1}}{n_{r-1}},\ldots,\frac{k_{r-1}}{n_{r-1}}}_{n_{r-1}},k_r \Bigr)
\in I_{SO(2n)}^{\pm 1}
}$, or\\
$\displaystyle{
\mu=\Bigl(\underbrace{\frac{k_1}{n_1},\ldots,\frac{k_1}{n_1}}_{n_1}, \ldots,
\underbrace{\frac{k_{r-1}}{n_{r-1}},\ldots,\frac{k_{r-1}}{n_{r-1}} }_{n_{r-1}},
\underbrace{\frac{k_r}{n_r},\ldots,\frac{k_r}{n_r}}_{n_r-1},\pm \frac{k_r}{n_r} \Bigr)
\in I_{SO(2n)}^{\pm 1}
}$, then
$$
P_t^{SO(2n)}\left(\ymSOE{\ell}{0}_\mu \right)
=\prod_{j=1}^r P_t^{U(n_j)}
\Bigl(X_{\mathrm{YM}}^{\ell,0}(U(n_i))_{-\frac{k_j}{n_j},\ldots,-\frac{k_j}{n_j}}\Bigr).
$$
\item[(ii)] If
$\displaystyle{
\mu=\Bigl(\underbrace{\frac{k_1}{n_1},\ldots,\frac{k_1}{n_1} }_{n_1}, \ldots,
\underbrace{\frac{k_{r-1}}{n_{r-1}},\ldots,\frac{k_{r-1}}{n_{r-1}} }_{n_{r-1}},
\underbrace{0,\ldots,0}_{n_r}\Bigr)\in I_{SO(2n)}^0
}$, then
\begin{eqnarray*}
&& P_t^{SO(2n)}\left(\ymSOE{\ell}{0}^{\pm 1}_\mu\right)\\
&=& \prod_{j=1}^{r-1} P_t^{U(n_j)}
\Bigl(X_{\mathrm{YM}}^{\ell,0}(U(n_j))_{-\frac{k_j}{n_j},\ldots,-\frac{k_j}{n_j}}\Bigr) \cdot
P_t^{SO(2n_r)}\left(X_{\mathrm{flat}}^{\ell,0}
(SO(2n_r))^{\pm  (-1)^{k_1+\cdots+k_{r-1}}}\right).
\end{eqnarray*}
\end{enumerate}
\end{thm}

\subsection{Equivariant Poincar\'{e} series}
\label{sec:SOeven_poincare}
Recall from Section \ref{sec:SOevenC}:
\begin{eqnarray*}
&& \Delta=\{\alpha_i=\theta_i-\theta_{i+1}\mid i=1,\ldots, n-1\}\cup
\{ \alpha_n=\theta_{n-1}+ \theta_n\}\\
&& \Delta^\vee=\{ \alpha_i^\vee =e_i - e_{i+1}\mid i=1,\ldots,n-1\} \cup
\{\alpha_n^\vee = e_{n-1}+e_n\} \\
&& \pi_1(H)=\bigoplus_{i=1}^n \bZ e_i,\quad
\Lambda=\bigoplus_{i=1}^{n-1} \bZ(e_i-e_{i+1}) \oplus \bZ(e_{n-1}+e_n),\\
&&\pi_1(SO(2n))= \langle e_n \rangle \cong \bZ/2\bZ
\end{eqnarray*}

We now apply Theorem \ref{thm:Pt} to the case $G_\bR=SO(2n)$.
\begin{eqnarray*}
&& \varpi_{\alpha_i}= \theta_1+\cdots + \theta_i,
\quad i=1,\ldots, n-2\\
&&\varpi_{\alpha_{n-1}}=\frac{1}{2}(\theta_1+\cdots+\theta_{n-1}-\theta_n),\quad
\varpi_{\alpha_n}=\frac{1}{2}(\theta_1+\cdots+\theta_{n-1}+\theta_n) \\
&&
\varpi_{\alpha_i}(ke_n)=\left\{
\begin{array}{ll}
0 & i\leq n-2\\
-k/2 & i=n-1\\
k/2 & i=n
\end{array} \right.
\end{eqnarray*}
We have the following four cases:\\

Case 1. $\alpha_{n-1}, \alpha_n\in I$: $n_r=1$

\begin{eqnarray*}
&& I= \{ \alpha_{n_1}, \alpha_{n_1+n_2},\ldots, \alpha_{n_1+\cdots + n_{r-2}},\alpha_{n-1},
\alpha_n  \}\\
&& L^I= GL(n_1,\bC)\times \cdots \times  GL(n_{r-1},\bC)\times GL(1,\bC), \quad n_1+\cdots
+n_{r-1}+1=n \\
&& \dim_\bC \fz_{L^I}-\dim_\bC \fz_{SO(2n,\bC)}=r,\quad
\dim_\bC U^I = \sum_{1\leq i<j \leq r}n_i n_j +\frac{n(n-1)}{2}\\
&&\rho^I =\frac{1}{2}\sum_{i=1}^r
\biggl(n-2\sum_{j=1}^{i} n_j +n_i \biggr)
\biggl(\sum_{j=1}^{n_i}\theta_{n_1+\cdots+n_{i-1}+j}\biggr)
 +\frac{n-1}{2}(\theta_1+\cdots +\theta_n)\\
&&\langle \rho^I,  \alpha_{n_1+\cdots+n_i}^\vee
\rangle = \frac{n_i+n_{i+1}}{2},\quad \textup{ for }i=1,\dots,r-2,
\\&&
\langle \rho^I, \alpha_{n-1}^\vee\rangle
=\langle \rho^I, \alpha_n^\vee \rangle= \frac{n_{r-1}+1}{2}
\end{eqnarray*}
\\

Case 2. $\alpha_{n-1}\in I$, $\alpha_n\notin I$: $n_r>1$

\begin{eqnarray*}
&& I= \{ \alpha_{n_1}, \alpha_{n_1+n_2},\ldots, \alpha_{n_1+\cdots + n_{r-1}}, \alpha_{n-1}  \}\\
&& L^I= GL(n_1,\bC)\times \cdots \times  GL(n_{r},\bC),
\quad n_1+\cdots +n_r=n \\
&& \dim_\bC \fz_{L^I}-\dim_\bC \fz_{SO(2n,\bC)}=r,\quad \dim_\bC U^I
= \sum_{1\leq i<j\leq r}n_i n_j +\frac{n(n-1)}{2}\\ && \rho^I
=\frac{1}{2}\sum_{i=1}^r \biggl(n-2\sum_{j=1}^{i} n_j +n_i \biggr)
\biggl(\sum_{j=1}^{n_i}\theta_{n_1+\cdots+n_{i-1}+j}\biggr)+
\frac{n-1}{2}(\sum_{j=1}^{n}\theta_j) -(n_r-1)\theta_n
\end{eqnarray*}

$$
\langle \rho^I,  \alpha_{n_1+\cdots+n_i}^\vee \rangle = \frac{n_i+n_{i+1}}{2}
\textup{ for }i=1,\dots,r-1,\quad
\langle \rho^I, \alpha_{n-1}^\vee\rangle
= n_r-1
$$
\\

Case 3. $\alpha_{n-1}\notin I$, $\alpha_n\in I$: $n_r>1$
\begin{eqnarray*}
&& I=\{ \alpha_{n_1}, \alpha_{n_1+n_2},\ldots, \alpha_{n_1+n_2+\cdots+n_{r-1}},\alpha_n\}\\
&& L^I=GL(n_1,\bC)\times \cdots \times GL(n_r,\bC),\quad n_1+\cdots+n_r =n\\
&& \dim_\bC \fz_{L^I}-\dim_\bC \fz_{SO(2n,\bC)}=r,\quad \dim_\bC U^I = \sum_{1\leq i<j\leq r}n_i n_j
+\frac{n(n-1)}{2}\\
&& \rho^I =\frac{1}{2}\sum_{i=1}^r
\biggl(n-2\sum_{j=1}^{i} n_j +n_i \biggr)
\biggl(\sum_{j=1}^{n_i}\theta_{n_1+\cdots+n_{i-1}+j}\biggr)
+\frac{n-1}{2}(\theta_1+\cdots +\theta_n)
\\&&
\langle \rho^I,  \alpha_{n_1+\cdots+n_i}^\vee \rangle = \frac{n_i+n_{i+1}}{2}
\textup{ for }i=1,\dots,r-1,\quad
\langle \rho^I, \alpha_n^\vee\rangle
= n_r-1
\end{eqnarray*}
\\

Case 4. $\alpha_{n-1}\notin I$, $\alpha_n\notin I$: $n_r>1$
\begin{eqnarray*}
&& I=\{ \alpha_{n_1}, \alpha_{n_1+n_2},\ldots, \alpha_{n_1+n_2+\cdots+n_{r-1}} \}\\
&& L^I=GL(n_1,\bC)\times \cdots \times GL(n_{r-1},\bC)\times SO(2n_r),
\quad n_1+\cdots+n_r =n\\
&& \dim_\bC \fz_{L^I}-\dim_\bC \fz_{SO(2n,\bC)}=r-1,\\
&& \dim_\bC U^I = \sum_{1\leq i<j\leq r}n_i n_j
+\frac{n(n-1)-n_r(n_r-1)}{2}\\
&& \rho^I =\frac{1}{2}\sum_{i=1}^r
\biggl(n-2\sum_{j=1}^{i} n_j +n_i \biggr)
\biggl(\sum_{j=1}^{n_i}\theta_{n_1+\cdots+n_{i-1}+j}\biggr) \\
&& \quad\quad+
\frac{n-1}{2}(\theta_1+\cdots +\theta_{n_1+\cdots+n_{r-1}})
+ \frac{n-n_r}{2}(\theta_{n_1+\cdots+n_{r-1}+1}+\cdots +\theta_n)\\
&&
\langle \rho^I,  \alpha_{n_1+\cdots+n_i}^\vee \rangle = \frac{n_i+n_{i+1}}{2},\quad
\textup{ for }i=1,\dots,r-2,\\
&& \langle \rho^I, \alpha_{n_1+\cdots+n_{r-1}}^\vee\rangle
= \frac{n_{r-1}+2n_r-1}{2}
\end{eqnarray*}

We have the following closed formula for the $SO(2n)$-equivariant Poincar\'{e} series
of the representation variety of flat $SO(2n)$-connections over $\Si^\ell_0$:

\begin{thm}\label{thm:PtSOeven} $n\geq 2$
\begin{eqnarray*}
&& P_t^{SO(2n)}(X_{\mathrm{flat}}^{\ell,0}(SO(2n))^{(-1)^k} ) =\\
&&  \sum_{r=2}^n\sum_{\tiny \begin{array}{c}n_1,\ldots, n_r\in\bZ_{>0}\\ \sum n_j=n, n_r=1 \end{array}}
(-1)^r \prod_{i=1}^r \frac{\prod_{j=1}^{n_i} (1+t^{2j-1})^{2\ell} }
{ (1-t^{2n_i})\prod_{j=1}^{n_i-1}(1-t^{2j})^2 }\\
&& \cdot
 \frac{t^{(\ell-1)(2\sum_{i<j} n_i n_j+n(n-1))} }
{\left[\prod_{i=1}^{r-1}(1-t^{2(n_i+n_{i+1})} ) \right](1-t^{2(n_{r-1}+1)}) }
t^{2\sum_{i=1}^{r-2}(n_i+ n_{i+1})+ 4(n_{r-1}+1)\langle k/2\rangle }\\
&&+ \sum_{r=1}^{n-1}
\sum_{\tiny \begin{array}{c}n_1,\ldots, n_r\in\bZ_{>0}\\ \sum n_j=n, n_r>1 \end{array}}
\left(2(-1)^r \prod_{i=1}^r \frac{\prod_{j=1}^{n_i} (1+t^{2j-1})^{2\ell} }
{ (1-t^{2n_i})\prod_{j=1}^{n_i-1}(1-t^{2j})^2 } \right.\\
&& \cdot \frac{t^{(\ell-1)(2\sum_{i<j} n_i n_j+n(n-1))} }
{\left[\prod_{i=1}^{r-1}(1-t^{2(n_i+n_{i+1} )} )\right](1-t^{4(n_r-1)}) }
\cdot t^{2\sum_{i=1}^{r-1}(n_i +n_{i+1})+4 (n_r-1)\langle k/2\rangle}\\
&& +(-1)^{r-1} \prod_{i=1}^{r-1} \frac{\prod_{j=1}^{n_i} (1+t^{2j-1})^{2\ell} }
{ (1-t^{2n_i})\prod_{j=1}^{n_i-1}(1-t^{2j})^2 }
\cdot\frac{(1+t^{2n_r-1})^{2\ell}\prod_{j=1}^{n_r-1}(1+t^{4j-1})^{2\ell}}
{(1-t^{2n_r-2})(1-t^{2n_r})\prod_{j=1}^{2n_r-2}(1-t^{2j})}\\
&&\left. \cdot
 \frac{t^{(\ell-1)(2\sum_{i<j} n_i n_j+n(n-1)-n_r(n_r-1))} }
{\left[\prod_{i=1}^{r-2}(1-t^{2(n_i+n_{i+1})} ) \right](1-\epsilon(r)t^{2(n_{r-1}+2n_r-1)}) }
t^{2\sum_{i=1}^{r-2}(n_i+ n_{i+1})+ 2\epsilon(r)(n_{r-1}+2n_r-1) }\right)
\end{eqnarray*}
where
$$
\epsilon(r)=\left\{\begin{array}{ll}0 & r=1\\ 1 & r>1 \end{array}\right.
$$
\end{thm}
\begin{rem}
For $n\geq 2$, we have
$$
 P_t^{SO(2n)}(X_{\mathrm{flat}}^{\ell,0}(SO(2n))^{+1} )
= P_t^{Spin(2n)}(X_{\mathrm{flat}}^{\ell,0}(Spin(2n)) ),
$$
so Theorem \ref{thm:PtSOeven} also gives a formula
for $X_{\mathrm{flat}}^{\ell,0}(Spin(2n))$.
\end{rem}

\begin{ex}
\begin{eqnarray*}
&& P^{SO(4)}_t(X_{\mathrm{flat}}^{\ell,0}(SO(4))^{+1})
= P^{Spin(4)}_t(X_{\mathrm{flat}}^{\ell,0}(Spin(4)))\\
&=& \frac{(1+t)^{4\ell} t^{4\ell+ 4} }{(1-t^2)^2(1-t^4)^2}
-2\frac{(1+t)^{2\ell}(1+t^3)^{2\ell}t^{2\ell+2}}{(1-t^2)^2(1-t^4)^2}
+\frac{(1+t^3)^{4\ell}}{(1-t^2)^2(1-t^4)^2}\\
&=&\frac{1}{(1-t^2)^2 (1-t^4)^2}
\left( (1+t^3)^{4\ell}-2t^{2\ell+2}(1+t)^{2\ell}(1+t^3)^{2\ell}+ t^{4\ell+4}(1+t)^{4\ell}\right)\\
&& P^{SO(4)}_t(X_{\mathrm{flat}}^{\ell,0}(SO(4))^{-1})\\
&=& \frac{(1+t)^{4\ell} t^{4\ell+ 2} }{(1-t^2)^2(1-t^4)^2}
-2\frac{(1+t)^{2\ell}(1+t^3)^{2\ell}t^{2\ell}}{(1-t^2)^2(1-t^4)^2}
+\frac{(1+t^3)^{4\ell}}{(1-t^2)^2(1-t^4)^2}\\
&=&\frac{(1+t)^{2\ell}}{(1-t^2)^2 (1-t^4)^2}
\left( (1+t^3)^{4\ell}-2t^{2\ell}(1+t)^{2\ell}(1+t^3)^{2\ell}+ t^{4\ell}(1+t)^{4\ell}\right)
\end{eqnarray*}
Note that $Spin(4)=SU(2)\times SU(2)$, so
$$
 P^{Spin(4)}_t(X_{\mathrm{flat}}^{\ell,0}(Spin(4)))
= \left(P^{SU(2)}_t(X_{\mathrm{flat}}^{\ell,0}(SU(2)))\right)^2
$$
as expected, where $P^{SU(2)}_t(X_{\mathrm{flat}}^{\ell,0}(SU(2)))$ is calculated
in Example \ref{UtwoSUtwo}.
\end{ex}

\begin{ex}
\begin{eqnarray*}
&& P^{SO(6)}_t(X_{\mathrm{flat}}^{\ell,0}(SO(6))^{+1})
= P^{Spin(6)}_t(X_{\mathrm{flat}}^{\ell,0}(Spin(6)))\\
&=&\frac{(1+t)^{4\ell}(1+t^3)^{2\ell} t^{10\ell+2} }
{(1-t^2)^3(1-t^4)(1-t^6)^2}
-\frac{(1+t)^{6\ell} t^{12\ell}}{(1-t^2)^3 (1-t^4)^3}\\
&&  -2\frac{(1+t)^{2\ell} (1+t^3)^{2\ell}(1+t^5)^{2\ell} t^{6\ell+2}}
{(1-t^2)^2 (1-t^4)^2(1-t^6)(1-t^8)}
 +2\frac{(1+t)^{4\ell}(1+t^3)^{2\ell} t^{10\ell}}
{(1-t^2)^3(1-t^4)^2(1-t^6)}\\
&& +\frac{(1+t^3)^{2\ell}(1+t^5)^{2\ell}(1+t^7)^{2\ell}}
{(1-t^2)(1-t^4)^2(1-t^6)^2(1-t^8)}
-\frac{(1+t)^{2\ell}(1+t^3)^{4\ell} t^{8\ell}}
{(1-t^2)^3(1-t^4)^2(1-t^8)}\\
&&P^{SO(6)}_t(X_{\mathrm{flat}}^{\ell,0}(SO(6))^{-1})\\
&=&\frac{(1+t)^{4\ell}(1+t^3)^{2\ell} t^{10\ell-4}}
{(1-t^2)^3(1-t^4)(1-t^6)^2}
-\frac{(1+t)^{6\ell} t^{12\ell-4}}{(1-t^2)^3 (1-t^4)^3}\\
&& -2\frac{(1+t)^{2\ell} (1+t^3)^{2\ell}(1+t^5)^{2\ell} t^{6\ell-2}}
{(1-t^2)^2 (1-t^4)^2(1-t^6)(1-t^8)}
 +2\frac{(1+t)^{4\ell}(1+t^3)^{2\ell} t^{10\ell-2}}
{(1-t^2)^3(1-t^4)^2(1-t^6)}\\
&& +\frac{(1+t^3)^{2\ell}(1+t^5)^{2\ell}(1+t^7)^{2\ell}}
{(1-t^2)(1-t^4)^2(1-t^6)^2(1-t^8)}
 -\frac{(1+t)^{2\ell}(1+t^3)^{4\ell} t^{8\ell}}
{(1-t^2)^3(1-t^4)^2(1-t^8)}
\end{eqnarray*}
Note that $Spin(6)=SU(4)$, so
$$
P^{Spin(6)}_t(X_{\mathrm{flat}}^{\ell,0}(Spin(6)))
=P^{SU(4)}_t(X_{\mathrm{flat}}^{\ell,0}(SU(4)))
$$
as expected, where $P^{SU(4)}_t(X_{\mathrm{flat}}^{\ell,0}(SU(4)))$
is calculated in Example \ref{SUthreeSUfour}.
\end{ex}

\subsection{$SO(4m+2)$-connections on nonorientable surfaces}
\label{sec:SO4m+2_nonorientable} In this subsection, we consider
$SO(2n)$ where $n=2m+1$ is odd, so that
$$
\BC_0^\tau
=\{ \sqrt{-1}\diag(t_1 J,\ldots,t_{2m} J, 0J )\mid t_1\geq \cdots\geq t_{2m}\geq 0 \}.
$$
Any $\mu\in\BC_0^\tau$ is of the form
$$
\mu=\sqrt{-1}\diag(\lambda_1 J_{n_1},\ldots, \lambda_{r-1} J_{n_{r-1}},0 J_{n_r}),
$$
where $\lambda_1>\cdots > \lambda_{r-1}>0$ and $n_r>0$. We have
$$
SO(2n)_{X_\mu}\cong \Phi(U(n_1))\times\cdots\times\Phi(U(n_{r-1}))\times SO(2n_r),
$$
where $X_\mu=-2\pi\sqrt{-1}\mu$.

Given $\mu\in \BC_0^\tau$, let
$$
\ep_\mu =\diag(H_{n-n_r}, (-1)^{(n-n_r)} I_1, I_{2n_r-1}).
$$
Then $\Ad(\ep_\mu) X_\mu = -X_\mu$. Note that $n_r\geq 1$.

Suppose that
$$
(\ab,\ep_\mu c',X_\mu/2)\in X_{\mathrm{YM}}^{\ell,1}(SO(2n)).
$$
Then
$$
\exp(X_\mu/2)\ep_\mu c' \ep_\mu c'=\pab
$$
where $a_i,~b_i,~c'\in \Phi(U(n_1))\times\cdots \times  \Phi(U(n_{r-1}))\times SO(2n_r)$.

Let $L:\bR^{2(n-n_r)}\to \bC^{n-n_r}$ be defined as in Section
\ref{sec:SOodd_orientable}, and let
\begin{eqnarray*}
X_\mu' &=& L\circ\left(2\pi\diag(\lambda_1 J_{n_1},\ldots,\lambda_{r-1} J_{n_{r-1}})\right)
\circ L^{-1}\\
&=& 2\pi\sqrt{-1}\diag(\lambda_1I_{n_1},\cdots,\lambda_{r-1}I_{n_{r-1}})\in
\mathfrak{u}(n_1)\times\cdots\times\mathfrak{u}(n_{r-1}).
\end{eqnarray*}
Then the condition on $X'_\mu$ is
$$
\exp(X'_\mu/2)\bar{c'} c'=\pab \in SU(n_1)\times\cdots\times SU(n_{r-1})
$$
where $a_i,~b_i,~c'\in U(n_1)\times\cdots \times  U(n_{r-1})$,
and $\bar{c'}$ is the complex conjugate
of $c'$. In order that this is nonempty, we need
$1=\det(e^{\pi\sqrt{-1}\lambda_j}I_{n_j})$, i.e.,
\begin{equation}
\lambda_j=\frac{2k_j}{n_j}, \quad k_j \in \bZ, \quad j=1,\ldots,r-1.
\end{equation}

Similarly, suppose that
$(\ab,d, \ep_\mu c',X_\mu/2)\in X_{\mathrm{YM}}^{\ell,2}(SO(2n)$. Then
\[
\exp(X_\mu/2)(\ep_\mu c') d(\ep_\mu c')^{-1}d =\pab
\]
where $a_i,~b_i,~d,~c'\in \Phi(U(n_1))\times\cdots \times  \Phi(U(n_{r-1}))\times SO(2n_r)$.
The condition on $X'_\mu$ is
$$
\exp(X'_\mu/2)\bar{c'}\bar{d}\bar{c'}^{-1}d \in
SU(n_1)\times \cdots \times SU(n_{r-1}).
$$
Again, we need
$$
\lambda_j=\frac{2k_j}{n_j},\quad  k_j \in\bZ, \quad j=1,\ldots,r-1.
$$

We conclude that for nonorientable surfaces,
$$
\mu= \sqrt{-1}\diag\Bigl(\frac{2k_1}{n_1} J_{n_1},\ldots, \frac{2k_{r-1}}{n_{r-1}} J_{n_{r-1}},
0 J_{n_r}\Bigr),
\textup{where }
k_j \in \bZ, ~\frac{k_i}{n_i}>\frac{k_{i+1}}{n_{i+1}}>0, n_r>0.
$$

For each $\mu$, the $\ep_\mu$-reduced representation varieties are
\begin{eqnarray*}
V_{\mathrm{YM}}^{\ell,1}(SO(2n))_\mu &=&  \{ (\ab,c')\in SO(2n)_{X_\mu}^{2\ell+1}\mid\\
&& \quad \pab=\exp(\frac{X_\mu}{2})\ep_\mu c' \ep_\mu c' \}, \\
 V_{\mathrm{YM}}^{\ell,2}(SO(2n))_\mu&=& \{ (\ab,d,c')\in SO(2n)_{X_\mu}^{2\ell+2}\mid\\
&& \quad \pab=\exp(\frac{X_\mu}{2})\ep_\mu c' d(\ep_\mu c')^{-1}d \}.
\end{eqnarray*}

For $i= 1,\cdots,\ell$, write
\begin{eqnarray*}
&& a_i=\diag(A^i_1,\ldots, A^i_r),\quad
b_i=\diag(B^i_1,\ldots, B^i_r),\\
&& c'=\diag(C_1,\ldots,C_r),\quad
d =\diag(D_1,\ldots,D_r),
\end{eqnarray*}
where $A^i_j, ~B^i_j,~C_j,~D_j \in \Phi(U(n_{j}))$ for $j=1,\ldots,r-1$,
and $A^i_r, ~B^i_r,~C_r,~D_r\in SO(2n_r)$.

\paragraph{$i=1$.} Define $V_j$ as in \eqref{eqn:VjRP}, and define
\begin{equation}\label{eqn:VzeroRP}
V_r=
\Bigl\{(A^1_r,B^1_r,\ldots, A^\ell_r,B^\ell_r,C_r)
\in SO(2n_r)^{2\ell+1}\mid \prod_{i=1}^\ell[A^i_r,B^i_r]= (\epsilon C_r)^2
  \Bigr \},
\end{equation}
where $\ep=\diag((-1)^{n-n_r}, I_{2n_r-1})$, $\det(\ep)=(-1)^{n-n_r}$.
Let $C_r'=\ep C_r$. We see that
\begin{eqnarray*}
V_r &\cong&
\Bigl\{(A^1_r,B^1_r,\ldots, A^\ell_r,B^\ell_r,C_r')
\in SO(2n_r)^{2\ell}\times O(2n_r)\mid \\
&& \prod_{i=1}^\ell[A^i_r,B^i_r]= (C_r')^2, \det(C_r')=(-1)^{n-n_r} \Bigr \}\\
&\cong & V^{\ell,1}_{O(2n_r),(-1)^{n-n_r}}
\end{eqnarray*}
where
$V^{\ell,1}_{O(n),\pm 1}$ is the twisted representation variety defined in (\ref{eqn:Oone})
of Section \ref{sec:twistedO}. $V^{\ell,1}_{O(n),\pm 1}$ is nonempty if $\ell\geq 2$. We have shown
that $V^{\ell,1}_{O(n),\pm 1}$ is disconnected with two components
$V^{\ell,1,+ 1}_{O(n),\pm 1}$ and $V^{\ell,1,-1}_{O(n),\pm 1}$
if $\ell\geq 2$ and
$n>2$ (Proposition \ref{thm:tO}).

We have
$$
V_{\mathrm{YM}}^{\ell,1}(SO(2n))_\mu =
\prod_{j=1}^r V_j.
$$
We define a  $U(n_j)$-action on $V_j=\tV^{\ell,1}_{n_j,-k_j}$ by \eqref{eqn:actI}
of Section \ref{sec:twisted}, and an $SO(2n_r)$-action on $V_r=V^{\ell,1}_{O(2n_r),(-1)^{n-n_r}}$
by \eqref{eqn:OactI} of Section \ref{sec:twistedO}. Then we have a homeomorphism
$$
V_{\mathrm{YM}}^{\ell,1}(SO(2n))_\mu/SO(2n)_{X_\mu}\cong
\prod_{j=1}^{r-1} (V_j/U(n_j)) \times V_r/SO(2n_r)
$$
and a homotopy equivalence
$$
V_{\mathrm{YM}}^{\ell,1}(SO(2n))_\mu ^{h SO(2n)_{X_\mu} }\sim
\prod_{j=1}^{r-1} V_j^{h U(n_j)} \times {V_r}^{h SO(2n_r)}.
$$

\paragraph{$i=2$}
Define $V_j$ as in \eqref{eqn:VjKlein}, and define
\begin{equation}\label{eqn:VzeroKlein}
V_r= \Bigl\{(A^1_r,B^1_r,\ldots, A^\ell_r,B^\ell_r,D_r,C_r)
\in SO(2n_r)^{2\ell+2}\mid \prod_{i=1}^\ell[A^i_r,B^i_r]=\epsilon C_r D_r
(\epsilon C_r)^{-1}D_r
  \Bigr \}
\end{equation}
where $\ep=\diag((-1)^{n-n_r}I_1, I_{2n_r-1})$, $\det(\ep)=(-1)^{n-n_r}$.
Let $C_r'=\ep C_r$. We see that
\begin{eqnarray*}
V_r &\cong&
\Bigl\{(A^1_r,B^1_r,\ldots, A^\ell_r,B^\ell_r,D_r,C_r')
\in SO(2n_r)^{2\ell+1}\times O(2n_r)\mid \\
&& \prod_{i=1}^\ell[A^i_r,B^i_r]= C_r' D_r C_r'^{-1} D_r, \det(C_r')=(-1)^{n-n_r} \Bigr \}\\
&\cong & V^{\ell,2}_{O(2n_r),(-1)^{n-n_r}}
\end{eqnarray*}
where
$V^{\ell,2}_{O(n),\pm 1}$ is the twisted representation variety defined in (\ref{eqn:Otwo})
of Section \ref{sec:twistedO}. $V^{\ell,2}_{O(n),\pm 1}$ is nonempty if $\ell\geq 4$. We have shown
that $V^{\ell,2}_{O(n),\pm 1}$ is disconnected with two components
$V^{\ell,1,+1}_{O(n),\pm 1}$ and $V^{\ell,1,-1}_{O(n),\pm 1}$ if $\ell\geq 4$ and
$n>2$ (Proposition \ref{thm:tO}).

We have
$$
V_{\mathrm{YM}}^{\ell,2}(SO(2n))_\mu =
\prod_{j=1}^r V_j.
$$
We define a $U(n_j)$-action on $V_j=\tV^{\ell,2}_{n_j,-k_j}$ by
\eqref{eqn:actII} of Section \ref{sec:twisted}, and an
$SO(2n_r)$-action on $V_r=V^{\ell,2}_{O(2n_r),(-1)^{n-n_r}}$ by
\eqref{eqn:OactII} of Section \ref{sec:twistedO}. Then we have a
homeomorphism
$$
V_{\mathrm{YM}}^{\ell,2}(SO(2n))_\mu/SO(2n)_{X_\mu}\cong
\prod_{j=1}^{r-1} (V_j/U(n_j)) \times V_r/SO(2n_r)
$$
and a homotopy equivalence
$$
{V_{\mathrm{YM}}^{\ell,2}(SO(2n))_\mu}^{h SO(2n)_{X_\mu} }\sim
\prod_{j=1}^{r-1} {V_j}^{h U(n_j)} \times {V_r}^{h SO(2n_r)}.
$$

Note that,
$V^{\ell,i}_{O(2),-1}\cong \tV^{\ell,i}_{1,0}\cong U(1)^{2\ell+i}$ is connected as
mentioned in Section \ref{sec:twistedO}.

We have seen that $V_{\mathrm{YM}}^{\ell,i}(SO(2n))_\mu$ is disconnected with two connected
components if $\ell\geq 2i$ and $n_r\geq 1$ (notice that when $n_r=1$, $n-n_r=2m$ is even).
To determine the underlying topological $SO(2n)$-bundle $P$ for each component,
we consider four special cases.

\paragraph{\em Case 1.} Assuming that $n_r>1$, we
consider special points
$$
(\ab,c)\in \VymSOE{\ell}{1}_\mu,\quad
(\ab,d,c)\in \VymSOE{\ell}{2}_\mu,
$$
where
\begin{eqnarray*}
&& a_i=\diag(A^i_1,\ldots, A^i_{r-1},I_{2n_r}),
\quad b_i=\diag(B^i_1,\ldots, B^i_{r-1},I_{2n_r}),\\
&&c=\ep_\mu=\diag(H_{n-n_r},(-1)^{n-n_r}I_1, I_{2n_r-1}),
\quad d=I_{2n}.
\end{eqnarray*}
Let $\ep_1=\diag((-1)^{n-n_r}I_1, I_{2n_r-1})$. Then
\begin{eqnarray*}
&& (A^i_j,B^i_j,\ldots,A^i_j, B^i_j) \in X^{\ell,0}_{\mathrm{YM}}
(U(n_j))_{-\frac{k_j}{n_j},\ldots,-\frac{k_j}{n_j}},\quad j=1,\ldots, r-1,\\
&& (I_{2n_r},\ldots, I_{2n_r},\ep_1) \in V^{\ell,1,(-1)^{n-n_r}}_{O(2n_r), (-1)^{n-n_r}},
(I_{2n_r},\ldots, I_{2n_r},I_{2n_r},\ep_1) \in V^{\ell,2,1}_{O(2n_r),(-1)^{n-n_r}}.
\end{eqnarray*}
We have $P=P_1\times P_2$, where $P_1$ is an $SO(2(n-n_r)+1)$-bundle, and
$P_2$ is an $SO(2n_r-1)$-bundle with trivial holonomies $I_{2n_r-1}$. We have
$$
w_2(P)=w_2(P_1) = k_1+\cdots+k_{r-1} + i\frac{(n-n_r)(n-n_r+1)}{2}\ (\mathrm{mod}\ 2),
$$
where the second equality follows from the argument in
Section \ref{sec:SOodd_nonorientable}.

\paragraph{\em Case 2.} Assuming that $n_r>1$, as in {\em Case 1},
we consider special points
$$
(\ab,c)\in \VymSOE{\ell}{1}_\mu,\quad (\ab,d,c)\in\VymSOE{\ell}{2}
$$
where
\begin{eqnarray*}
&& a_i=\diag(A^i_1,\ldots, A^i_{r-1},I_{2n_r}),
\quad b_i=\diag(B^i_1,\ldots, B^i_{r-1},I_{2n_r}),\\
&&c=\diag(H_{n-n_r},(-1)^{(n-n_r)}I_1,-I_2, I_{2n_r-3}),
\quad d=\diag(I_{2(n-n_r)+1}, -I_2, I_{2n_r-3}).
\end{eqnarray*}

Let $\ep_1=\diag((-1)^{n-n_r}I_1,-I_2, I_{2n_r-3})$,
$\ep_2=\diag(I_1,-I_2, I_{2n_r-3})$, and $\ep=\diag(-I_2, I_{2n_r-3})$. Then
\begin{eqnarray*}
&& (A^i_j,B^i_j,\ldots,A^i_j, B^i_j)\in X^{\ell,0}_{\mathrm{YM}}
(U(n_j))_{-\frac{k_j}{n_j},\ldots,-\frac{k_j}{n_j}},\quad j=1,\ldots, r-1,\\
&& (I_{2n_r},\ldots, I_{2n_r},\ep_1) \in V^{\ell,1,-(-1)^{n-n_r}}_{O(2n_r), (-1)^{n-n_r}},
\quad
(I_{2n_r},\ldots, I_{2n_r},\ep_2,\ep_1) \in V^{\ell,2,-1}_{O(2n_r),(-1)^{n-n_r}}.
\end{eqnarray*}
We have $P=P_1\times P_2$, where $P_1$ is an $SO(2(n-n_r)+1)$-bundle, and
$P_2$ is an $SO(2n_r-1)$-bundle with holonomies $a_i=b_i=I_{2n_r-1}$, $c=d=\ep$.
We can choose the lifting of $d$ and $c$ as $\tilde{d}=\tilde{c}=e_1e_2$ and
$\tilde{c}^2=\tilde{c}\tilde{d}\tilde{c}^{-1}\tilde{d}=-1$. Thus we have
\begin{eqnarray*}
w_2(P_1)= k_1+\cdots+k_{r-1} + i\frac{(n-n_r)(n-n_r+1)}{2}\ (\mathrm{mod}\ 2),\quad
w_2(P_2)=1\ (\mathrm{mod}\ 2),
\end{eqnarray*}
so
$$
w_2(P)=w_2(P_1)+w_2(P_2)= k_1+\cdots+k_{r-1} + i\frac{(n-n_r)(n-n_r+1)}{2} + 1\
(\mathrm{mod}\ 2).
$$

\paragraph{\em Case 3.} Assuming that $n_r=1$ so that $n-n_r = 2m$ is even,
we consider special points
$$
(\ab,c)\in \VymSOE{\ell}{1}_\mu,\quad (\ab,d,c)\in \VymSOE{\ell}{2}_\mu,
$$
where
\begin{eqnarray*}
&& a_i=\diag(A^i_1,\ldots, A^i_{r-1},I_2),
\quad b_i=\diag(B^i_1,\ldots, B^i_{r-1},I_2),\\
&&c= \diag(H_{2m},-I_2),\quad d=\diag(I_{4m}, -I_2).
\end{eqnarray*}
Then
\begin{eqnarray*}
&& (A^i_j,B^i_j,\ldots,A^i_j, B^i_j)\in X^{\ell,0}_{\mathrm{YM}}
(U(n_j))_{-\frac{k_j}{n_j},\ldots,-\frac{k_j}{n_j}},\quad j=1,\ldots, r-1,\\
&& (I_2,\ldots, I_2,-I_2) \in V^{\ell,1,-1}_{O(2), +1},
\quad
(I_2,\ldots, I_2,-I_2,-I_2) \in V^{\ell,2,-1}_{O(2),+1}.
\end{eqnarray*}
We have $P=P_1\times P_2$, where $P_1$ is an $SO(4m)$-bundle with holonomies $d=I_{4m}$
and $c=H_{2m}$ with lifting $\tilde{c}=e_2e_4\cdots e_{4m}$, and
$P_2$ is an $SO(2)$-bundle with holonomies $a_i=b_i=I_2$ and $c=d=-I_2$ with lifting
$\tilde{d}=\tilde{c}=e_1e_2$. Then we have
$$
w_2(P_1)=k_1+\cdots +k_{r-1} +im\ (\mathrm{mod}\ 2),\quad
w_2(P_2)=1\ (\mathrm{mod}\ 2),
$$
so
$$
w_2(P)= k_1+\cdots+k_{r-1} + im + 1.
$$

\paragraph{\em Case 4.} Assuming that $n_r=1$ as in {\em Case 3},
we consider special points
$$
(\ab,c)\in \VymSOE{\ell}{1}_\mu,\quad (\ab,d,c)\in \VymSOE{\ell}{2}_\mu,
$$
where
\begin{eqnarray*}
&& a_i=\diag(A^i_1,\ldots, A^i_{r-1},I_2),
\quad b_i=\diag(B^i_1,\ldots, B^i_{r-1},I_2),\\
&&c= \diag(H_{2m},I_2),\quad d=I_{2n}.
\end{eqnarray*}
Then
\begin{eqnarray*}
&& (A^i_j,B^i_j,\ldots,A^i_j, B^i_j)\in X^{\ell,0}_{\mathrm{YM}}
(U(n_j))_{-\frac{k_j}{n_j},\ldots,-\frac{k_j}{n_j}},\quad j=1,\ldots, r-1,\\
&& (I_2,\ldots, I_2,I_2) \in V^{\ell,1,1}_{O(2), +1},
\quad
(I_2,\ldots,I_2,I_2,I_2) \in V^{\ell,2,1}_{O(2),+1}.
\end{eqnarray*}
We have $P=P_1\times P_2$, where $P_1$ is an $SO(4m)$-bundle with holonomies $d=I_{4m}$
and $c=H_{2m}$ with lifting $\tilde{c}=e_2e_4\cdots e_{4m}$, and
$P_2$ is an $SO(2)$-bundle with trivial holonomies $I_2$. Then we have
$$
w_2(P)=w_2(P_1)=k_1+\cdots +k_{r-1} +im\ (\mathrm{mod}\ 2).
$$
To summarize, when $n=2m+1$, we have
$$
\VymSOE{\ell}{i}_\mu^\pm  =\prod_{j=1}^{r-1} V_j \times
V_{O(2n_r),(-1)^{n-n_r}  }^{\ell, i,
\pm (-1)^{k_1+\cdots + k_{r-1} + i\frac{(n-n_r)(n-n_r-1)}{2}}  },
$$
where $\VymSOE{\ell}{i}_\mu^\pm$ is the $\ep_\mu$-reduced version of
$X_{\mathrm{YM}}^{\ell,i}(SO(2n))^{\pm 1}_\mu$. Note that
$$
i\frac{(n-n_r)(n-n_r-1)}{2}\equiv i(m+\frac{n_r(n_r-1)}{2}) \ (\mathrm{mod}\ 2),
\quad n-n_r \equiv n_r-1\ (\mathrm{mod}\ 2).
$$

To simplify the notation, we write
$$
\mu=(\mu_1,\ldots,\mu_{2m},0)=\Bigl(\underbrace{\frac{2k_1}{n_1},\ldots,\frac{2k_1}{n_1} }_{n_1}, \ldots,
\underbrace{\frac{2k_{r-1}}{n_{r-1}},\ldots,\frac{2k_{r-1}}{n_{r-1}} }_{n_{r-1}},
\underbrace{0,\ldots,0}_{n_r}\Bigr)
$$
instead of
$$
\sqrt{-1}\diag\Bigl(\frac{2k_1}{n_1}J_{n_1},\ldots, \frac{2k_{r-1}}{n_{r-1}} J_{n_{r-1}}, 0 J_{n_r}\Bigr).
$$
Let
\begin{eqnarray*}
\hat{I}_{SO(4m+2)}&=&\Bigl \{
\mu=\Bigl(\underbrace{\frac{2k_1}{n_1},\ldots,\frac{2k_1}{n_1}}_{n_1},
\ldots,\underbrace{\frac{2k_{r-1}}{n_{r-1}},\ldots,\frac{2k_{r-1}}{n_{r-1}}}_{n_{r-1}},
\underbrace{0,\ldots,0}_{n_r}\Bigr) \bigl|
n_j\in \bZ_{>0},\\&&\quad n_1+\cdots+n_r=n=2m+1,\ k_j\in \bZ,\
\frac{k_1}{n_1}>\cdots >\frac{k_{r-1}}{n_{r-1}}>0\Bigr\}
\end{eqnarray*}
Recall that the twisted moduli spaces for $U(n)$ are defined by $\tilde{\cM}^{\ell,i}_{n,k}=\tV^{\ell,i}_{n,k}/U(n)$,
where $i=1,2$. Also we define the twisted moduli spaces for $SO(n)$ by
$$
\cM^{\ell,i,\pm 1}_{O(n),\pm1}=V^{\ell,i,\pm 1}_{O(n),\pm1}/SO(n),\quad
\textup{where } i=1,2.
$$

\begin{pro}\label{thm:muSO4m+2_nonorientable}
Suppose that  $\ell\geq  2i$, where $i=1,2$.
Let
\begin{equation}\label{eqn:muSO4m+2_nonorientable}
\mu=\Bigl(\underbrace{\frac{2k_1}{n_1},\ldots,\frac{2k_1}{n_1} }_{n_1}, \ldots,
\underbrace{\frac{2k_{r-1}}{n_{r-1}},\ldots,\frac{2k_{r-1}}{n_{r-1}} }_{n_{r-1}},
\underbrace{0,\ldots,0}_{n_r}\Bigr)\in \hat{I}_{SO(4m+2)}.
\end{equation}
Then $X_{\mathrm{YM}}^{\ell,i}(SO(4m+2))_\mu$ has two connected components
(from both bundles over $\Si^\ell_i$)
$$
X_{\mathrm{YM}}^{\ell,i}(SO(4m+2))_\mu^{+1},\quad \textup{and}\quad
X_{\mathrm{YM}}^{\ell,i}(SO(4m+2))_\mu^{-1}.
$$
We have a homeomorphism
\begin{eqnarray*}
&&X_{\mathrm{YM}}^{\ell,i}(SO(4m+2))_\mu^{\pm 1} /SO(4m+2)\\
&\cong & \prod_{j=1}^{r-1} \tilde{\cM}^{\ell,i}_{n_j,-k_j}\times
\cM^{\ell,i,\pm (-1)^{k_1+\cdots+k_{r-1}+im+ i\frac{(n_r)(n_r-1)}{2} }}_{O(2n_r),
(-1)^{n_r-1}}
\end{eqnarray*}
and a homotopy equivalence
\begin{eqnarray*}
&& {X_{\mathrm{YM}}^{\ell,i}(SO(4m+2))_\mu^{\pm 1} }^{h SO(4m+2) }\\
&\sim&\prod_{j=1}^{r-1}\Bigl(\tV^{\ell,i}_{n_j,-k_j} \Bigr)^{h U(n_j)}\times
 \Bigl(V^{\ell,i,\pm (-1)^{k_1+\cdots+k_{r-1}+im+ i\frac{(n_r)(n_r-1)}{2} }}_{O(2n_r),
(-1)^{n_r-1}}\Bigr)^{h SO(2n_r)}.
\end{eqnarray*}
\end{pro}

\begin{pro}Suppose that $\ell\geq 2i$, where $i=1,2$.
The connected components of
$X_{\mathrm{YM}}^{\ell,i}(SO(4m+2))^{\pm 1}$ are
$$
\{X_{\mathrm{YM}}^{\ell,i}(SO(4m+2))_\mu^{\pm 1}\mid \mu\in \hat{I}_{SO(4m+2)}\}.
$$
\end{pro}

Notice that, the set
$\{\mu=\sqrt{-1}\diag(\mu_1J,\ldots,\mu_{2m}J,0 J) \mid (\mu_1,\ldots,\mu_{2m},0)\in\hat{I}_{SO(4m+2)}\}
$ is a {\em proper} subset of $\{\mu\in(\Xi^I_+)^\tau \mid I\subseteq\Delta,\tau(I)=I\}$
as mentioned in Section \ref{sec:connected}.

The following is an immediate consequence of
Proposition \ref{thm:muSO4m+2_nonorientable}.
\begin{thm}\label{thm:PtmuSO4m+2_nonorientable}
  Suppose that $\ell\geq 2i$, where $i=1,2$,
and let $\mu$ be as in \eqref{eqn:muSO4m+2_nonorientable}. Then
\begin{eqnarray*}
&& P_t^{SO(4m+2)}\left(X_{\mathrm{YM}}^{\ell,i}(SO(4m+2))_\mu^{\pm 1}\right)\\
&=&\prod_{j=1}^{r-1} P_t^{U(n_j)}(\tV^{\ell,i}_{n_j,-k_j})\cdot
 P_t^{SO(2n_r)}\left(V^{\ell,i, \pm (-1)^{k_1+\cdots+k_{r-1}+
im+i\frac{n_r(n_r-1)}{2}}}_{O(2n_r),(-1)^{n_r-1} }\right).
\end{eqnarray*}
\end{thm}

\subsection{$SO(4m)$-connections on nonorientable surfaces}
\label{sec:SO4m_nonorientable}

In this subsection, we consider $SO(2n)$ where $n=2m$ is even, so
that $\BC_0^\tau=\BC_0$. There are four cases.

\paragraph{Case 1. $t_{n-1}>|t_n|$, $n_r=1$}
$$
\mu= \sqrt{-1}\diag(\lambda_1 J_{n_1},\ldots, \lambda_{r-1} J_{n_{r-1}}, \lambda_r J),
$$
where $\lambda_1 >\cdots >\lambda_{r-1} > |\lambda_r|\geq 0$. Thus
$$
SO(4m)_{X_\mu}\cong \Phi(U(n_1))\times\cdots\times
\Phi(U(n_{r-1}))\times \Phi(U(n_r)).
$$
where $X_\mu=-2\pi\sqrt{-1}\mu$.

Let $\ep=H_{2m}$. Suppose that $(\ab,\ep c',X_\mu/2)\in X_{\mathrm{YM}}^{\ell,1}(SO(4m))$. Then
$$
\exp(X_\mu/2)\epsilon c' \epsilon c'=\pab
$$
where $a_i,~b_i,~c'\in \Phi(U(n_1))\times\cdots \times \Phi( U(n_r))$.

Let $L:\bR^{2n}\to \bC^{n}$ be defined as in Section
\ref{sec:SOodd_orientable}, and let
\begin{eqnarray*}
X_\mu' &=& L\circ X_\mu \circ L^{-1}\\
&=&2\pi\sqrt{-1}\diag(\lambda_1 I_{n_1},\cdots,\lambda_rI_{n_r})\in
\mathfrak{u}(n_1)\times\cdots\times\mathfrak{u}(n_r).
\end{eqnarray*}
Then the condition on $X'_\mu$ is
$$
\exp(X'_\mu/2)\bar{c'} c'=\pab \in SU(n_1)\times\cdots\times SU(n_{r-1})\times \{I_2\}
$$
where $a_i,~b_i,~c'\in U(n_1)\times\cdots \times  U(n_r)$, and $\bar{c'}$ is the
complex conjugate of $c'$.
In order that this is nonempty, we need
$1=\det(e^{\pi\sqrt{-1}\lambda_j}I_{n_j})$, i.e.,
$$
\lambda_j=\frac{2k_j}{n_j}, \quad k_j \in \bZ, \quad j=1,\ldots,r.
$$

Similarly, suppose that $(\ab,d, \ep c',X_\mu/2)\in X_{\mathrm{YM}}^{\ell,2}(SO(4m))$. Then
$$
\exp(X_\mu/2)(\ep c')d (\epsilon c')^{-1} d=\pab
$$
where $a_i,~b_i,~d,~c'\in \Phi(U(n_1))\times\cdots \times \Phi( U(n_r))$.
The condition on $X'_\mu$ is
$$
\exp(X'_\mu/2)\bar{c'}\bar{d} \bar{c'}^{-1}d =\pab \in SU(n_1)\times
\cdots\times SU(n_{r-1})\times \{I_2\},
$$
where $a_i,~b_i,~d,~c'\in U(n_1)\times\cdots \times  U(n_r)$, and $\bar{c'}$ is the
complex conjugate of $c'$. Again, we need
$$
\lambda_j=\frac{2k_j}{n_j}, \quad k_j\in\bZ, \quad j=1,\ldots,r.
$$

We conclude that for nonorientable surfaces
$$
\mu=\sqrt{-1}\diag\Bigl(\frac{2k_1}{n_1} J_{n_1},\ldots,
\frac{2k_{r-1}}{n_{r-1}} J_{n_{r-1}},2k_r J\Bigr),
 ~k_j \in\bZ, ~\frac{k_1}{n_1}>\cdots>\frac{k_{r-1}}{n_{r-1}}>|k_r|\geq 0.
$$

For each $\mu$, define $\ep$-reduced representation varieties
\begin{equation} \label{eqn:nonorient-SO4mRP}
\begin{aligned}
V_{\mathrm{YM}}^{\ell,1}(SO(4m))_\mu =& \{ (\ab,c')\in SO(4m)_{X_\mu}^{2\ell+1}\mid\\
& \quad \pab=\exp(\frac{X_\mu}{2}) \epsilon c' \epsilon c'\},
\end{aligned}
\end{equation}
\begin{equation}\label{eqn:nonorient-SO4mK}
\begin{aligned}
V_{\mathrm{YM}}^{\ell,2}(SO(4m))_\mu =& \{ (\ab,d,c')\in SO(4m)_{X_\mu}^{2\ell+2}\mid\\
& \quad \pab=\exp(\frac{X_\mu}{2}) \epsilon c'd (\epsilon c')^{-1}d\}.
\end{aligned}
\end{equation}

For $i= 1,\ldots,\ell$, write
\begin{eqnarray*}
&& a_i=\diag(A^i_1,\ldots, A^i_r),\quad
b_i=\diag(B^i_1,\ldots, B^i_r),\\
&& c'=\diag(C_1,\ldots,C_r),\quad
d =\diag(D_1,\ldots,D_r),
\end{eqnarray*}
where $A^i_j, ~B^i_j ,~C_j,~D_j\in \Phi(U(n_{j}))$.

Define $V_j$ as in \eqref{eqn:VjRP} when $i=1$,
and as in \eqref{eqn:VjKlein} when $i=2$.
Then
$V_j\cong \tV^{\ell,i}_{n_j,-k_j}$ is connected, and
$$
V_{\mathrm{YM}}^{\ell,1}(SO(4m))_\mu =
\prod_{j=1}^r V_j.
$$
Thus $V_{\mathrm{YM}}^{\ell,i}(SO(4m))_\mu$ is connected, and
it corresponds to connections on a fixed topological $SO(4m)$-bundle $P$.
By the argument in Section \ref{sec:SOodd_nonorientable},
$$
w_2(P)=k_1+\cdots+k_r + i\frac{2m(2m+1)}{2}
=k_1+\cdots+k_r +im \quad(\mathrm{mod}\ 2).
$$

Let $U(n_j)$ acts on $V_j\cong \tV^{\ell,i}_{n_j,-k_j}$ by
\eqref{eqn:actI} and \eqref{eqn:actII} in Section \ref{sec:twisted}
when $i=1$ and when $i=2$, respectively. Then we have a
homeomorphism
\begin{equation}\label{eqn:nonorient-SO4m-homeo}
V_{\mathrm{YM}}^{\ell,i}(SO(4m))_\mu/SO(4m)_{X_\mu}\cong \prod_{j=1}^r (V_j/U(n_j))
\end{equation}
and a homotopy equivalence
\begin{equation}\label{eqn:nonorient-SO4m-homotopy}
{V_{\mathrm{YM}}^{\ell,i}(SO(4m))_\mu }^{h SO(4m)_{X_\mu} }\cong
\prod_{j=1}^r {V_j}^{ hU(n_j)}.
\end{equation}

\paragraph{Case 2. $t_{n-1}=-t_n>0$, $n_r>1$}
$$
\mu = \sqrt{-1}\diag(\lambda_1 J_{n_1},\ldots,\lambda_{r-1} J_{n_{r-1}}, \lambda_{r} J_{n_r-1}, -\lambda_r J),
$$
where $\lambda_1 >\cdots >\lambda_r >0$,
Thus
$SO(4m)_{X_\mu}\cong \Phi(U(n_1))\times\cdots\times
\Phi(U(n_{r-1}))\times \Phi'(U(n_r))$, where
$\Phi:U(k)\hookrightarrow SO(2k)$ is the standard embedding,
and $\Phi':U(k)\hookrightarrow SO(2k)$ is defined as
in Section \ref{sec:SOeven_orientable}.

Let $\ep=H_{2m}$. Suppose that $(\ab,\ep c',X_\mu/2)\in X_{\mathrm{YM}}^{\ell,1}(SO(4m))$. Then
$$
\exp(X_\mu/2)\ep c' \ep c'=\pab
$$
where $a_i,~b_i,~c'\in \Phi(U(n_1))\times\cdots \times \Phi(U(n_{r-1})\times \Phi'( U(n_r)) $.

Let $L\oplus L':\bR^{2(n-n_r)}\oplus\bR^{2n_r}\rightarrow
\bC^{n-n_r}\oplus\bC^{n_r}$,
and let
$$
X'_\mu = (L\oplus L')\circ  X_\mu \circ(L\otimes L')^{-1}
=2\pi\sqrt{-1}\diag(\lambda_1I_{n_1},\ldots,\lambda_rI_{n_r})\in
\mathfrak{u}(n_1)\times\cdots\times\mathfrak{u}(n_r).
$$
Then the condition on $X'_\mu$ is
$$
\exp(X'_\mu/2)\bar{c'} c'=\pab \in SU(n_1)\times\cdots\times SU(n_r)$$
where $a_i,~b_i,~c'\in U(n_1)\times\cdots \times  U(n_r)$, and $\bar{c'}$ is the
complex conjugate of $c'$. In order that this is nonempty, we need
$1=\det(e^{\pi\sqrt{-1}\lambda_j}I_{n_j})$, i.e.,
$$
\lambda_j=\frac{2k_j}{n_j},\quad k_j\in\bZ,\quad j=1,\ldots,r.
$$

Similarly, suppose that $(\ab,d,\ep c',X_\mu/2)\in X_{\mathrm{YM}}^{\ell,2}(SO(4m))$. Then
$$
\exp(X_\mu/2)(\ep c')d (\ep c')^{-1}d=\pab,
$$
where $a_i,~b_i,~d,~c'\in \Phi(U(n_1))\times\cdots \times \Phi(U(n_{r-1}))
\times\Phi'( U(n_r))$.
The condition on $X'_\mu$ is
$$
\exp(X'_\mu/2)\bar{c'}\bar{d}\bar{c'}^{-1}d=\pab \in SU(n_1)\times\cdots\times SU(n_r),
$$
where $a_i,~b_i,~d,~c'\in U(n_1)\times \cdots \times U(n_r)$ and $\bar{d}$ is the complex
conjugate of $d$.

Again, we need
$$
\lambda_j=\frac{2k_j}{n_j},\quad k_j\in\bZ, \quad j=1,\ldots, r.
$$

We conclude that for nonorientable surfaces
$$
\mu=\sqrt{-1}\diag(\frac{2k_1}{n_1}J_{n_1},\ldots,\frac{2k_{r-1}}{n_{r-1}} J_{n_{r-1}},
\frac{2k_r}{n_r}J_{n_r-1},
-\frac{2k_r}{n_r}J),\ k_j\in \bZ,\
\frac{k_1}{n_1}>\cdots>\frac{k_r}{n_r}>0.
$$

For each $\mu$, define $\ep$-reduced representation varieties as in
\eqref{eqn:nonorient-SO4mRP} and \eqref{eqn:nonorient-SO4mK}.
For $i= 1,\ldots,\ell$, write
\begin{eqnarray*}
&& a_i=\diag(A^i_1,\ldots, A^i_r),\quad
b_i=\diag(B^i_1,\ldots, B^i_r),\\
&& c'=\diag(C_1,\ldots,C_r),\quad
d =\diag(D_1,\ldots,D_r),
\end{eqnarray*}
where $A^i_j, ~B^i_j ,~C_j,~D_j\in \Phi(U(n_j))$ for $j=1,\cdots,r-1$, and
$A^i_r, ~B^i_r,~C_r,~D_r \in \Phi'(U(n_r))$.

\paragraph{$i=1$}
For $j=1,\ldots,r-1$, define $V_j$ as in \eqref{eqn:VjRP}. Define
\begin{eqnarray*}
V_r& \stackrel{\Phi'}{\cong}&
\Bigl\{(A^1_r,B^1_r,\ldots, A^\ell_r,B^\ell_r,C_r)
\in U(n_r)^{2\ell+1}\mid \prod_{i=1}^\ell[A^i_r,B^i_r]=
e^{\frac{2\pi\sqrt{-1}k_r}{n_r}}I_{n_r}\bar{C_r}C_r
  \Bigr \}\\
&\cong &\tV^{\ell,1}_{n_r,-k_r}.
\end{eqnarray*}
Then $V_{\mathrm{YM}}^{\ell,1}(SO(4m))_\mu =
\prod_{j=1}^r V_j$.

\paragraph{$i=2$}
For $j=1,\ldots,r-1$, define $V_j$ as in \eqref{eqn:VjKlein}. Define
\begin{eqnarray*}
\lefteqn{V_r=
\Bigl\{(A^1_r,B^1_r,\ldots, A^\ell_r,B^\ell_r,D_r,C_r)
\in U(n_r)^{2\ell+2}\mid}\\&&\quad\quad \prod_{i=1}^\ell[A^i_r,B^i_r]=
e^{\frac{2\pi\sqrt{-1}k_r}{n_r}}I_{n_r} \bar{C_r} \bar{D_r}\bar{C}^{-1}_rD_r
  \Bigr\}\\&\cong& \tV^{\ell,2}_{n_r,-k_r}.
\end{eqnarray*}
Then
$V_{\mathrm{YM}}^{\ell,2}(SO(4m))_\mu =
\prod_{j=1}^r V_j$.

Thus $V_{\mathrm{YM}}^{\ell,i}(SO(4m))_\mu$ is also connected,
so it corresponds to a fixed topological $SO(4m)$-bundle
$P$. As in Case 1,
$$
w_2(P)=k_1+\cdots+k_r +im \quad (\mathrm{mod}\ 2).
$$
We also have a homeomorphism \eqref{eqn:nonorient-SO4m-homeo}
and a homotopy equivalence \eqref{eqn:nonorient-SO4m-homotopy}.

\paragraph{Case 3. $t_{n-1}=t_n> 0$, $n_r>1$}
$$
\mu =\sqrt{-1}\diag(\lambda_1 J_{n_1},\ldots, \lambda_r J_{n_r}),
$$
where $\lambda_1>\cdots >\lambda_r >0$.
Let $X_\mu=-2\pi\sqrt{-1}\mu$ as before. Then
$$
SO(2n)_\mu =SO(2n)_{X_\mu}\cong
\Phi(U(n_1))\times\cdots\times\Phi(U(n_r)).
$$

Let $\ep=H_{2m}$ as in Example \ref{ISO(2n)}.
Suppose that $(\ab,\ep c',X_\mu/2)\in X_{\mathrm{YM}}^{\ell,1}(SO(4m))$.
Then
$$
\exp(X_\mu/2)\ep c' \ep c'=\pab
$$
where $a_i,~b_i,~c'\in \Phi(U(n_1))\times\cdots \times  \Phi(U(n_r))$.

Let $L:\bR^{2n}\to \bC^n$ be defined as in Section
\ref{sec:SOodd_orientable}, and let
\begin{eqnarray*}
X_\mu'&=&L\circ X_\mu \circ L^{-1}\\
&=& 2\pi\sqrt{-1}\diag(\lambda_1 I_{n_1},\ldots,\lambda_r I_{n_r})
\in \mathfrak{u}(n_1)\times\cdots\times\mathfrak{u}(n_r).
\end{eqnarray*}
Then the condition on $X'_\mu$ is
$$
\exp(X'_\mu/2)\bar{c'} c'=\pab \in SU(n_1)\times\cdots\times SU(n_r),
$$
where $a_i,~b_i,~c'\in U(n_1)\times\cdots \times  U(n_r)$, and $\bar{c'}$ is the complex conjugate
of $c'$. In order that this is nonempty, we need
$1=\det(e^{\pi\sqrt{-1}\lambda_j}I_{n_j})$, i.e.,
$$
\lambda_j=\frac{2k_j}{n_j}, \quad k_j \in \bZ, \quad j=1,\ldots,r.
$$

Similarly, suppose that
$(\ab,d, \ep c',X_\mu/2)\in X_{\mathrm{YM}}^{\ell,2}(SO(4m))$. Then
\[
\exp(X_\mu/2)(\ep c') d(\ep c')^{-1}d =\pab
\]
where $a_i,~b_i,~d,~c'\in \Phi(U(n_1))\times\cdots \times  \Phi(U(n_r))$.
The condition on $X'_\mu$ is
$$
\exp(X'_\mu/2)\bar{c'}\bar{d}\bar{c'}^{-1}d \in
SU(n_1)\times \cdots \times SU(n_r),
$$
where $d,~c' \in U(n_1)\times \cdots \times U(n_r)$, and $\bar{d}$ is the complex conjugate
of $d$.
Again, we need
$$
\lambda_j=\frac{2k_j}{n_j},\quad  k_j \in \bZ, \quad j=1,\ldots,r.
$$
We conclude that for nonorientable surfaces,
$$
\mu= \sqrt{-1}\diag\Bigl(\frac{2k_1}{n_1} J_{n_1},\ldots, \frac{2k_r}{n_r} J_{n_r}\Bigr),
\quad k_j\in\bZ,\quad  \frac{k_1}{n_1}>\cdots \frac{k_r}{n_r}>0.
$$

For each $\mu$, we define the $\ep$-reduced representation varieties as in
\eqref{eqn:nonorient-SO4mRP} and \eqref{eqn:nonorient-SO4mK}
when $i=1$ and when $i=2$, respectively;
we define $V_j$ as in \eqref{eqn:VjRP} and \eqref{eqn:VjKlein}
when $i=1$ and when $i=2$, respectively.  Then
$$
V_{\mathrm{YM}}^{\ell,i}(SO(4m))_\mu =
\prod_{j=1}^r V_j.
$$
Again, $V_{\mathrm{YM}}^{\ell,i}(SO(4m))_\mu$ is connected,
so it corresponds to a fixed topological $SO(4m)$-bundle $P$,
and
$$
w_2(P)= k_1+\cdots +k_r + im\quad(\mathrm{mod}\ 2).
$$
We also have a homeomorphism \eqref{eqn:nonorient-SO4m-homeo} and a
homotopy equivalence \eqref{eqn:nonorient-SO4m-homotopy}.

\paragraph{Case 4. $t_{n-1}=t_n= 0$, $n_r>1$.}
$$
\mu =\sqrt{-1}\diag(\lambda_1 J_{n_1},\ldots, \lambda_{r-1} J_{n_{r-1}}, 0 J_{n_r}),
$$
where $\lambda_1>\cdots >\lambda_{r-1} >0$.
Let $X_\mu=-2\pi\sqrt{-1}\mu$ as before. Then
$$
SO(2n)_\mu =SO(2n)_{X_\mu}\cong
\Phi(U(n_1))\times\cdots\times\Phi(U(n_{r-1}))\times SO(2n_r).
$$
Let $\ep_\mu=\diag(H_{2m-n_r}, (-1)^{n_r}I_1, I_{2n_r-1})$.
Consider $(\ab,\ep_\mu c',X_\mu/2)\in X_{\mathrm{YM}}^{\ell,1}(SO(4m))$.
Then
$$
\exp(X_\mu/2)\ep_\mu c' \ep_\mu c'=\pab
$$
where $a_i,~b_i,~c'\in \Phi(U(n_1))\times\cdots \times  \Phi(U(n_{r-1}))\times SO(2n_r)$.

Let $L:\bR^{2(n-n_r)}\to \bC^{n-n_r}$ be defined as in Section
\ref{sec:SOodd_orientable}, and let
\begin{eqnarray*}
X_\mu'&=&L\circ \left( 2\pi\diag(\lambda_1 J_{n_1},\ldots,\lambda_{r-1} J_{n_{r-1}})\right) \circ L^{-1}\\
&=& 2\pi\sqrt{-1}\diag(\lambda_1 I_{n_1},\ldots,\lambda_{r-1} I_{n_{r-1}})\in
\mathfrak{u}(n_1)\times\cdots\times\mathfrak{u}(n_{r-1}).
\end{eqnarray*}
Then the condition on $X'_\mu$ is
$$
\exp(X'_\mu/2)\bar{c'} c'=\pab \in SU(n_1)\times\cdots\times SU(n_{r-1}),
$$
where $a_i,~b_i,~c'\in U(n_1)\times\cdots \times  U(n_{r-1})$, and $\bar{c'}$ is the complex conjugate
of $c'$. In order that this is nonempty, we need
$1=\det(e^{\pi\sqrt{-1}\lambda_j}I_{n_j})$, i.e.,
$$
\lambda_j=\frac{2k_j}{n_j}, \quad k_j \in \bZ,\quad j=1,\ldots,r-1.
$$

Similarly, suppose that
$(\ab,d, \ep_\mu c',X_\mu/2)\in X_{\mathrm{YM}}^{\ell,2}(SO(4m))$. Then
\[
\exp(X_\mu/2)(\ep_\mu c') d(\ep_\mu c')^{-1}d =\pab,
\]
where $a_i,~b_i,~d,~c'\in \Phi(U(n_1))\times\cdots \times  \Phi(U(n_{r-1}))\times SO(2n_r)$.
The condition on $X'_\mu$ is
$$
\exp(X'_\mu/2)\bar{c'}\bar{d}\bar{c'}^{-1}d \in
SU(n_1)\times \cdots \times SU(n_{r-1}),
$$
where $d,~c' \in U(n_1)\times \cdots \times U(n_{r-1})$, and $\bar{d}$ is the complex conjugate
of $d$.
Again, we need
$$
\lambda_j=\frac{2k_j}{n_j},\quad  k_j \in \bZ,\quad j=1,\ldots,r-1.
$$

We conclude that for nonorientable surfaces,
$$
\mu= \sqrt{-1}\diag\Bigl(\frac{2k_1}{n_1} J_{n_1},\ldots, \frac{2k_{r-1}}{n_{r-1}} J_{n_{r-1}},0 J_{n_r}\Bigr),
\quad k_j\in\bZ,\quad  \frac{k_1}{n_1}>\cdots >\frac{k_{r-1}}{n_{r-1}}>0.
$$

For each $\mu$, define $\ep_\mu$-reduced representation varieties
\begin{eqnarray*}
 V_{\mathrm{YM}}^{\ell,1}(SO(4m))_\mu &=&  \{ (\ab,c')\in SO(4m)_{X_\mu}^{2\ell+1}\mid\\
&&\quad \pab=\exp(X_\mu/2)\ep_\mu c' \ep_\mu c' \},\\
 V_{\mathrm{YM}}^{\ell,2}(SO(4m))_\mu &=&
\{ (\ab,d,c')\in SO(4m)_{X_\mu}^{2\ell+2}\mid \\
&& \quad \pab=\exp(X_\mu/2)\ep_\mu c' d(\ep_\mu c')^{-1}d \}.
\end{eqnarray*}

\paragraph{$i=1$}
For $j=1,\ldots,r-1$, define $V_j$ as in \eqref{eqn:VjRP}. Define
\begin{equation}
V_r=
\Bigl\{(A^1_r,B^1_r,\ldots, A^\ell_r,B^\ell_r,C_r)
\in SO(2n_r)^{2\ell+1}\mid \prod_{i=1}^\ell[A^i_r,B^i_r]= (\epsilon C_r)^2
  \Bigr \},
\end{equation}
where $\ep=\diag((-1)^{n_r}I_1, I_{2n_r-1})$, $\det(\ep)=(-1)^{n_r}$.
Let $C_0'=\ep C_0$. We see that
\begin{eqnarray*}
V_r &\cong&
\Bigl\{(A^1_r,B^1_r,\ldots, A^\ell_r,B^\ell_r,C_r')
\in SO(2n_r)^{2\ell}\times O(2n_r)\mid \\
&& \prod_{i=1}^\ell[A^i_r,B^i_r]= (C_r')^2, \det(C_r')=(-1)^{n_r} \Bigr \}\\
&\cong & V^{\ell,1}_{O(2n_r),(-1)^{n_r}}
\end{eqnarray*}
where
$V^{\ell,1}_{O(n),\pm 1}$ is the twisted representation variety defined in (\ref{eqn:Oone})
of Section \ref{sec:twistedO}. $V^{\ell,1}_{O(n),\pm 1}$ is nonempty if $\ell\geq 2$. We have shown
that $V^{\ell,1}_{O(n),\pm 1}$ is disconnected with two components
$V^{\ell,1,+ 1}_{O(n),\pm 1}$ and $V^{\ell,1,-1}_{O(n),\pm 1}$
if $\ell\geq 2$ and
$n>2$ (Proposition \ref{thm:tO}). Then
$$
V_{\mathrm{YM}}^{\ell,1}(SO(4m))_\mu =
\prod_{j=1}^r V_j.
$$

\paragraph{$i=2$}
For $j=1,\ldots,r-1$, define $V_j$ as in \eqref{eqn:VjKlein}. Define
\begin{equation}
V_r= \Bigl\{(A^1_r,B^1_r,\ldots, A^\ell_r,B^\ell_r,D_r,C_r)
\in SO(2n_r)^{2\ell+2}\mid \prod_{i=1}^\ell[A^i_r,B^i_r]=\epsilon C_r D_r
(\epsilon C_r)^{-1}D_r  \Bigr \},
\end{equation}
where $\ep=\diag((-1)^{n_r}I_1, I_{2n_r-1})$, $\det(\ep)=(-1)^{n_r}$.
Let $C_r'=\ep C_r$. We see that
\begin{eqnarray*}
V_r &\cong&
\Bigl\{(A^1_r,B^1_r,\ldots, A^\ell_r,B^\ell_r,D_r,C_r')
\in SO(2n_r)^{2\ell+1}\times O(2n_r)\mid \\
&& \prod_{i=1}^\ell[A^i_r,B^i_r]= C_r' D_r C_r'^{-1} D_r, \det(C_r')=(-1)^{n_r} \Bigr \}\\
&\cong & V^{\ell,2}_{O(2n_r),(-1)^{n_r}}
\end{eqnarray*}
where
$V^{\ell,2}_{O(n),\pm 1}$ is the twisted representation variety defined in (\ref{eqn:Otwo})
of Section \ref{sec:twistedO}. $V^{\ell,2}_{O(n),\pm 1}$ is nonempty if $\ell\geq 4$. We have shown
that $V^{\ell,2}_{O(n),\pm 1}$ is disconnected with two components
$V^{\ell,2,+1}_{O(n),\pm 1}$ and $V^{\ell,2,-1}_{O(n),\pm 1}$ if $\ell\geq 4$ and
$n>2$ (Proposition \ref{thm:tO}). Then
$$
V_{\mathrm{YM}}^{\ell,2}(SO(4m))_\mu =
\prod_{j=1}^r V_j.
$$

Thus $V_{\mathrm{YM}}^{\ell,i}(SO(4m))_\mu$ is disconnected with two connected
components if $\ell\geq 2i$ and $n_r>1$ (because $V_r$ is).
By the argument in Section \ref{sec:SO4m+2_nonorientable},
$$
V_{\mathrm{YM}}^{\ell,i}(SO(4m))_\mu^{\pm 1} =
\prod_{j=1}^{r-1}V_j \times
V^{\ell,i,\pm (-1)^{k_1+\cdots+k_{r-1}+ i \frac{(n-n_r)(n-n_r-1)}{2}} }_{O(2n_r),(-1)^{n_r}}.
$$
Note that $\ \displaystyle{
i \frac{(n-n_r)(n-n_r-1)}{2} \equiv
i(m+\frac{n_r(n_r+1)}{2})}$.

Let  $U(n_j)$ act on $V_j=\tV^{\ell,i}_{n_j,-k_j}$ by \eqref{eqn:actI} and \eqref{eqn:actII}
of Section \ref{sec:twisted} when $i=1$ and when $i=2$, respectively;
let $SO(2n_r)$ act on $V_r=V^{\ell,i}_{O(2n_r),(-1)^{n_r}}$ by
\eqref{eqn:OactI} and \eqref{eqn:OactII} in Section \ref{sec:twistedO}
when $i=1$ and $i=2$, respectively. Then we  have
a homeomorphism
$$
V_{\mathrm{YM}}^{\ell,i}(SO(4m))_\mu/SO(4m)_{X_\mu}
\cong \prod_{j=1}^{r-1} (V_j/U(n_j)) \times V_r/SO(2n_r),
$$
and a homotopy equivalence
$$
{V_{\mathrm{YM}}^{\ell,i}(SO(4m))_\mu}^{h SO(4m)_{X_\mu} }\sim
\prod_{j=1}^{r-1} {V_j}^{h U(n_j)} \times {V_r}^{h SO(2n_r)}.
$$

To simplify the notation, we write
$$
\mu=(\mu_1,\ldots,\mu_{2m})=\Bigl(\underbrace{\frac{2k_1}{n_1},\ldots,\frac{2k_1}{n_1} }_{n_1}, \ldots,
\underbrace{\frac{2k_r}{n_r},\ldots,\frac{2k_r}{n_r} }_{n_r}\Bigr)
$$
instead of
$$
\sqrt{-1}\diag\Bigl(\frac{2k_1}{n_1}J_{n_1},\ldots, \frac{2k_r}{n_r} J_{n_r}\Bigr).
$$
Let
\begin{eqnarray*}
&&\hat{I}_{SO(4m)}^{\pm 1}=\Bigl \{
\mu=\Bigl(\underbrace{\frac{2k_1}{n_1},\ldots,\frac{2k_1}{n_1}}_{n_1}, \ldots,
\underbrace{\frac{2k_{r-1}}{n_{r-1}},\ldots,\frac{2k_{r-1}}{n_{r-1}}}_{n_{r-1}} , 2k_r\Bigr)\Bigl|
n_j \in \bZ_{>0},\  k_j\in\bZ, \\
&& \quad\quad n_1+\cdots+n_{r-1}+1=n,
\frac{k_1}{n_1}>\cdots >\frac{k_{r-1}}{n_{r-1}}> |k_r|,\ (-1)^{k_1+\cdots+k_r+im}=\pm 1 \Bigr\}\\
&&\bigcup \Bigl \{
\mu=\Bigl(\underbrace{\frac{2k_1}{n_1},\ldots,\frac{2k_1}{n_1}}_{n_1}, \ldots,
\underbrace{\frac{2k_{r-1}}{n_{r-1}},\ldots,\frac{2k_{r-1}}{n_{r-1}} }_{n_{r-1}},
\underbrace{\frac{2k_r}{n_r},\ldots,\frac{2k_r}{n_r}}_{n_r-1},\pm \frac{2k_r}{n_r} \Bigr)\Bigl|\
n_j \in \bZ_{>0}, \\
&& n_r>1,\  n_1+\cdots+n_r=n,\  k_j\in\bZ,\
\frac{k_1}{n_1}>\cdots >\frac{k_r}{n_r}>0,\  (-1)^{k_1+\cdots+k_r+im}=\pm 1 \Bigr\},
\end{eqnarray*}
\begin{eqnarray*}
&&\hat{I}_{SO(4m)}^{0}=\Bigl \{
\mu=\Bigl(\underbrace{\frac{2k_1}{n_1},\ldots,\frac{2k_1}{n_1}}_{n_1}, \ldots,
\underbrace{\frac{2k_{r-1}}{n_{r-1}},\ldots,\frac{2k_{r-1}}{n_{r-1}}}_{n_{r-1}},
\underbrace{0,\ldots,0}_{n_r}\Bigr)\Bigl|\ n_j\in\bZ_{>0},  \\
&& \quad\quad n_r>1,\ n_1+\cdots+n_r=n,\  k_j\in\bZ,\
\frac{k_1}{n_1}>\cdots >\frac{k_{r-1}}{n_{r-1}}>0  \Bigr\}.
\end{eqnarray*}

\begin{pro}\label{thm:muSO4m_nonorientable}
Suppose that $\ell\geq 2i$, where $i=1,2$.
\begin{enumerate}
\item[(i)] If
$\displaystyle{
\mu=\Bigl(\underbrace{\frac{2k_1}{n_1},\ldots,\frac{2k_1}{n_1}}_{n_1}, \ldots,
\underbrace{\frac{2k_{r-1}}{n_{r-1}},\ldots,\frac{2k_{r-1}}{n_{r-1}}}_{n_{r-1}},2k_r \Bigr)
\in \hat{I}_{SO(4m)}^{\pm 1}
}$, or
$$
\mu=\Bigl(\underbrace{\frac{2k_1}{n_1},\ldots,\frac{2k_1}{n_1}}_{n_1}, \ldots,
\underbrace{\frac{2k_{r-1}}{n_{r-1}},\ldots,\frac{2k_{r-1}}{n_{r-1}}}_{n_{r-1}},
\underbrace{\frac{2k_r}{n_r},\ldots,\frac{2k_r}{n_r}}_{n_r-1},\pm \frac{2k_r}{n_r} \Bigr)
\in \hat{I}_{SO(4m)}^{\pm 1},
$$
then
$$
X_{\mathrm{YM}}^{\ell,i}(SO(4m))_\mu=X _{\mathrm{YM}}^{\ell,i}(SO(4m))_\mu^{\pm 1}
$$
is nonempty and connected. We have a homeomorphism
$$
X_{\mathrm{YM}}^{\ell,i}(SO(4m))_\mu/SO(4m) \cong \prod_{j=1}^r \tilde{\cM}^{\ell,i}_{n_j,-k_j}
$$
and a homotopy equivalence
$$
{X_{\mathrm{YM}}^{\ell,i}(SO(4m))_\mu}^{h SO(4m)} \sim
\prod_{j=1}^r \bigl(\tV^{\ell,i}_{n_j,-k_j}\bigr)^{h U(n_j)}.
$$
\item[(ii)] If
$\displaystyle{
\mu=\Bigl(\underbrace{\frac{2k_1}{n_1},\ldots,\frac{2k_1}{n_1} }_{n_1}, \ldots,
\underbrace{\frac{2k_{r-1}}{n_{r-1}},\ldots,\frac{2k_{r-1}}{n_{r-1}} }_{n_{r-1}},
\underbrace{0,\ldots,0}_{n_r}\Bigr)\in \hat{I}_{SO(4m)}^0
}$,\\
then $X_{\mathrm{YM}}^{\ell,i}(SO(4m))_\mu$ has two connected components (from both bundles
over $\Si^\ell_i$)
$$
X_{\mathrm{YM}}^{\ell,i}(SO(4m))_\mu^{+1},\quad \textup{and}\quad
X_{\mathrm{YM}}^{\ell,i}(SO(4m))_\mu^{-1}.
$$
We have homeomorphisms
$$
X_{\mathrm{YM}}^{\ell,i}(SO(4m))_\mu^{\pm 1}/SO(4m)
\cong \prod_{j=1}^{r-1} \tilde{\cM}^{\ell,i}_{n_j,-k_j}
\times \cM^{\ell,i,\pm(-1)^{k_1+\cdots+k_{r-1}+im+i\frac{n_r(n_r+1)}{2}} }_{O(2n_r),(-1)^{n_r}}
$$
and homotopy equivalences
\begin{eqnarray*}
&&\Bigl(X_{\mathrm{YM}}^{\ell,i}(SO(4m))_\mu^{\pm 1}\Bigr)^{h SO(4m)}\\
&\sim&\prod_{j=1}^{r-1} \bigl(\tV^{\ell,i}_{n_j,-k_j}\bigr)^{h U(n_j)}
\times \Bigl(\tV^{\ell,i,\pm(-1)^{k_1+\cdots+k_{r-1}+im+i\frac{n_r(n_r+1)}{2}} }_{O(2n_r),(-1)^{n_r}}\Bigr)
^{h SO(2n_r)}.
\end{eqnarray*}
\end{enumerate}
\end{pro}

\begin{pro}Suppose that $\ell\geq 2i$, where $i=1,2$. The connected components
of $X_{\mathrm{YM}}^{\ell,i}(SO(4m))^{\pm 1}$ are
$$
 \{X_{\mathrm{YM}}^{\ell,i}(SO(4m))_\mu \mid \mu\in \hat{I}^{\pm 1}_{SO(4m)}\}\cup
\{X_{\mathrm{YM}}^{\ell,i}(SO(4m))_\mu^{\pm 1}\mid \mu\in \hat{I}^0_{SO(4m)}\}.
$$
\end{pro}

Notice that, the set
$\{\mu=\sqrt{-1}\diag(\mu_1 J,\ldots,\mu_{2m}J)\mid
(\mu_1,\ldots,\mu_{2m})\in\hat{I}^{\pm 1}_{SO(4m)}\cup\hat{I}^0_{SO(4m)}\}$
is a {\em proper} subset of $\{\mu\in(\Xi^I_+)^\tau \mid I\subseteq\Delta,\tau(I)=I\}$
as mentioned in Section \ref{sec:connected}.

The following is an immediate consequence of Proposition
\ref{thm:muSO4m_nonorientable}.
\begin{pro}\label{thm:PtmuSO4m_nonorientable}
Suppose that $\ell\geq 2i$, where $i=1,2$.
\begin{enumerate}
\item[(i)] If
$\displaystyle{
\mu=\Bigl(\underbrace{\frac{2k_1}{n_1},\ldots,\frac{2k_1}{n_1}}_{n_1}, \ldots,
\underbrace{\frac{2k_{r-1}}{n_{r-1}},\ldots,\frac{2k_{r-1}}{n_{r-1} }}_{n_{r-1}},2k_r \Bigr)
\in \hat{I}_{SO(4m)}^{\pm 1}
}$, or
$$
\mu=\Bigl(\underbrace{\frac{2k_1}{n_1},\ldots,\frac{2k_1}{n_1}}_{n_1}, \ldots,
\underbrace{\frac{2k_{r-1}}{n_{r-1}},\ldots,\frac{2k_{r-1}}{n_{r-1}}}_{n_{r-1}},
\underbrace{\frac{2k_r}{n_r},\ldots,\frac{2k_r}{n_r}}_{n_r-1},\pm \frac{2k_r}{n_r} \Bigr)
\in \hat{I}_{SO(4m)}^{\pm 1},
$$
then
$$
P_t^{SO(4m)}\left(X_{\mathrm{YM}}^{\ell,i}(SO(4m))_\mu\right)
=\prod_{j=1}^r P_t^{U(n_j)}(\tV^{\ell,i}_{n_j,-k_j}).
$$
\item[(ii)] If
$\displaystyle{
\mu=\Bigl(\underbrace{\frac{2k_1}{n_1},\ldots,\frac{2k_1}{n_1} }_{n_1}, \ldots,
\underbrace{\frac{2k_{r-1}}{n_{r-1}},\ldots,\frac{2k_{r-1}}{n_{r-1}} }_{n_{r-1}},
\underbrace{0,\ldots,0}_{n_r}\Bigr)\in \hat{I}_{SO(4m)}^0
}$,
then
\begin{eqnarray*}
&& P_t^{SO(4m)}\left(X_{\mathrm{YM}}^{\ell,i}(SO(4m))_\mu^{\pm 1}\right)\\
&=&\prod_{j=1}^{r-1} P_t^{U(n_j)}(\tV^{\ell,i}_{n_j,-k_j})\times
P_t^{SO(2n_r)}\left(V_{O(2n_r),(-1)^{n_r} }^{\ell,i,\pm
(-1)^{k_1+\cdots+k_{r-1}+im+i\frac{n_r(n_r+1)}{2}}}\right).
\end{eqnarray*}
\end{enumerate}
\end{pro}

\section{Yang-Mills $Sp(n)$-Connections}
\label{sec:Sp}

$$
Sp(n)=\left\{\left(\begin{array}{cc}A&-\bar{B}\\B&\bar{A}\\ \end{array}\right)\in U(2n)~\Bigr|~
A,B\in GL(n,\bC)\right\}
$$
The maximal torus of $Sp(n)$ consists of diagonal matrices of the form
$$
\diag(u_1,\ldots,u_n, u_1^{-1},\ldots, u_n^{-1}),
$$
where $u_1,\ldots,u_n\in U(1)$. The Lie algebra of the maximal torus consists of
diagonal matrices of the form
$$
-2\pi\sqrt{-1}\diag (t_1,\ldots,t_n,-t_1,\ldots,-t_n), \quad t_i \in\bR.
$$
The fundamental Weyl chamber is
$$
\BC_0=\{ \diag(t_1,\ldots,t_n,-t_1,\ldots,-t_n) \mid t_1\geq t_2\geq\cdots\geq t_n\geq 0\}.
$$

In this section, we assume
$$
n_1,\ldots,n_r\in\bZ_{>0},\quad n_1+\cdots +n_r=n.
$$

\subsection{$Sp(n)$-connections on orientable surfaces}
\label{sec:Sp_orientable}

Any $\mu\in \BC_0$ is of the form
$$
\mu=\diag(\lambda_1 I_{n_1},\ldots,\lambda_r I_{n_r},
-\lambda_1 I_{n_1},\ldots,-\lambda_r I_{n_r}),
$$
where $\lambda_1 >\cdots >\lambda_r \geq 0$. When $\lambda_r>0$,
$Sp(n)_{X_\mu}$ consists of matrices of the form
$$
\diag(M_1,\ldots,M_r, \overline{M}_1,\ldots, \overline{M}_r),
$$
where $M_j\in U(n_j)$. When $\lambda_r=0$, $Sp(n)_{X_\mu}$
consists of matrices of the form

$$
\left(\begin{array}{ccccccccc}
 M_1  &       &       &       &         &         &     &0     \\
      &\cdot  &       &       &         &         &     &      \\
      &       & M_{r-1} &         &         &         &     &      \\
      &       &       & M_r   & 0       &         &     &-\overline{N}_r\\
      &       &       & 0     &\overline{M}_1&         &     &     \\
      &       &       &       &         &\cdot    &     &     \\
      &       &       &       &         &         &\overline{M}_{r-1}&     \\
0     &       &       & N_r   &         &         &     &\overline{M}_r\\
\end{array}\right)
$$
where $M_j\in U(n_j)$ for $j=1,\ldots,r-1$, and
$$
S=\left( \begin{array}{cc}
M_r & -\overline{N}_r\\
N_r & \overline{M}_r
\end{array} \right)\in
Sp(n_r)\subset U(2n_r).
$$
So
$$
Sp(n)_{X_\mu}\cong \left\{\begin{array}{ll}
U(n_1)\times \cdots \times  U(n_r),& \lambda_r>0,\\
U(n_1)\times \cdots\times U(n_{r-1}) \times Sp(n_r), &\lambda_r=0.
\end{array}\right.
$$

Suppose that $(\ab,X_\mu)\in X_{\mathrm{YM}}^{\ell,0}(Sp(n))$. Then
$$
\exp(X_\mu)=\pab,
$$
where $\ab\in Sp(n)_{X_\mu}$. Then we have
$$
\exp(X_\mu)\in (Sp(n)_{X_\mu})_{ss}\cong
\left\{\begin{array}{ll}
SU(n_1)\times\cdots\times SU(n_r), & \lambda_r>0,\\
SU(n_r)\times \cdots\times SU(n_{r-1})\times Sp(n_r),& \lambda_r=0.
\end{array}\right.
$$
Thus
\begin{eqnarray*}
X_\mu&=&-2\pi\sqrt{-1}\diag\Bigl(\frac{k_1}{n_1} I_{n_1},\ldots,\frac{k_r}{n_r} I_{n_r},
-\frac{k_1}{n_1} I_{n_1},\ldots,-\frac{k_r}{n_r} I_{n_r} \Bigr),\\
\mu&=& \diag\Bigl(\frac{k_1}{n_1} I_{n_1},\ldots,\frac{k_r}{n_r} I_{n_r},
-\frac{k_1}{n_1} I_{n_1},\ldots,-\frac{k_r}{n_r} I_{n_r}\Bigr),
\end{eqnarray*}
where
$$
k_j\in\bZ, \quad
\frac{k_1}{n_1}>\cdots >\frac{k_r}{n_r} \geq 0.
$$
This agrees with Section \ref{sec:SpC}.

Recall for each $\mu$, the representation variety is
\begin{eqnarray*}
\VymSp{\ell}{0}_\mu&=&\{ (\ab)\in (Sp(n)_{X_\mu})^{2\ell}\mid
\pab=\exp(X_\mu) \}.
\end{eqnarray*}

Let $i=1,\ldots,\ell$. When $k_r>0$,  write
$$
a_i=\diag\left(A^i_1,\ldots,A^i_r,\bar{A^i_1},\ldots,\bar{A^i_r}\right),\quad
b_i=\diag\left(B^i_1,\ldots,B^i_r,\bar{B^i_1},\ldots,\bar{B^i_r}\right),
$$
where $A^i_j, B^i_j\in U(n_j)$. When $k_r=0$, write
$$
a_i=
\left(\begin{array}{ccccccccc}
 A^i &       &         & 0          \\
     & A^i_r &        &-\bar{E^i_r}\\
     &       &\bar{A^i}&            \\
0    & E^i_r &        &\bar{A^i_r}\\
\end{array}\right),\quad
b_i=
\left(\begin{array}{ccccccccc}
 B^i &       &         & 0          \\
     & B^i_r &        &-\bar{F^i_r}\\
     &      &\bar{B^i}&            \\
0    & F^i_r &        &\bar{B^i_r}\\
\end{array}\right),
$$
where
\begin{eqnarray*}
&& A^i=\diag\left(A^i_1,\ldots,A^i_{r-1}\right),\quad
B^i=\diag\left(B^i_1,\ldots,B^i_{r-1}\right),\\
&& A^i_j, B^i_j \in U(n_j), \quad j=1,\ldots,r-1,\\
&& P^i=\left( \begin{array}{cc}
A^i_r & -\bar{E^i_r}\\
E^i_r & \bar{A^i_r}
\end{array}\right),
Q^i=\left( \begin{array}{cc}
B^i_r & -\bar{F^i_r}\\
F^i_r & \bar{B^i_r}
\end{array}
\right)\in Sp(n_r)\subset U(2n_r).
\end{eqnarray*}

For $j=1,\ldots,r-1$, define
\begin{equation}\label{eqn:SpVj}
\begin{aligned}
V_j&=
\Bigl\{(A^1_j,B^1_j,\ldots, A^\ell_j,B^\ell_j)
\in U(n_j)^{2\ell}\mid \prod_{i=1}^\ell[A^i_j,B^i_j]= e^{-2\pi\sqrt{-1}\frac{k_j}{n_j}}I_{n_j}
  \Bigr \}\\
&\cong X_{\mathrm{YM}}^{\ell,0}(U(n_j))_{\frac{k_j}{n_j},\ldots,\frac{k_j}{n_j}}.
\end{aligned}
\end{equation}
When $k_r>0$, define $V_r$ by \eqref{eqn:SpVj}. When $k_r=0$, define
$$
 V_r= \Bigl\{(P^1_r,Q^1_r,\ldots, P^\ell_r,Q^\ell_r)
\in Sp(n_r)^{2\ell}\mid \prod_{i=1}^\ell[P^i_r,Q^i_r]= I_{n_r}
  \Bigr \}\cong X_{\mathrm{flat}}^{\ell,0}(Sp(n_r)).
$$
Then $V_1,\ldots, V_r$ are connected, and
$\VymSp{\ell}{0}_\mu =\prod^r_{j=1}V_j$.
We have homeomorphisms
$$
\VymSp{\ell}{0}_\mu/Sp(n)_{X_\mu}\cong
\begin{cases}
\prod_{j=1}^r (V_j/U(n_j)), & k_r>0, \\
\prod_{j=1}^{r-1} (V_j/U(n_j))\times V_r/Sp(n_r), & k_r=0,
\end{cases}
$$
and homotopy equivalences
$$
{\VymSp{\ell}{0}_\mu}^{h Sp(n)_{X_\mu}}\sim
\begin{cases}
\prod_{j=1}^r {V_j}^{h U(n_j)}, & k_r>0, \\
\prod_{j=1}^{r-1} {V_j}^{h U(n_j) }\times {V_r}^{h Sp(n_r)}, & k_r=0.
\end{cases}
$$

Recall that $Sp(n)$ is simply connected, so any principal $Sp(n)$-bundle
over an orientable or nonorientable surface is trivial. For $i=0,1,2$, let
$$
\cM(\Si^\ell_i,Sp(n))=X_{\mathrm{flat}}^{\ell,i}(Sp(n))/Sp(n)
$$
be the moduli space of gauge equivalence classes of flat $Sp(n)$-connections
on $\Si^\ell_i$.

To simplify the notation, we write
$$
\mu=(\mu_1,\ldots,\mu_n)=\Bigl(\underbrace{\frac{k_1}{n_1},\ldots,\frac{k_1}{n_1} }_{n_1}, \ldots,
\underbrace{\frac{k_r}{n_r},\ldots,\frac{k_r}{n_r} }_{n_r}\Bigr)
$$
instead of
$$
\diag\Bigl(\frac{k_1}{n_1}I_{n_1},\ldots, \frac{k_r}{n_r} I_{n_r},
-\frac{k_1}{n_1}I_{n_1},\ldots,-\frac{k_r}{n_r}I_{n_r}\Bigr).
$$
Let
\begin{eqnarray*}
I_{Sp(n)}&=&\Bigl \{
\mu=\Bigl(\underbrace{\frac{k_1}{n_1},\ldots,\frac{k_1}{n_1}}_{n_1}, \ldots,
\underbrace{\frac{k_r}{n_r},\ldots,\frac{k_r}{n_r}}_{n_r}\Bigr)
 \bigl| n_j\in \bZ_{>0},\\ &&\quad n_1+\cdots+n_r=n,\  k_j\in \bZ,\
\frac{k_1}{n_1}>\cdots >\frac{k_r}{n_r}\geq 0 \Bigr\}.
\end{eqnarray*}

From the discussion above, we conclude:
\begin{pro}\label{thm:muSp_orientable}
Suppose that $\ell\geq 1$.
Let
\begin{equation}\label{eqn:muSp_orientable}
\mu=\Bigl(\underbrace{\frac{k_1}{n_1},\ldots,\frac{k_1}{n_1} }_{n_1}, \ldots,
\underbrace{\frac{k_r}{n_r},\ldots,\frac{k_r}{n_r} }_{n_r}\Bigr)\in I_{Sp(n)}.
\end{equation}
Then
$$
X_{\mathrm{YM}}^{\ell,0}(Sp(n))_\mu/Sp(n) \cong
\begin{cases}
\prod_{i=1}^r\cM(\Si^\ell_0, P^{n_j,k_j}), & k_r>0,\\
\prod_{i=1}^{r-1}\cM(\Si^\ell_0, P^{n_j,k_j}) \times \cM(\Si^\ell_0, Sp(n_r)), & k_r=0.
\end{cases}
$$
In particular, $X_{\mathrm{YM}}^{\ell,0}(Sp(n))_\mu$ is nonempty and connected. We have homotopy
equivalences
\begin{eqnarray*}
&& {X_{\mathrm{YM}}^{\ell,0}(Sp(n))_\mu}^{h Sp(n)}\\
&\sim&\begin{cases}
\prod_{i=1}^r \Bigl(X^{\ell,0}_{\mathrm{YM}}(U(n_j))
_{ \frac{k_j}{n_j},\ldots,\frac{k_j}{n_j} }\Bigr)^{h U(n_j)}, & k_r>0,\\
\prod_{i=1}^{r-1}\Bigl(X^{\ell,0}_{\mathrm{YM}}(U(n_j))
_{\frac{k_j}{n_j}, \ldots,\frac{k_j}{n_j}}\Bigr)^{h U(n_j)}
 \times X_{\mathrm{flat}}^{\ell,0}( Sp(n)) ^{h Sp(n_r)}, & k_r=0.
\end{cases}
\end{eqnarray*}
\end{pro}

\begin{pro}Suppose that $\ell\geq 1$.
The connected components of the representation variety
$X_{\mathrm{YM}}^{\ell,0}(Sp(n))$ are
$$
\{ X_{\mathrm{YM}}^{\ell,0}(Sp(n))_\mu \mid \mu\in I_{Sp(n)}\}.
$$
\end{pro}

The following is an immediate consequence of
Proposition \ref{thm:muSp_orientable}.
\begin{thm}
  Suppose that $\ell\geq 1$, and let $\mu$ be as in \eqref{eqn:muSp_orientable}.
Then
\begin{eqnarray*}
&&P_t^{Sp(n)}\left(X_{\mathrm{YM}}^{\ell,0}(Sp(n))_\mu\right)\\
&=&\left\{\begin{array}{ll}
\prod_{j=1}^r P_t^{U(n_j)}
\Bigl(X_{\mathrm{YM}}^{\ell,0}(U(n_j))_{\frac{k_j}{n_j},\ldots,\frac{k_j}{n_j}}\Bigr), & k_r>0,\\
\prod_{j=1}^{r-1} P_t^{U(n_j)}
\Bigl(X_{\mathrm{YM}}^{\ell,0}(U(n_j))_{\frac{k_j}{n_j},\ldots,\frac{k_j}{n_j}}\Bigr)\cdot
P_t^{Sp(n_r)}\left(X_{\mathrm{flat}}^{\ell,0}(Sp(n_r))\right), & k_r=0.
\end{array}\right.
\end{eqnarray*}
\end{thm}

\subsection{Equivariant Poincar\'{e} series}
\label{sec:Sp_poincare} Recall from Section \ref{sec:SpC}:
\begin{eqnarray*}
&& \Delta=\{\alpha_i=\theta_i-\theta_{i+1}\mid i=1,\ldots, n-1\}\cup \{\alpha_n=2\theta_n\}\\
&& \Delta^\vee=\{\alpha_i^\vee = e_i-e_{i+1} \mid i=1,\ldots, n-1\}\cup \{\alpha_n^\vee =e_n \}\\
&& \pi_1(H)=\bigoplus_{i=1}^n\bZ e_i,\quad
\Lambda =\bigoplus_{i=1}^{n-1} \bZ(e_i-e_{i+1})\oplus \bZ e_n,
\quad \pi_1(Sp(n))=0
\end{eqnarray*}

We will apply Theorem \ref{thm:Pt} to the case $G_\bR=Sp(n)$.
$$
\varpi_{\alpha_i}=\theta_1 +\cdots +\theta_i
$$

Case 1. $\alpha_n\in I$:
\begin{eqnarray*}
&& I= \{ \alpha_{n_1}, \alpha_{n_1+n_2},\ldots, \alpha_{n_1+\cdots + n_{r-1}},\alpha_n  \}\\
&& L^I= GL(n_1,\bC)\times \cdots \times  GL(n_r,\bC), \quad n_1+\cdots +n_r=n \\
&& \dim_\bC \fz_{L^I}-\dim_\bC \fz_{Sp(n,\bC)}=r,\quad \dim_\bC U^I = \sum_{1\leq i<j\leq r}n_i n_j +\frac{n(n+1)}{2}\\
&&\rho^I =\frac{1}{2}\sum_{i=1}^r
\biggl(n-2\sum_{j=1}^{i} n_j +n_i \biggr)
\biggl(\sum_{j=1}^{n_i}\theta_{n_1+\cdots+ n_{i-1}+j}\biggr)
+\frac{n+1}{2}(\theta_1+\cdots +\theta_n)
\end{eqnarray*}
$$
\langle \rho^I, \alpha_{n_1+\cdots+n_i}^\vee \rangle = \frac{n_i+n_{i+1}}{2}
\textup{ for } i=1,\ldots,r-1,
\quad \langle \rho^I, \alpha_n^\vee \rangle = \frac{n_r+1}{2}
$$
Case 2. $\alpha_n\notin I$:
\begin{eqnarray*}
&& I= \{ \alpha_{n_1}, \alpha_{n_1+n_2},\ldots, \alpha_{n_1+\cdots + n_{r-1}}  \}\\
&& L^I= GL(n_1,\bC)\times \cdots \times  GL(n_{r-1},\bC)\times Sp(n_r,\bC),
\quad n_1+\cdots +n_r=n \\
&& \dim_\bC \fz_{L^I}-\dim_\bC \fz_{Sp(n,\bC)}=r-1,\\
&&\dim_\bC U^I = \sum_{1\leq i<j\leq r}n_i n_j
+\frac{n(n+1)-n_r(n_r+1)}{2}\\
&& \rho^I =\frac{1}{2}\sum_{i=1}^r
\biggl(n-2\sum_{j=1}^{i} n_j +n_i \biggr)
\biggl(\sum_{j=1}^{n_i}\theta_{n_1+\cdots + n_{i-1}+j}\biggr) \\
&& \quad\quad +
\frac{n+1}{2}(\theta_1+\cdots +\theta_{n_1+\cdots+n_{r-1}})
+ \frac{n-n_r}{2}(\theta_{n_1+\cdots+n_{r-1}+1}+\cdots +\theta_n)
\end{eqnarray*}
$$
\langle \rho^I,  \alpha_{n_1+\cdots+n_i}^\vee \rangle = \frac{n_i+n_{i+1}}{2}
\textup{ for }i=1,\ldots, r-2,\quad
\langle \rho^I, \alpha_{n_1+\cdots+n_{r-1}}^\vee \rangle =\frac{n_{r-1}+1}{2}+ n_r
$$

Then we have the closed formula for the equivariant poinar\'{e} series for the
moduli space of flat $Sp(n)$-connections:

\begin{thm}\label{thm:PtSp}
\begin{eqnarray*}
&& P_t^{Sp(n)}(X_{\mathrm{flat}}^{\ell,0}(Sp(n)) ) =\\
&&\sum_{r=1}^n\sum_{\tiny \begin{array}{c}n_1,\ldots, n_r\in\bZ_{>0}\\\sum n_j=n \end{array}}
\left( (-1)^r \prod_{i=1}^r \frac{\prod_{j=1}^{n_i} (1+t^{2j-1})^{2\ell} }
{ (1-t^{2n_i})\prod_{j=1}^{n_i-1}(1-t^{2j})^2 }\right. \\
&& \cdot \frac{t^{(\ell-1)(2\sum_{i<j} n_i n_j+n(n+1))} }
{\left[\prod_{i=1}^{r-1}(1-t^{2(n_i+n_{i+1} )} )\right](1-t^{2(n_r+1)}) }
\cdot t^{2\sum_{i=1}^{r-1}(n_i +n_{i+1})+2 (n_r+1)}\\
&& + (-1)^{r-1} \prod_{i=1}^{r-1} \frac{\prod_{j=1}^{n_i} (1+t^{2j-1})^{2\ell} }
{ (1-t^{2n_i})\prod_{j=1}^{n_i-1}(1-t^{2j})^2 }
\cdot\frac{\prod_{j=1}^{n_r}(1+t^{4j-1})^{2\ell}}{\prod_{j=1}^{2n_r}(1-t^{2j})}\\
&& \left.\cdot
 \frac{t^{(\ell-1)(2\sum_{i<j} n_i n_j+n(n+1)-n_r(n_r+1))} }
{\left[\prod_{i=1}^{r-2}(1-t^{2(n_i+n_{i+1})} ) \right](1-\epsilon(r)t^{2(n_{r-1}+2n_r+1)}) }
t^{2\sum_{i=1}^{r-2}(n_i+ n_{i+1})+ 2\epsilon(r)(n_{r-1}+2n_r+1) }
\right)
\end{eqnarray*}
where
$$
\epsilon(r)=\left\{\begin{array}{ll}0 & r=1\\ 1 & r>1 \end{array}\right.
$$
\end{thm}

\begin{ex}
$$
P^{Sp(1)}_t(X_{\mathrm{flat}}^{\ell,0}(Sp(1)))
= -\frac{(1+t)^{2\ell} t^{2\ell+2} }{(1-t^2)(1-t^4)} +\frac{(1+t^3)^{2\ell}}{(1-t^2)(1-t^4)}
$$
Note that $Sp(1)=SU(2)=Spin(3)$, so
$$
P^{Sp(1)}_t(X_{\mathrm{flat}}^{\ell,0}(Sp(1)))=
P^{SU(2)}_t(X_{\mathrm{flat}}^{\ell,0}(SU(2)))=
P^{Spin(3)}_t(X_{\mathrm{flat}}^{\ell,0}(Spin(3)))
$$
as expected, where $P^{SU(2)}_t(X_{\mathrm{flat}}^{\ell,0}(SU(2)))$
is calculated in Example \ref{UtwoSUtwo}, and that
$P^{Spin(3)}_t(X_{\mathrm{flat}}^{\ell,0}(Spin(3)))$
is calculated in Example \ref{SOthree}.
\end{ex}

\begin{ex}
\begin{eqnarray*}
&& P^{Sp(2)}_t(X_{\mathrm{flat}}^{\ell,0}(Sp(2)))\\
&=& -\frac{(1+t)^{2\ell}(1+t^3)^{2\ell} t^{6\ell} }{(1-t^2)^2(1-t^4)(1-t^6)}
 +\frac{(1+t)^{4\ell} t^{8\ell}}{(1-t^2)^2(1-t^4)^2}\\
&& +\frac{(1+t^3)^{2\ell}(1+t^7)^{2\ell}}{(1-t^2)(1-t^4)(1-t^6)(1-t^8)}
- \frac{(1+t)^{2\ell}(1+t^3)^{2\ell} t^{6\ell+2}}{(1-t^2)^2 (1-t^4)(1-t^8)}
\end{eqnarray*}
Note that $Sp(2)=Spin(5)$, so
$$
P^{Sp(2)}_t(X_{\mathrm{flat}}^{\ell,0}(Sp(2)))=
P^{Spin(5)}_t(X_{\mathrm{flat}}^{\ell,0}(Spin(5)))
$$
as expected, where $P^{Spin(5)}_t(X_{\mathrm{flat}}^{\ell,0}(Spin(5)))$
is calculated in Example \ref{SOfive}.
\end{ex}

\begin{ex}
\begin{eqnarray*}
&&P^{Sp(3)}_t(X_{\mathrm{flat}}^{\ell,0}(Sp(3)))\\
&=& -\frac{(1+t)^{2\ell}(1+t^3)^{2\ell}(1+t^5)^{2\ell}t^{12\ell-4}}
{(1-t^2)^2(1-t^4)^2(1-t^6)(1-t^8)}
+\frac{(1+t)^{4\ell}(1+t^3)^{2\ell}t^{16\ell-4}}{(1-t^2)^3(1-t^4)(1-t^6)^2}\\
&&+\frac{(1+t)^{4\ell}(1+t^3)^{2\ell}t^{16\ell-6}}{(1-t^2)^3(1-t^4)^2(1-t^6)}
-\frac{(1+t)^{6\ell}t^{18\ell-6}}{(1-t^2)^3(1-t^4)^3}\\
&&+\frac{(1+t^3)^{2\ell}(1+t^7)^{2\ell}(1+t^{11})^{2\ell}}
{(1-t^2)(1-t^4)(1-t^6)(1-t^8)(1-t^{10})(1-t^{12})}\\
&&-\frac{(1+t)^{2\ell}(1+t^3)^{2\ell}(1+t^7)^{2\ell}t^{10\ell+2}}{(1-t^2)^2(1-t^4)(1-t^6)(1-t^8)(1-t^{12})}\\
&&-\frac{(1+t)^{2\ell}(1+t^3)^{4\ell}t^{14\ell-4}}{(1-t^2)^3(1-t^4)^2(1-t^{10})}
+\frac{(1+t)^{4\ell}(1+t^3)^{2\ell}t^{16\ell-4}}{(1-t^2)^3(1-t^4)^2(1-t^8)}
\end{eqnarray*}
\end{ex}

\subsection{$Sp(n)$-connections on nonorientable surfaces}
\label{sec:Sp_nonorientable}
We have $\BC_0^\tau=\BC_0$. Any $\mu\in \BC_0^\tau$ is of the form
$$
\mu=\diag(\lambda_1 I_{n_1},\ldots,\lambda_r I_{n_r},
-\lambda_1 I_{n_1},\ldots,-\lambda_r I_{n_r}),
$$
where $\lambda_1 >\cdots >\lambda_r \geq 0$.
We have
$$
Sp(n)_{X_\mu}\cong \left\{\begin{array}{ll}
U(n_1)\times\cdots\times U(n_r),& \lambda_r>0,\\
U(n_1)\times \cdots\times U(n_{r-1})\times Sp(n_r), & \lambda_r=0.
\end{array}\right.
$$

Suppose that $(\ab,\ep c',X_\mu/2)\in X_{\mathrm{YM}}^{\ell,1}(Sp(n))$, where
$$
\ep=\left(\begin{array}{rr}0& -I_n\\ I_n & 0\end{array}\right)\in Sp(n)
$$
is defined as in Example \ref{ISp(n)}. Notice that here $\ep^2\neq1$.
Then
$$
\exp(X_\mu/2)\ep c' \ep c'=\pab
$$
where $a_i,~b_i,~c'\in Sp(n)_{X_\mu}$.
Note that $\ep c'\ep c'=-\bar{c'} c'$ where $\bar{c'}$ is the
complex conjugate of $c'$, so
\begin{equation*}
\exp(X_\mu/2)(-\bar{c'}c')\in (Sp(n)_{X_\mu})_{ss}\cong
\left\{ \begin{array}{ll}
SU(n_1)\times\cdots\times SU(n_r), & \lambda_r>0,\\
SU(n_1)\times \cdots\times SU(n_{r-1})\times Sp(n_r),& \lambda_r=0.
\end{array}\right.
\end{equation*}
In order that this is nonempty, we need
$1=\det(-e^{\pi\sqrt{-1}\lambda_j}I_{n_j})$, i.e.,
$$
\lambda_j=\frac{2k_j}{n_j}-1,\quad k_j\in \bZ,\quad j=1,\ldots,r.
$$

Similarly, suppose that $(\ab,d, \ep c',X_\mu/2)\in X_{\mathrm{YM}}^{\ell,2}(Sp(n))$.
Then
$$
\exp(X_\mu/2)(\ep c')d(\ep c')^{-1} d=\pab,
$$
or equivalently,
\begin{eqnarray*}
\lefteqn{ \exp(X_\mu/2)(-\bar{c'}\bar{d}\bar{c'}^{-1}d) \in (Sp(n)_{X_\mu})_{ss}\cong}\\
&&\left\{ \begin{array}{ll}
SU(n_1)\times\cdots\times SU(n_r), & \lambda_r>0,\\
SU(n_1)\times \cdots\times SU(n_{r-1})\times Sp(n_r),& \lambda_r=0.
\end{array}\right.
\end{eqnarray*}
Again, we need
$$
\lambda_j=\frac{2k_j}{n_j}-1,\quad k_j\in \bZ,\quad j=1,\ldots,r.
$$

We conclude that for nonorientable surfaces, either
$$
\mu=\diag\Bigl( \bigl(\frac{2k_1}{n_1}-1\bigr)I_{n_1},\ldots,
\bigl(\frac{2k_r}{n_r}-1\bigr)I_{n_r},
- \bigl(\frac{2k_1}{n_1}-1\bigr)I_{n_1},\ldots,
-\bigl(\frac{2k_r}{n_r}-1\bigr)I_{n_r}\Bigr),
$$
where
$$
k_j\in\bZ,\quad \frac{k_1}{n_1}>\cdots>\frac{k_r}{n_r}> \frac{1}{2},
$$
or
\begin{eqnarray*}\lefteqn{
\mu=\diag\Bigl( \bigl(\frac{2k_1}{n_1}-1\bigr)I_{n_1},\ldots,
\bigl(\frac{2k_{r-1}}{n_{r-1}}-1\bigr)I_{n_{r-1}},0I_{n_r},}\\&&\quad\quad
- \bigl(\frac{2k_1}{n_1}-1\bigr)I_{n_1},\ldots,
-\bigl(\frac{2k_{r-1}}{n_{r-1}}-1\bigr)I_{n_{r-1}}, 0 I_{n_r}\Bigr),
\end{eqnarray*}
where
$$
k_j\in\bZ,\quad \frac{k_1}{n_1}>\cdots>\frac{k_{r-1}}{n_{r-1}}> \frac{1}{2}.
$$

Recall the for each $\mu$, the $\ep$-reduced representation varieties are
\begin{eqnarray*}
\VymSp{\ell}{1}_\mu &=& \{ (\ab,c')\in Sp(n)_{X_\mu}^{2\ell+1}\mid \\
&&\quad \pab=\exp(X_\mu/2)\epsilon c' \epsilon c' \},\\
 \VymSp{\ell}{2}_\mu &=& \{ (\ab,d,c')\in Sp(n)_{X_\mu}^{2\ell+2}\mid\\
&&\quad \pab=\exp(X_\mu/2)\ep c' d (\ep c')^{-1}d \}.
\end{eqnarray*}

Let $i= 1,\cdots,\ell$. When $\lambda_r>0$, write
\begin{eqnarray*}
&& a_i=\diag\left( A^i_1,\dots,A^i_r, \bar{A}^i_1,\ldots,\bar{A}^i_r\right),\quad
   b_i=\diag\left( B^i_1,\dots,B^i_r, \bar{B}^i_1,\ldots,\bar{B}^i_r\right),\\
&& c'=\diag\left( C_1,\dots,C_r, \bar{C}_1,\ldots,\bar{C}_r\right),\quad
   d=\diag\left( D_1,\dots,D_r, \bar{D}_1,\ldots,\bar{D}_r\right),
\end{eqnarray*}
where $A^i_j, B^i_j, C_j,B_j\in U(n_j)$. When $\lambda_r=0$, write
\begin{eqnarray*}
&&
a_i=
\left(\begin{array}{ccccccccc}
 A^i &       &           &0           \\
     & A^i_r &           &-\bar{E}^i_r\\
     &       &\bar{A}^i  &            \\
0    & E^i_r &           &\bar{A}^i_r
\end{array}\right),\quad
b_i=
\left(\begin{array}{ccccccccc}
 B^i &       &           &0           \\
     & B^i_r &           &-\bar{F}^i_r\\
     &       &\bar{A}^i  &            \\
0    & F^i_r &           &\bar{B}^i_r
\end{array}\right),\\
&& c'=
\left(\begin{array}{ccccccccc}
 C &       &           &0           \\
     & C_r &           &-\bar{H}_r\\
     &       &\bar{C}  &            \\
0    & H_r &           &\bar{C}_r
\end{array}\right),\quad
d=
\left(\begin{array}{ccccccccc}
 D &       &           &0           \\
     & D_r &           &-\bar{G}_r\\
     &       &\bar{D}  &            \\
0    & G_r &           &\bar{D}_r
\end{array}\right),
\end{eqnarray*}
where
\begin{eqnarray*}
&& A^i=\diag\left(A^i_1,\ldots,A^i_{r-1}\right),\quad
 B^i=\diag\left(B^i_1,\ldots,B^i_{r-1}\right),\\
&&  C=\diag \left(C_1,\ldots,C_{r-1}\right),\quad
 D=\diag \left(D_1,\ldots,D_{r-1}\right),\\
&& A^i_j, B^i_j, C_j, D_j\in U(n_j),\quad j=1,\ldots,r-1,\\
&& P^i=\left( \begin{array}{cc}
A^i_r & -\bar{E}^i_r\\
E^i_r & \bar{A}^i_r
\end{array}
\right),
Q^i=\left( \begin{array}{cc}
B^i_r & -\bar{F}^i_r\\
F^i_r & \bar{B}^i_r
\end{array}
\right)\in Sp(n_r)\subset U(2n_r),\\
&& S^i=\left( \begin{array}{cc}
C_r & -\bar{H}_r\\
H_r & \bar{C}_r
\end{array}
\right),
R^i=\left( \begin{array}{cc}
D_r & -\bar{G}_r\\
G_r & \bar{D}_r
\end{array}
\right)\in Sp(n_r)\subset U(2n_r).
\end{eqnarray*}

Let $\ep=\left( \begin{array}{cc}
0 & -I_{n_r}\\
I_{n_r} & 0
\end{array}
\right)\in Sp(n_r)$.
For $j=1, \ldots,r-1$, define
\begin{equation}\label{eqn:SpVjRP}
\begin{aligned}
V_j&=
\Bigl\{(A^1_j,B^1_j,\ldots, A^\ell_j,B^\ell_j,C_j)
\in U(n_j)^{2\ell+1}\mid \prod_{i=1}^\ell[A^i_j,B^i_j]= e^{-2\pi\sqrt{-1}\frac{k_j}{n_j}}
(\overline{C}_jC_j)
  \Bigr \} \\
  &\cong \tV^{\ell,1}_{n_j,k_j},
\end{aligned}\end{equation}
where $\tV^{\ell,1}_{n_j,k_j}$ is the twisted representation
variety defined in \eqref{eqn:twistI} of  Section \ref{sec:twisted}.
$\tV^{\ell,1}_{n_j,k_j}$ is nonempty if $\ell\geq 1$. We have
shown that $\tV^{\ell,1}_{n_j,k_j}$ is connected if $\ell\geq 2$
(Proposition \ref{thm:tV}).
When $\lambda_r>0$, define $V_r$ by \eqref{eqn:SpVjRP}. When $\lambda_r=0$, define
\begin{eqnarray*}
V_r&=&
\Bigl\{(P^1_r,Q^1_r,\ldots, P^\ell_r,Q^\ell_r,S_r)
\in Sp(n_r)^{2\ell+1}\mid \prod_{i=1}^\ell[P^i_r,Q^i_r]= (\epsilon S_r)^2
  \Bigr \} \\
&\stackrel{(S_r'=\epsilon S_r)}{\cong}&
\Bigl\{(P^1_r,Q^1_r,\ldots, P^\ell_r,Q^\ell_r,S_r')
\in Sp(n_r)^{2\ell+1}\mid \prod_{i=1}^\ell[P^i_r,Q^i_r]= (S_r')^2
  \Bigr \}\\
  &\cong& X_{\mathrm{flat}}^{\ell,1}(Sp(n_0)).
\end{eqnarray*}
Then $\VymSp{\ell}{1}_\mu = \prod_{j=1}^r V_j$.

Similarly, for $j=1,\ldots,r-1$, define
\begin{equation}\label{eqn:SpVjKlein}
\begin{aligned}
V_j&=
\Bigl\{(A^1_j,B^1_j,\ldots, A^\ell_j,B^\ell_j,D_j,C_j)
\in (U(n_j))^{2\ell+2}\mid \\
& \prod_{i=1}^\ell[A^i_j,B^i_j]=
\exp^{-2\pi\sqrt{-1}\frac{k_j}{n_j}}I_{n_j} \bar{C}_j \bar{D}_j\bar{C}^{-1}_jD_j
  \Bigr \}\quad\cong\quad \tV^{\ell,2}_{n_j,k_j},
\end{aligned}
\end{equation}
where $\tV^{\ell,2}_{n_j,k_j}$ is the twisted representation
variety defined in \eqref{eqn:twistII} of Section \ref{sec:twisted}.
$\tV^{\ell,2}_{n_j,k_j}$ is nonempty if $\ell\geq 1$. We have
shown that $\tV^{\ell,2}_{n_j,k_j}$ is connected if $\ell\geq 4$
(Proposition \ref{thm:tV}). When $\lambda_r>0$, define $V_r$ by
\eqref{eqn:SpVjKlein}. When $\lambda_r=0$, define
\begin{eqnarray*}
&&V_r= \Bigl\{(P^1_r,Q^1_r,\ldots, P^\ell_r,Q^\ell_r,R_r,S_r)
\in Sp(n_r)^{2\ell+2}\mid \prod_{i=1}^\ell[P^i_r,Q^i_r]=
\epsilon S_r R_r (\epsilon S_r)^{-1} R_r \Bigr \}\\
&&\stackrel{(S_r'=\epsilon S_r)}{\cong}
\Bigl\{(P^1_r,Q^1_r,\ldots, P^\ell_r,Q^\ell_r,R_r,S_r')
\in Sp(n_r)^{2\ell+2}\mid \prod_{i=1}^\ell[P^i_r,Q^i_r]=
 S_r' R_r (S_r')^{-1} R_r \Bigr \}\\
&&\quad\cong X_{\mathrm{flat}}^{\ell,2}(Sp(n_r)).
\end{eqnarray*}
Then $\VymSp{\ell}{2}_\mu =\prod_{j=1}^r V_j$.

Let $U(n_j)$ act on $V_j=\tV^{\ell,i}_{n_j,k_j}$ by \eqref{eqn:actI}
and \eqref{eqn:actII} in Section \ref{sec:twisted} when $i=1$ and
when $i=2$, respectively. Then we have homeomorphisms
$$
\VymSp{\ell}{i}_\mu/Sp(n)_{X_\mu}\cong
\begin{cases}
\prod_{j=1}^r (V_j/U(n_j)), & \lambda_r>0, \\
\prod_{j=1}^{r-1} (V_j/U(n_j))\times V_r/Sp(n_r), & \lambda_r=0,
\end{cases}
$$
and homotopy equivalences
$$
{\VymSp{\ell}{i}_\mu}^{Sp(n)_{X_\mu} }\sim
\begin{cases}
\prod_{j=1}^r {V_j}^{h U(n_j)}, & \lambda_r>0, \\
\prod_{j=1}^{r-1} {V_j}^{h U(n_j)}\times {V_r}^{h Sp(n_r)}, & \lambda_r=0.
\end{cases}
$$

To simplify the notation, we write
\begin{equation}\label{eqn:muSp_positive}
\mu=(\mu_1,\ldots,\mu_n)=\Bigl(\underbrace{\frac{2k_1}{n_1}-1,\ldots,\frac{2k_1}{n_1}-1 }_{n_1}, \ldots,
\underbrace{\frac{2k_r}{n_r}-1,\ldots,\frac{2k_r}{n_r}-1 }_{n_r}\Bigr)
\end{equation}
instead of
$$
\diag\Bigl(\bigl(\frac{2k_1}{n_1}-1\bigr)I_{n_1},\ldots,
\bigl(\frac{2k_r}{n_r}-1\bigr) I_{n_r},
-\bigl(\frac{2k_1}{n_1}-1\bigr)I_{n_1},\ldots,-\bigl(\frac{2k_r}{n_r}-1)\bigr)I_{n_r}\Bigr),
$$
and write
\begin{equation}\label{eqn:muSp_zero}
\mu=(\mu_1,\ldots,\mu_n)=\Bigl(\underbrace{\frac{2k_1}{n_1}-1,\ldots,\frac{2k_1}{n_1}-1 }_{n_1}, \ldots,
\underbrace{\frac{2k_{r-1}}{n_{r-1}}-1,\ldots,\frac{2k_{r-1}}{n_{r-1}}-1 }_{n_{r-1}},
\underbrace{0,\ldots,0}_{n_r}\Bigr)
\end{equation}
instead of
\begin{eqnarray*}\lefteqn{
\diag\Bigl(\bigl(\frac{2k_1}{n_1}-1\bigr)I_{n_1},\ldots,
\bigl(\frac{2k_{r-1}}{n_{r-1}}-1\bigr) I_{n_{r-1}}, 0I_{n_r},}\\&&\quad\quad
-\bigl(\frac{2k_1}{n_1}-1\bigr)I_{n_1},\ldots,-\bigl(\frac{2k_{r-1}}{n_{r-1}}-1)\bigr)I_{n_{r-1}},
0I_{n_r}\Bigr).
\end{eqnarray*}

Let
\begin{eqnarray*}
\hat{I}_{Sp(n)}&=&\Bigl \{
\mu=\Bigl(\underbrace{\frac{2k_1}{n_1}-1,\ldots,\frac{2k_1}{n_1}-1}_{n_1}, \ldots,
\underbrace{\frac{2k_r}{n_r}-1,\ldots,\frac{2k_r}{n_r}-1}_{n_r}\Bigr)\bigl|\ n_j\in \bZ_{>0},
\\&&
\quad  \  n_1+\cdots+n_r=n,\ k_j\in\bZ,\
\frac{k_1}{n_1}>\cdots >\frac{k_r}{n_r}>\frac{1}{2} \Bigr\}\\
& \bigcup&\Bigl \{
\mu=\Bigl(\underbrace{\frac{2k_1}{n_1}-1,\ldots,\frac{2k_1}{n_1}-1}_{n_1}, \ldots,
\underbrace{\frac{2k_{r-1}}{n_{r-1}}-1,\ldots,\frac{2k_{r-1}}{n_{r-1}}-1}_{n_{r-1}},
\underbrace{0,\ldots,0}_{n_r}\Bigr) \bigl| \\&&\quad
n_j\in \bZ_{>0}, \  n_1+\cdots+n_r=n,\ k_j\in\bZ,\
\frac{k_1}{n_1}>\cdots >\frac{k_{r-1}}{n_{r-1}}>\frac{1}{2} \Bigr\}
\end{eqnarray*}

\begin{pro}\label{thm:muSp_nonorientable}
Suppose that $\ell\geq 2i$, where $i=1,2$, and let
$\mu\in \hat{I}_{Sp(n)}$.
\begin{enumerate}
\item[(i)] If $\mu$ is of the form \eqref{eqn:muSp_positive}, then
$$
X_{\mathrm{YM}}^{\ell,i}(Sp(n))_\mu/Sp(n) \cong \prod_{j=1}^r
\tilde{\cM}^{\ell,i}_{n_j,k_j}.
$$
We have a homotopy equivalence
$$
{X_{\mathrm{YM}}^{\ell,i}(Sp(n))_\mu}^{h Sp(n)}\sim
 \prod_{j=1}^r\bigl(\tV^{\ell,i}_{n_j,k_j}\bigr)^{h U(n_j)}.
$$
\item[(ii)]  If $\mu$ is of the form \eqref{eqn:muSp_zero}, then
$$
X_{\mathrm{YM}}^{\ell,i}(Sp(n))_\mu/Sp(n) \cong \prod_{j=1}^{r-1}
\tilde{\cM}^{\ell,i}_{n_j,k_j}\times \cM(\Si^\ell_i, Sp(n_r)).
$$
We have a homotopy equivalence
$$
{X_{\mathrm{YM}}^{\ell,i}(Sp(n))_\mu}^{h Sp(n)}\sim
 \prod_{j=1}^{r-1}
\bigl(\tV^{\ell,i}_{n_j,k_j}\bigr)^{h U(n_j)}\times X_{\mathrm{flat}}^{\ell,i}( Sp(n_r))^{h Sp(n_r)}.
$$
\end{enumerate}
In particular, $X_{\mathrm{YM}}^{\ell,i}(Sp(n))_\mu$ is nonempty and connected.
\end{pro}

\begin{pro}Suppose that $\ell\geq 2i$, where $i=1,2$.
The connected components of
$X_{\mathrm{YM}}^{\ell,i}(Sp(n))$ are
$$
\{ X_{\mathrm{YM}}^{\ell,i}(Sp(n))_\mu \mid \mu\in \hat{I}_{Sp(n)}\}.
$$
\end{pro}

Notice that, the set
$\{\mu=\diag(\mu_1,\ldots,\mu_n,-\mu_1,\ldots,-\mu_n)|(\mu_1,\ldots,\mu_n)\in\hat{I}_{Sp(n)}\}
$ is a {\em proper} subset of $\{\mu\in(\Xi^I_+)^\tau|I\subseteq\Delta,\tau(I)=I\}$
as mentioned in Section \ref{sec:connected}.

The following is an immediate consequence of
Proposition \ref{thm:muSp_nonorientable}.
\begin{thm}\label{thm:PtmuSp_nonorientable}
 Suppose that $\ell\geq 2i$, where $i=1,2$, and let $\mu\in \hat{I}_{Sp(n)}$.
\begin{enumerate}
\item[(i)] If $\mu$ is of the form \eqref{eqn:muSp_positive}, then
$$
P_t^{Sp(n)}\left(X_{\mathrm{YM}}^{\ell,i}(Sp(n))_\mu\right)
= \prod_{j=1}^r P_t^{U(n_j)}(\tV^{\ell,i}_{n_j,k_j}).
$$
\item[(ii)] If $\mu$ is of the form \eqref{eqn:muSp_zero}, then
$$
P_t^{Sp(n)}\left(X_{\mathrm{YM}}^{\ell,i}(Sp(n))_\mu\right)
= \prod_{j=1}^{r-1} P_t^{U(n_j)}(\tV^{\ell,i}_{n_j,k_j})\cdot
P_t^{Sp(n_r)}\left(X_{\mathrm{flat}}^{\ell,i}(Sp(n_r))\right).
$$
\end{enumerate}
\end{thm}

\begin{appendix}
\section{Remarks on Laumon-Rapoport Formula} \label{sec:LR}

In this appendix, we explain how to use the argument in \cite{LR}
to obtain Theorem \ref{thm:Pt}, which is a slightly modified version of
\cite[Theorem 3.4]{LR}. We work over $\bC$.

\subsection{Notation}
The following is  a correspondence between the notation in \cite{FM}
(which we followed closely in Section \ref{sec:algebraic_geometry})
and that in \cite{LR}.
\\

\begin{center}
\begin{tabular}{||c|c|c||}\hline
                           & \cite{LR}    & \cite{FM} \\ \hline
\begin{tabular}{c}
minimal parabolic \\
subgroup (Borel)
\end{tabular}               & $P_0$        &         \\ \hline
Cartan of $G$               & $M_0$        & $H$       \\ \hline
parabolic subgroup          & $P=M_P N_P$  & $P=LU$    \\ \hline
Levi subgroup               & $M_P$        & $L$       \\ \hline
unipotent radical           & $N_P$        & $U$       \\ \hline
\begin{tabular}{c}
center of the \\
Levi subgroup
\end{tabular}               & $Z_P$        & $Z(L)$    \\ \hline
\begin{tabular}{c}
connected center \\
of $M_P$
\end{tabular}               & $A_P$        & $Z(L)_0$  \\ \hline
 \begin{tabular}{c} \\ \end{tabular}
                            & $A_P'\subset M_{P,\mathrm{ab}}$ & $L/[L,L]=Z(L)_0/Z(L)_0\cap[L,L]$ \\ \hline
                            & $X_*(A_P)$   & $\pi_1(Z(L)_0)$ \\ \hline
                   & $X_*(A_P')$  & $\pi_1(H)/\hat{\Lambda}_L =\pi_1(L/[L,L])$ \\ \hline
                   & $X_*(A_{P_0}')$ & $\pi_1(H)$\\ \hline
                   & $\fa_0=\fa_{P_0}$ & $\fh_\bR$\\ \hline
                   & \begin{tabular}{c}
                   $\fa_P =\bR \otimes X_*(A_P)$\\
                   $\quad =\bR\otimes X_*(A_P')$ \end{tabular}  & $(\fz_L)_\bR$ \\ \hline
                   & \begin{tabular}{c}
                   $\fa_G =\bR \otimes X_*(A_G)$\\
                   $\quad =\bR\otimes X_*(A_G')$ \end{tabular} & $(\fz_G)_\bR \cong \fh_\bR/V^*$ \\ \hline
& $\fa^G_0 =\fa^G_{P_0}\subset \fa_0$  &  $V^* =\Lambda\otimes \bR  \subset \fh_\bR$ \\ \hline
root system           & $\Phi_0=\Phi_{P_0} \subset \fa_0^\vee$   & $R\subset \fh_\bR^*$ \\ \hline
set of positive roots & $\Phi_0^+=\Phi_{P_0}^+\subset \Phi_0$    & $R^+\subset R$ \\ \hline
set of simple roots   & $\Delta_0 =\Delta_{P_0}\subset \Phi_0^+$ & $\Delta\subset R^+$ \\ \hline
coroot lattice of $G$ & $\bigoplus_{\alpha\in \Delta_0}\bZ \alpha^\vee$ &
$\Lambda = \bigoplus_{\alpha\in \Delta}\bZ \alpha^\vee \subset \pi_1(H)$ \\ \hline
\end{tabular}\end{center}

In this appendix, we will closely follow the notation in \cite{LR}. We will not repeat
most of the definitions in \cite{LR}.

Following \cite{LR}, if $P\subset Q\subset R$
are three parabolic subgroups of $G$, there are canonical splittings
$\fa_P=\fa^Q_P\oplus\fa^R_Q\oplus\fa_R$ and
$\fa_P^* =\fa_P^{Q*}\oplus \fa_Q^{R*} \oplus \fa_R^*$.
Given $H\in\fa_P$, we denote by $[H]^Q$, $[H]_Q^R$, and $[H]_R$ the canonical projections
of $H$ onto $\fa^Q_P$, $\fa^R_Q$, and $\fa_R$, respectively. The components of
$\beta\in \fa_P^*$ in $\fa_P^{Q*}$, $\fa_Q^{R*}$, and $\fa_R^*$ are $\beta|_{\fa^Q_P}$,
$\beta|_{\fa^R_Q}$, and $\beta|_{\fa_R}$, respectively.
Given $\alpha\in \Delta_P=\Delta_P^G\subset \fa_P^{G*}$, let $\tal$ denote
the unique element in $\Delta_0\subset \fa_{P_0}^{G*}$ such that $\tal|_{\fa_P^G}=\alpha$.
Then $\tal^\vee\in \fa_{P_0}^G$ and
$\alpha^\vee = [\tal^\vee]_P \in \fa_P^G$.
The subset $I^P$ of the set of simple roots in \cite{FM} corresponds to $\Delta_P=\Delta_P^G$
in the following way:
\begin{eqnarray*}
I^P&=& \{ \tilde{\alpha}\mid \alpha\in \Delta_P \} \subset \Delta_0 \subset \fa_{P_0}^{G*}\\
\Delta_P &=& \{\beta|_{\fa_P}\mid \beta\in I^P\}\subset \fa_P^{G*}\\
\Delta_{P_0}^P &=& \{\beta|_{\fa_{P_0}^P}\mid \beta \in \Delta_0\setminus I^P\} \subset \fa_{P_0}^{P*}\\
\Delta_{P}^{Q} &=& \{\beta|_{\fa_{P}^Q}\mid \beta \in I^P\setminus I^Q\}\subset \fa_{P}^{Q*}
\end{eqnarray*}

We continue the table of correspondence between notations in \cite{LR} and \cite{FM}:
\\
\begin{center}
\begin{tabular}{||c|c|c||}\hline
                     & \cite{LM} & \cite{FM} \\ \hline
                     & $X_*(A_P')$ & $\pi_1(H)/\hat{\Lambda}_L$ \\ \hline
                     & $\Lambda^G_P =X_*(A_P')/\bigoplus_{\alpha \in \Delta^G_P}\bZ \alpha^\vee$
                     & $\pi_1(H) / (\hat{\Lambda}_L\oplus  \bigoplus_{\alpha\in I^P} \bZ\alpha^\vee)$ \\ \hline
\begin{tabular}{c}
Topological type \\
 of $G$-bundle
\end{tabular}
                     & $\Lambda^G_{P_0}$  & $\pi_1(H)/\Lambda =\pi_1(G)$\\ \hline
\begin{tabular}{c}
Topological type \\
 of $M_P$-bundle
\end{tabular}
                     & $\Lambda^P_{P_0}$ & $\pi_1(H)/\Lambda_L =\pi_1(L)$ \\ \hline
\end{tabular}\end{center}


Given a parabolic subgroup $P$ of $G$, the topological
type of an $M_P$-bundle is given by $\nu_P\in \Lambda^P_{P_0}\cong \pi_1(M_P)$.
The slope of an $M_P$-bundle is given by $\nu_P'\in X_*(A_P')$.
The commutative diagram in Section \ref{sec:ABpts} can be rewritten as follows:
$$
\begin{CD}
@. 0     @.  0\\
@. @VVV @VVV\\
0 @>>> \hLambda_P/\Lambda_P @>{j_{ss}}>> \hLambda/\Lambda
@>{ \oplus_{\alpha\in \Delta^G_P}\varpi_\tal }>> \oplus_{\alpha\in \Delta^G_P}\bQ/\bZ \\
@.  @VV{i_P}V @VV{i_G}V @|\\
@.  \Lambda^P_{P_0} @>{[\cdot]_G}>> \Lambda^G_{P_0}
@>{\oplus_{\alpha\in \Delta^G_P}\varpi_\tal}>> \oplus_{\alpha\in \Delta^G_P}\bQ/\bZ \\
@. @VV{p_P}V @VV{p_G}V\\
@. X_*(A_P') @>{[\cdot]_G}>>  X_*(A_G') \\
@. @VVV @VVV\\
@. 0 @. 0 \\
\end{CD}
$$

Here $\Lambda_P =\oplus_{\alpha\in \Delta_{P_0}^P} \bZ\alpha^\vee\subset
 X_*(A_{P_0}')$, and $\hLambda_P$  is the saturation of $\Lambda_P$ in $X_*(A'_{P_0})$.
Let $\nu_P'$ and $\nu_G'$ denote the projections $p_P(\nu_P)$ and $p_G(\nu_G)$, respectively.

Recall that $\{\varpi_\alpha \mid \alpha\in \Delta_0\}$ is a basis of
the real vector space $\fa_{P_0}^{G*}$ which is dual to the basis
$\{ \alpha^\vee \mid \alpha\in\Delta_0\}$ of $\fa_{P_0}^G$. Given $\alpha\in \Delta_0$,
we extend $\varpi_\alpha:\fa_{P_0}^G\to \bR$ to $\varpi_\alpha:\fa_0=\fa_{P_0}^G \oplus \fa_G \to \bR$ by
zero on $\fa_G$. Then $\varpi_\alpha$ takes integral values on
$\oplus_{\alpha\in \Delta_0} \bZ \alpha^\vee\subset \fa_{P_0}^G \subset \fa_0$,
and takes rational values on $X_*(A_{P_0}') \subset \fa_0$. So it induces a map
$$
\varpi_\alpha: \Lambda_{P_0}^Q=X_*(A_{P_0}')\left/\bigoplus_{\alpha\in \Delta_{P_0}^Q}\bZ\alpha^\vee\right.
\to \bQ/\bZ
$$
where $Q$ is any parabolic subgroup of $G$. More explicitly,
given $\nu_Q \in \Lambda_{P_0}^Q$, let $X\in X_*(A'_{P_0})$ be a representative
of $\nu_Q$. Then $\varpi_\alpha(\nu_Q)= \varpi_\alpha(X)+\bZ$.

\subsection{Inversion formulas}
Let $A$ be a fixed topological abelian group.
In \cite{LR}, Laumon and Rapoport introduced the notion of $\widehat\Gamma$-converging
functions and $\Gamma$-converging functions from $\fP$ to $A$, where
$$
\fP=\{ (P,\nu_P')\mid P\in \cP, \nu'_P \in X_*(A_P') \}.
$$
We will introduce similar notion for
functions from $\fT$ to $A$, where
$$
\fT=\{ (P,\nu_P)\mid P\in \cP, \nu_P \in \Lambda^P_{P_0} \}.
$$

\begin{df} Let
$\fT=\{ (P,\nu_P)\mid P\in \cP, \nu_P \in \Lambda^P_{P_0} \}$,
and let $A$ be a fixed topological abelian group.
A function $a: \fT\to A$ is
{\em $\widehat\Gamma$--converging} if for each standard parabolic subgroup
$P\subset Q$ of $G$ and each $\nu_Q\in \Lambda^Q_{P_0}$, the finite sum
$$
\sum_{\scriptstyle \nu_P\in \Lambda^P_{P_0}\atop\scriptstyle [\nu_P]_Q=\nu_Q}
\widehat\Gamma_P^Q([\nu_P']^Q,T)a(P,\nu_P)
$$
admits a limit as $T\in \fa_P^{Q+}$ goes to infinity. If this is the case,
we shall denote this limit by
$$
\sum_{\scriptstyle \nu_P\in \Lambda^P_{P_0}\atop\scriptstyle [\nu_P]_Q=\nu_Q}
\widehat\tau_P^Q([\nu_P']^Q)a(P,\nu_P)\ .
$$
A function $b: \fT \to A$ is {\em $\Gamma$--converging} if for each
standard parabolic subgroup $P\subset Q$ of $G$ and
each $\nu_Q\in \Gamma^Q_{P_0}$, the finite sum
$$
\sum_{\scriptstyle \nu_P\in \Lambda^P_{P_0}\atop\scriptstyle [\nu_P]_Q=\nu_Q}
\Gamma_P^Q([\nu_P']^Q,T)b(P,\nu_P)
$$
admits a limit as $T\in \fa_P^{Q+}$ goes to infinity.
If this is the case, we shall denote this limit by
$$
\sum_{\scriptstyle \nu_P\in \Lambda^P_{P_0}\atop\scriptstyle [\nu_P]_Q=\nu_Q}
\tau_P^Q([\nu_P']^Q)b(P,\nu_P)\ .
$$
\end{df}

The following inversion formula is an analogue of
\cite[Theorem 2.1]{LR}.
\begin{thm} \label{thm:inversionI}
For each $\widehat\Gamma$--converging function
$a:\fT\rightarrow A$, there exists a
unique $\Gamma$--converging function $b:\fT\rightarrow A$ such that,
for each $(Q,\nu_Q)\in \fT$, we have
$$
a(Q,\nu_Q)=\sum_{\scriptstyle P\in \cP\atop\scriptstyle P\subset Q}
\sum_{\scriptstyle \nu_P\in \Lambda^P_{P_0}\atop\scriptstyle [\nu_P]_Q=\nu_Q}
\tau_P^Q([\nu_P']^Q)b(P,\nu_P)\ .
$$

The function $b$ is given by the following formula~:
for each $(Q,\nu_Q)\in \fT$, we have
$$
b(Q,\nu_Q)=
\sum_{\scriptstyle P\in \cP\atop
\scriptstyle P\subset Q}(-1)^{{\rm dim}(\fa_P^Q)}
\sum_{\scriptstyle \nu_P\in \Lambda^P_{P_0}\atop\scriptstyle [\nu_P]_Q=\nu_Q}
\widehat\tau_P^Q([\nu_P']^Q)a(P,\nu_P)\ .
$$
\end{thm}

Theorem \ref{thm:inversionI} is an easy consequence of the following
two lemmas:

\begin{lm}[Langlands] \label{thm:langlands}
For any standard parabolic subgroups $P\subset R$ of $G$ and any
$H\in \fa^R_P$, we have
\begin{equation}\label{eqn:QR}
\sum_{P\subset Q\subset R}(-1)^{\dim(\fa^R_Q)}
\tau^Q_P([H]^Q)\widehat\tau^R_Q([H]_Q)=\delta^R_P
\end{equation}
and
\begin{equation}\label{eqn:PQ}
\sum_{P\subset Q\subset R} (-1)^{\dim(\fa^Q_P)}
\widehat\tau^Q_P([H]^Q)\tau^R_Q([H]_Q)=\delta^R_P.
\end{equation}
\end{lm}

\begin{lm}[Arthur] \label{thm:arthur}
If $T\in \fa^{R+}_P\subset\, ^+\fa^R_P$, the function
$H\mapsto \Gamma^R_P(H,T)$
$($resp. $H \mapsto \hat{\Gamma}^R_P(H,T)\,)$ is the characteristic
function of the bounded subset
$$
\{ H\in \fa^R_P \mid \langle \alpha, H\rangle >0,
\langle \varpi_\alpha, H\rangle
\leq \langle \varpi_\alpha, T \rangle,
\forall \alpha\in \Delta^R_P \} \subset \fa^{R+}_P
$$
$($resp.
$$
\{ H\in \fa^R_P \mid \langle \varpi_\alpha^R, H\rangle >0,
\langle \alpha, H\rangle
\leq \langle \alpha, T \rangle,
\forall \alpha\in \Delta^R_P \} \subset\, ^+\fa^R_P \, )
$$
of $\fa^R_P$.
\end{lm}

\begin{proof}[Proof of Theorem \ref{thm:inversionI}]
\begin{eqnarray*}
& & \sum_{\scriptstyle Q\in \cP\atop\scriptstyle Q\subset R}
\sum_{\scriptstyle \nu_Q\in \Lambda^Q_{P_0}\atop\scriptstyle [\nu_Q]_R=\nu_R}
\tau_Q^R([\nu_Q']^R)b(Q,\nu_Q)\\
&=& \sum_{\scriptstyle Q\in \cP\atop\scriptstyle Q\subset R}
\sum_{\scriptstyle \nu_Q\in \Lambda^Q_{P_0}\atop\scriptstyle [\nu_Q]_R=\nu_R}
\tau_Q^R([\nu_Q']^R)
\sum_{\scriptstyle P\in \cP\atop
\scriptstyle P\subset Q}(-1)^{{\rm dim}(\fa_P^Q)}
\sum_{\scriptstyle \nu_P\in \Lambda^P_{P_0}\atop\scriptstyle [\nu_P]_Q=\nu_Q}
\widehat\tau_P^Q([\nu_P']^Q)a(P,\nu_P)\\
&=& \sum_{\scriptstyle P\subset Q \subset R}
\sum_{\scriptstyle \nu_P\in \Lambda^P_{P_0}\atop\scriptstyle [\nu_P]_R=\nu_R}
(-1)^{\dim(\fa_P^Q)} \widehat\tau_P^Q([\nu_P']^Q)\tau_Q^R([\nu_P']_Q^R)
a(P,\nu_P)
\end{eqnarray*}
For fixed $\nu_P$, we have
\begin{eqnarray*}
\lefteqn{
\sum_{\scriptstyle P\subset Q \subset R}
 (-1)^{\dim(\fa_P^Q)} \widehat\tau_P^Q([\nu_P']^Q)\tau_Q^R([\nu_P']_Q^R)=}\\&&
 \sum_{\scriptstyle P\subset Q \subset R}
 (-1)^{\dim(\fa_P^Q)} \widehat\tau_P^Q([[\nu_P']^R]^Q)\tau_Q^R([[\nu_P']^R]_Q) = \delta_P^R
\end{eqnarray*}
where the last equality follows from \eqref{eqn:PQ} in Lemma \ref{thm:langlands}. So
$$\sum_{\scriptstyle Q\in \cP\atop\scriptstyle Q\subset R}
\sum_{\scriptstyle \nu_Q\in \Lambda^Q_{P_0}\atop\scriptstyle [\nu_Q]_R=\nu_R}
\tau_Q^R([\nu_Q']^R)b(Q,\nu_Q)\\
=  \sum_{\scriptstyle P\in \cP\atop\scriptstyle P\subset R}
\sum_{\scriptstyle \nu_P\in \Lambda^P_{P_0}\atop\scriptstyle [\nu_P]_R=\nu_R}
\delta_P^R a(P,\nu_P)
=a(R,\nu_R)
$$
\end{proof}

Now we consider a special case of Theorem \ref{thm:inversionI}.
For any $P\in \cP$, fix $n_P\in \bZ_{\geq 0}$ and
$\ep^P_0 \in \fa_0^{P*}\subset \fa_0^*$ such that
for any standard parabolic subgroups $P\subset Q$ of $G$,
$$
n_P \geq n_Q,\quad
(\ep^Q_0 -\ep^P_0)\bigr|_{\fa_0^P}=0,
\quad \langle \ep^Q_P,\alpha^\vee\rangle\in \bZ_{>0}
\quad \forall \alpha\in \Delta_P^Q,
$$
where $\ep_P^Q=(\ep^Q_0-\ep^P_0)\bigr|_{\fa^Q_P}$.
(Here we use $\ep_P^Q$ instead of $\delta_P^Q$, which is used
in \cite{LR}, to avoid confusion with the $\delta_P^R$ in
Lemma \ref{thm:langlands}.)

We have the following analogue of  \cite[Lemma 2.3]{LR}:
\begin{lm}
For each $(Q,\nu_Q)\in \fT$ and each standard parabolic subgroup
$P\subset Q$
of $G$, we have
$$\displaylines{
\qquad\sum_{\scriptstyle \nu_P\in \Lambda^P_{P_0}\atop
\scriptstyle [\nu_P]_Q=\nu_Q}
\widehat\tau_P^Q([\nu_P']^Q)t^{\langle\ep_P^Q,[\nu_P']^Q\rangle}
=\Bigl(\prod_{\alpha\in\Delta_P^Q}
{1\over 1-t^{\langle\ep_P^Q,\alpha^\vee\rangle}}\Bigl)
t^{\sum_{\alpha\in\Delta_P^Q}\langle\ep_P^Q,\alpha^\vee\rangle
\langle\varpi_\tal (\nu_Q)\rangle}\,,\qquad}
$$
where, for each $\mu \in \bR/\bZ$, $\langle \mu \rangle\in \bR$
is the unique representative of the class $\mu$ such that
$0<\langle \mu \rangle\leq 1$.
\end{lm}

Notice that $<\cdot,\cdot>$ denotes the pairing between dual spaces, while
$<\cdot>$ denotes the unique representative in $(0,1]$ of the class $\cdot\in \bR/\bZ$.

\begin{proof}
Given
$$
\nu_Q\in \Lambda_{P_0}^Q = X_*(A_{P_0}')\left/\bigoplus_{\alpha\in\Delta_{P_0}^Q} \bZ\alpha^\vee\right.,
$$
we choose a representative $X_0 \in X_*(A_{P_0}')\subset \fa_0$ of $\nu_Q$.
Let
$$
\tilde{S}=\Bigl\{ X_0 + \sum_{\alpha\in \Delta^Q_P} m_\alpha \tal^\vee
\, \Bigl|\, m_\alpha\in  \bZ \Bigr. \Bigr\} \subset  X_*(A_{P_0}') .
$$
Then the natural projection
$$
X_*(A_{P_0}')\to  \Lambda_{P_0}^P
= X_*(A_{P_0}')\left/\bigoplus_{\alpha\in\Delta_{P_0}^P} \bZ\alpha^\vee\right.
$$
restricts to a bijection
$$
j: \tilde{S} \stackrel{\cong}{\longrightarrow}
S=\{ \nu_P\in \Lambda_{P_0}^P \mid [\nu_P]_Q =\nu_Q \}.
$$
Let
$$
\displaylines{
f(t)=\sum_{\scriptstyle \nu_P\in \Lambda^P_{P_0}\atop
\scriptstyle [\nu_P]_Q=\nu_Q}
\widehat\tau_P^Q([\nu_P']^Q)t^{\langle\ep_P^Q,[\nu_P']^Q\rangle}.}
$$
Then
$$
f(t)=\sum_{\nu_P\in S}\widehat\tau_P^Q([\nu_P']^Q)t^{\langle\ep_P^Q,[\nu_P']^Q\rangle}
=\sum_{\nu_P\in S_+} t^{\langle\ep_P^Q,[\nu_P']^Q\rangle},
$$
where
$$
S_+=\{ \nu_P \in S \mid
\langle \varpi^Q_\alpha, [\nu_P']^Q \rangle >0\ \forall \alpha \in \Delta_P^Q \}.
$$
Let $\tilde{S}_+=j^{-1}(S_+)$. Then
$$
\tilde{S}_+ =\Bigl\{ X_0 + \sum_{\alpha\in \Delta_P^Q} m_\alpha \tal^\vee
\,\Bigl|\, m_\alpha \in \bZ, \varpi_\tal(X_0)+ m_\alpha >0
\ \forall\alpha\in \Delta_P^Q \bigr\}.
$$
So
\begin{eqnarray*}
f(t)&=& \sum_{\nu_P\in S_+} t^{\langle \ep_P^Q, [\nu_P']^Q \rangle }
=\sum_{X \in \tilde{S}_+}t^{\langle \ep_P^Q, [X]_P^Q \rangle }
= \sum_{X\in \tilde{S}_+} \prod_{\alpha\in \Delta_P^Q}
t^{\langle \ep_P^Q,\alpha^\vee\rangle \langle \varpi^Q_\alpha,[X]_P^Q\rangle}\\
&=& \prod_{\alpha\in \Delta_P^Q}
\sum_{\scriptstyle m_\alpha\in \bZ\atop
\scriptstyle \varpi_\tal(X_0)+m_\alpha>0}
t^{\langle \ep_P^Q,\alpha^\vee\rangle (\varpi_\tal(X_0)+m_\alpha)}
\end{eqnarray*}

Note that $\langle \varpi_\tal, \nu_Q\rangle=\varpi_\tal(\nu_Q)
=\varpi_\tal(X_0) + \bZ \in \bQ/\bZ$
for all $\alpha\in \Delta^Q_P$.
As in the proof of \cite[Lemma 2.3]{LR}, for $p\in \bZ_{>0}$
and $x\in \bR$, we have
$$
\sum_{\scriptstyle m\in \bZ\atop \scriptstyle x+m>0}t^{p(x+m)}
=\frac{t^{p\langle x+\bZ \rangle}}{1-t^p}.
$$
Thus,
$$
f(t)= \prod_{\alpha\in \Delta_P^Q}
\frac{t^{\langle \ep_P^Q,\alpha^\vee\rangle\langle \varpi_\tal(\nu_Q)\rangle} }
{1-t^{\langle \ep_P^Q,\alpha^\vee\rangle} }
$$
\end{proof}

Set $m(P,\nu_P')= n_P + \langle \ep^G_P,\nu_P'\rangle$. We have now concluded with
the following inversion formula, which
is a slightly modified version of \cite[Theorem 2.4]{LR}.
\begin{thm}\label{thm:inversionII}
Given $a_0:\cP\to A$, there exists a unique function $b_0:\fT\rightarrow A$
which satisfies the relation
$$
a_0(Q)=\sum_{\scriptstyle P\in \cP \atop\scriptstyle P\subset Q}
\sum_{\scriptstyle \nu_P\in \Lambda^P_{P_0}\atop\scriptstyle [\nu_P]_Q=\nu_Q}
\tau_P^Q([\nu_P']^Q)b_0(P,\nu_P)t^{m(P,\nu_P')-m(Q,\nu_Q')}\,,
$$
for each $(Q,\nu_Q)\in\fT$. This function is given by
\begin{eqnarray*}\lefteqn{
\quad b_0(Q,\nu_Q)=}\\&&
\sum_{\scriptstyle P\in \cP \atop
\scriptstyle P\subset Q}
(-1)^{\dim(\fa_P^Q)}a_0(P)t^{n_P-n_Q}
\Bigl(\prod_{\alpha\in\Delta_P^Q}
{1\over 1-t^{\langle\ep_P^Q,\alpha^\vee\rangle}}\Bigl)\cdot
t^{\sum_{\alpha\in\Delta_P^Q}\langle\ep_P^Q,\alpha^\vee\rangle
\langle\varpi_\tal (\nu_Q)\rangle}
\in A\,,
\end{eqnarray*}
for each $(Q,\nu_Q)\in\fT$.
\end{thm}

\subsection{Inversion of the Atiyah-Bott recursion relation}

Let $\cC(G,\nu_G)$ be the space of complex structures
on a $C^\infty$ principal $G$-bundle over a Riemann
surface of genus $g\geq 2$ with topological type $\nu_G \in \Lambda^G_{P_0} \cong\pi_1(G)$.
Let $\cC^{ss}(G,\nu_G)\subset \cC(G,\nu_G)$ be the semi-stable stratum.
Let $P_t(G,\nu_G)$ and $P_t^{ss}(G,\nu_G)$ be the $\cG$-equivariant
Poincar\'{e} series of $\cC(G,\nu_G)$ and $\cC^{ss}(G,\nu_G)$, respectively.
Let $\cC(G,P,\nu_P)\subset \cC(G,\nu_G)$ be the stratum which corresponds to
$(P,\nu_P)\in \fT$, where $[\nu_P]_G =\nu_G$. Then the real codimension $m(P,\nu_P')$
of the stratum
$\cC(G,P,\nu_P)$ is equal to
$$
 2\dim(N_P)(g-1) + 4 \langle\rho_P^G,\nu_P'\rangle,
$$
where $N_P$ is the unipotent radical of $P$ and
$$
\rho_P^G=\frac{1}{2}\sum_{\alpha\in \Phi_P^{G+} } \alpha\in \fa^{G*}_P
\subset \fa_P^*.
$$
Clearly $m(G,\nu_G') =0$.

With the above notation, the Atiyah-Bott recursion relation can be stated as follows:
\begin{thm}[Atiyah-Bott] \label{thm:ABrecursion}
The stratification of $\cM(G,\nu_G)$ by the $\cM(G,P,\nu_P)$
is perfect modulo torsion, so that for the Poincar\'{e} series, we have
\begin{equation}\label{eqn:ab}
P_t(G,\nu_G)=\sum_{P\in \cP}
\sum_{\scriptstyle \nu_P\in \Lambda^P_{P_0}\atop
\scriptstyle [\nu_P]_G=\nu_G}
\tau_P^G([\nu_P']^G)t^{m(P,\nu_P')}P^{ss}_t(M_P,\nu_P).
\end{equation}
\end{thm}

Note that Theorem \ref{thm:ABrecursion} and
\cite[Theorem 3.2]{LR} are slightly different
when $G_{ss}$ is not simply connected.

\begin{thm}[{\cite[Theorem 3.3]{LR} }]\label{thm:LRgeneralG}
For any $\nu_G\in \Lambda^G_{P_0}$, we have
$$
P_t(G,\nu_G)=\Bigl({(1+t)^{2g}\over 1-t^2}\Bigr)^{{\rm dim}(\fa_G)}
\prod_{i=1}^{{\rm dim}(\fa_0^G)}
{(1+t^{2d_i(G)-1})^{2g}\over (1-t^{2d_i(G)-2})(1-t^{2d_i(G)})}\ .
$$
In particular, $P_t(G,\nu_G)$ does not depend on $\nu_G$.
\end{thm}

Note that in both Theorem \ref{thm:ABrecursion} and Theorem \ref{thm:LRgeneralG}, we
may replace $G$ by the Levi component $M_P$ of a parabolic subgroup $P$.

To invert the recursion relation \eqref{eqn:ab}, we
apply Theorem \ref{thm:inversionII}, with
$$
a_0(P)=P_t(M_P,\nu_P),
~b_0(P,\nu_P)=P_t^{ss}(M_P,\nu_P),
~n_P=2\dim(N_P)(g-1),
~\ep_P^G =4\rho_P^G.
$$
We obtain
\begin{thm}\label{thm:LRformula}
For any $\nu_G\in \Lambda^G_{P_0}$, we have
\begin{eqnarray*}
\lefteqn{
P_t^{ss}(G,\nu_G)=}\\&& \sum_{P\in \cP}
(-1)^{{\rm dim}(\fa_P^G)}
\Bigl({(1+t)^{2g}\over 1-t^2}\Bigr)^{{\rm dim}(\fa_P)}
\Bigl(\prod_{i=1}^{{\rm dim}(\fa_0^P)}
{(1+t^{2d_i(M_P)-1})^{2g}\over (1-t^{2d_i(M_P)-2})
(1-t^{2d_i(M_P)})}\Bigr)\\
&&  \cdot t^{2{\rm dim}(N_P)(g-1)}
\Bigl(\prod_{\alpha\in\Delta^G_P}
{1\over 1-t^{4\langle\rho^G_P,\alpha^\vee\rangle}}\Bigl)
\cdot t^{4\sum_{\alpha\in\Delta^G_P}\langle\rho^G_P,\alpha^\vee\rangle
\langle\varpi_\tal (\nu_G)\rangle} \in \bQ(t).
\end{eqnarray*}
\end{thm}
This is exactly Theorem \ref{thm:Pt}.

\end{appendix}

\end{document}